%% file: script.tex
\documentclass[aos, noinfoline]{my-imsart}
 \usepackage[margin=1.3in]{geometry}

\RequirePackage[numbers]{natbib}
\usepackage{appendix}
\usepackage{amsmath,amsfonts,amsthm,amssymb,amscd}
\usepackage{etex,etoolbox}
\usepackage{environ}
\usepackage{mathrsfs}
\usepackage{bbm}
\usepackage{bm}
\usepackage{mathtools}
\usepackage{cancel}
\usepackage{listings}
\usepackage{booktabs}
\usepackage{hdiv}
\usepackage[usenames,dvipsnames]{xcolor}
\usepackage{hyperref}
\usepackage[inline]{enumitem}

\usepackage{pdfsync}

\newcommand\myshade{85}
\colorlet{mylinkcolor}{Orchid}
\colorlet{mycitecolor}{YellowOrange}
\colorlet{myurlcolor}{Aquamarine}

\hypersetup{
  linkcolor  = mylinkcolor!\myshade!black,
  citecolor  = mycitecolor!\myshade!black,
  urlcolor   = myurlcolor!\myshade!black,
  colorlinks = true,
}

\DeclareMathOperator*{\argmin}{arg\,min}

\makeatletter
\providecommand{\@fourthoffour}[4]{#4}
\def\fixstatement#1{%
  \AtEndEnvironment{#1}{%
    \xdef\pat@label{\expandafter\expandafter\expandafter
      \@fourthoffour\csname#1\endcsname\space\@currentlabel}}}

\globtoksblk\prooftoks{1000}
\newcounter{proofcount}

\NewEnviron{proofatend}[1]{%
  \edef\next{%
    \noexpand#1%
    \unexpanded\expandafter{\BODY}}%
  \global\toks\numexpr\prooftoks+\value{proofcount}\relax=\expandafter{\next}
  \stepcounter{proofcount}}

\def\proofs{%
  \count@=\z@
  \loop
    \the\toks\numexpr\prooftoks+\count@\relax
    \ifnum\count@<\value{proofcount}%
    \advance\count@\@ne
  \repeat}
\makeatother

\newcommand{\myappendix}[1]{%
  \begin{proofatend}
    #1
  \end{proofatend}
}%
\usepackage{color}

\newif\ifgray
\graytrue
\grayfalse
\IfFileExists{BackgroundColor.tex}{\input{BackgroundColor}}{}
\IfFileExists{../BackgroundColor.tex}{\input{../BackgroundColor}}{}
\ifgray
\pagecolor{black!30!white}
\fi

\newcommand{\new}[1]{{\color{black} #1}}

\begin{document}
%%%%%%%%%%%%%%%%%%%%%
% Front-matter
\newcommand{\TitleMain}{Inference for high-dimensional \\ instrumental variables regression}
\newcommand{\RunningTitle}{High-dimensional IV regression}
\newcommand{\TitleAdd}{}
\newcommand{\front}{%
  \begin{frontmatter}
  \title{\TitleAdd \TitleMain}
  \runtitle{\RunningTitle}

  \begin{aug}

  \author{\fnms{David} \snm{Gold}\thanksref{t2,m1}\ead[label=e1]{}},
  \author{\fnms{Johannes} \snm{Lederer}\thanksref{t2,m1}\ead[label=e2]{ledererj@uw.edu}},\and
  \author{\fnms{Jing} \snm{Tao}\thanksref{t2,m1}\ead[label=e3]{jingtao@uw.edu}}
  \thankstext{t2}{First version: July, 2017. This version: November, 2019. The authors acknowledge support through Cloud Credits for Research by Amazon and through the University of Washington Royalty Research Fund.}
%  \ead[label=u1,url]{http://www.johanneslederer.com}
  \address[a]{University of Washington \\
  Department of Statistics \\
  Box 354322 \\
  Seattle, WA 98195-4322 \\
  Homepage: \href{https://johanneslederer.com}{johanneslederer.com}\\
 \printead{e2} \\
 \phantom{E-mail:\ } \printead*{e1}}
  \address[b]{University of Washington \\
  Department of Economics \\
  305 Savery Hall \\
  Box 353330 \\
  Seattle, WA 98195-3330 \\
 \printead{e3}}
  \affiliation{University of Washington\thanksmark{m1}}

  \end{aug}

  \begin{abstract}~
This paper concerns statistical inference for the components of a high-dimensional regression parameter despite possible endogeneity of each regressor.
Given a first-stage linear model for the endogenous regressors and a second-stage linear model for the dependent variable, we develop a novel adaptation of the parametric one-step update to a generic second-stage estimator.
We provide conditions under which the scaled update is asymptotically normal.
We then introduce a two-stage Lasso procedure and show that the second-stage Lasso estimator satisfies the aforementioned conditions.
Using these results, we construct asymptotically valid confidence intervals for the components of the second-stage regression coefficients. We complement our asymptotic theory with simulation studies, which demonstrate the performance of our method in finite samples.
  \end{abstract}

  \begin{keyword}
  \kwd{High-dimensional inference}\kwd{instrumental variables}\kwd{de-biasing}%this is the correct capitalization
  \end{keyword}

  \end{frontmatter}
}%
\front

%%%%%%%%%%%%%%%%%%%%%%%%%%%%%%%%%%%%%%%%%%%%%%%%%%%%%%%%%%%%%%%%%%%%%%%%%%%%%%%
% Body
\input{introduction/introduction}
\input{introduction/contribution}
\input{introduction/related}
\input{introduction/organization}
\input{introduction/notation}
\input{two-stage/two-stage}
\input{inference/inference}
\input{inference/endogeneity}
\input{inference/inverse}
\input{inference/normality}
\input{estimation/estimation}
\input{estimation/two-stage}
\input{estimation/bounds/bounds}
\input{estimation/bounds/second}
\input{inference/remainders}
\input{experiments/experiments}
\input{experiments/design}
\input{experiments/results}
\input{conclusion/conclusion}
% \input{discussion/discussion}

%\appendix
\section*{Acknowledgements}
We are grateful to the Editors and Reviewers for their valuable suggestions. We sincerely thank Alex Belloni, Mehmet Caner, Denis Chetverikov, Yanqin Fan, Amit Gandhi, Mohamed Hebiri, Joseph Salmon, Alexandre Tsybakov, and Jon Wellner for their insightful comments and Donghui Mai for her help in setting up the cloud computing.

%%%%%%%%%%%%%%%%%%%%%
% End

%%%%%%%%%%%%%%%%%%%%%%%%%%%%%%%%%%%%%%%%%%%%%%%%%%%%%%%%%%%%%%%%%%%%%%%%%%%%%%%
% Appendix
% \appendix
% \input{proofs/estimation}

%%%%%%%%%%%%%%%%%%%%%%%%%%%%%%%%%%%%%%%%%%%%%%%%%%%%%%%%%%%%%%%%%%%%%%%%%%%%%%%
% References
% \bibliographystyle{plain}

% Supplementary Material
%\begin{supplement}[id=suppA]
%  \sname{Appendix}
%  \stitle{Proofs and requisite materials}
%  \slink[doi]{COMPLETED BY THE TYPESETTER}
% \sdatatype{.pdf}
%  \sdescription{This supplement contains all proofs and requisite materials for the main text.}
%\end{supplement}
\bibliographystyle{imsart-number}
\bibliography{refs/refs}

\newpage
\begin{appendices}
\begin{center}
\textbf{\large APPENDIX: \\
PROOFS AND REQUISITE MATERIALS}
\end{center}
%technical lemmas: proofs/proofs
\input{proofs/proofs}
\end{appendices}

\end{document}

%% file: introduction/introduction.tex
\section{Introduction}
\subsection{Overview}
High-dimensional estimation has been extensively studied and is now ubiquitous in the data-intensive sciences \cite{BvdG11,G14,HTW15}.
High-dimensional inference, on the other hand, is much less developed.
In particular, although considerable progress has been made for inference in standard high-dimensional regression \cite{JM14,Ning17,vdGetal14,ZZ14}, much less is known for more complex models.

In this paper, we extend the study of high-dimensional inference to the \emph{linear instrumental variables (IV) model}.
To motivate the linear IV model, we consider the ordinary linear model
\begin{equation*}
\dy\equals\dX\rb + \du\,,
\end{equation*}
where~$\dy$ is the vector of responses,~$\dX$ is the matrix of regressors,~$\rb$ is the regression vector, and~$\du$ is a vector of random disturbances.
Standard inference for~$\rb$ using ordinary least squares is valid only if~$\E[\du|\dX]=0$.
However, this assumption is easily violated in practice in view of selection biases, omitted variables, measurement errors, and many other challenges common to data collection.
Hence, it is often more reasonable to allow for~$\E[\du|\dX]\neq 0$ and instead assume that~$\dX$ can be modeled based on observable variables~$\dZ$ that satisfy~$\E[\du|\dZ]=0$.
As standard in the econometric literature, we then call the regressors in~$\dX$ \emph{endogenous}, because they can be correlated with~$\du$, and we call the instrumental variables~$\dZ$ \emph{exogenous}, because~$\E[\du|\dZ]=0$~\cite{stock2003}.

Inference for such models in low-dimensional settings, where the number of samples is much larger than both the number of regressors in~$\dX$ and the number of regressors in~$\dZ$, has been extensively studied and put to use in economic applications and beyond~\cite{AP09}.
In the era of Big Data, however, low-dimensional settings often do not apply, because one wants to allow for flexible parameter combinations, or because many variables are measured in the first place.
We are thus interested in inference for double high-dimensional settings with the number of samples dominated by both the number of regressors in~$\dX$ and the number of regressors in~$\dZ$.

The method we propose is based on a novel adaptation of the parametric one-step update to a generic two-stage estimation procedure.
In parametric models, the one-step update~$\rbdb$ to an initial estimator~$\rbh$ is one Newton-Raphson step in the direction of a solution to the empirical analogue of the score equations.
This approach is similar to those of \cite{JM14, vdGetal14}, who de-bias the Lasso, and \cite{ZZ14}, who use a low-dimensional projection technique, to obtain asymptotic pivots for the low-dimensional components of the high-dimensional linear regression models when endogeneity is absent.
The present work extends this approach to two-stage estimation when both stages are high-dimensional and the regressors of interest may each be endogenous.
To adapt the one-step update to handle endogeneity of~$\dX$, we (i) choose the update as a step towards the solution to the empirical analogue of a valid moment condition and (ii) apply the update to a generic second-stage estimator~$\rbh$ that depends on the predicted conditional means~$\hat\E[\dX|\dZ]$.
The resultant estimator decomposes into a main term and four remainder terms, which contrasts with the single remainder term in the case of the de-biased Lasso in the ordinary linear model.
% The de-biased estimator decomposes into a main term linear in the noise and a remainder term that, under certain growth conditions on the model parameters, is asymptotically negligible.

We present high-level conditions under which the updated estimator yields asymptotic pivots for the components of~$\rb$, and we show, as an example, how these conditions may be satisfied by a two-stage Lasso estimation routine.
We assume a sub-Gaussian regime throughout for the noise elements and instrumental variables in order to support flexibility of distributional assumptions.
The main challenges of establishing the example are due to (i)~the involved structure of the remainder terms, whose control requires a variety of concentration bounds and lead to extensive proofs and (ii)~the estimation of the population precision matrix of the conditional means~$\E[\dX|\dZ]$, since these are not observed directly.

%% file: introduction/contribution.tex
\subsection{Our contributions}\label{contribution}

Our primary contribution is to develop a method with which to conduct statistical inference for the components~$\rbj$, $j=1,\dots,\px$ of a high-dimensional regression vector~$\rb$ despite endogeneity of the respective regressors.
We develop a novel adaptation of the one-step update and high-level conditions under which the updated estimator yields asymptotically Gaussian pivots for each~$\rbj$.
Our rigorous demonstration of conditions under which such inference is possible in the doubly high-dimensional setting differentiates the present paper from similar works such as \cite{BCH11}, who develop inferential methods for IV models with high-dimensional instruments and low-dimensional endogenous regressors, and \cite{Zhu15}, who work under a doubly high-dimensional regime but focus primarily on bounding the error of a two-stage Lasso estimator such as that of Section~\ref{estimation}.

A related contribution concerns sparse inverse covariance matrix estimation.
The updated estimator~$\rbdb$ depends on an estimate of the inverse covariance matrix~$\rT$ of the conditional means~$\E[\dx_i|\dz_i]$.
However, we do not observe these conditional means directly, and must base our estimate of~$\rT$ on the predictions~$\hat\E[\dx_i|\dz_i]$.
For this, we use a modifcation of the the CLIME estimator~$\rTh$ of~\cite{CLL11}.
Our paper is the first one to use such an estimator in the context of instrumental variable selection,
and we do novel work to account for the prediction step in deriving probabilistic guarantees for the estimator's performance.

Another contribution is to show that the updated second-stage Lasso estimator studied in Section~\ref{estimation} satisfies the high-level conditions in Section \ref{inference} and therefore supports inference for the~$\rbdbj$.
To show as much, we develop probabilistic bounds for the second-stage~$\ell_1$ estimation error, and we use these bounds to show asymptotic negligibility of the four remainder terms described in the previous section.
We also demonstrate the feasibility of the compatibility condition in the second-stage regression, thereby justifying the practical use of the second-stage rates.

A majority of  the proofs factor nicely into deterministic and stochastic components.
This also allows future analysts easily to combine the generic bounds contained in Section~\ref{proofs:estimation} of the Appendix with concentration results for specific error and design matrix distribution regimes and thereby derive the growth conditions required for good asymptotic behavior of the updated second-stage estimator under a variety of models.

%% file: introduction/related.tex
\subsection{Related work}\label{related}
Our work is related to the recent research on inference for high-dimensional linear instrumental variables models such as~\cite{BCCH12, belloni2018, FL14,GT14, Cheng&Liao2015, neykov15}.
\cite{BCCH12} use the Lasso to obtain representations of the optimal IVs of \cite{A74, A77, H82} for models in which the conditional mean of the response is linear in a small and fixed number of endogenous variables and show that the second-stage estimator is $\sqrt{n}$-consistent.
Following the line of the seminal work by \cite{BCCH12}, some recent papers propose different novel procedures to select many instruments when the number of second-stage regressors remains to be fixed or low-dimensional (for instance, \cite{Hansen&Kozbur2014RJIVE, Caner&Fan2015hybrid, Fan&Zhong2018}, among others).
\cite{Cheng&Liao2015} propose a Lasso procedure to select valid and relevant moments for the GMM estimation when the number of moments increases with sample size.
However, and in contrast with the present paper, the dimensions of moments and parameters of interest are both smaller than the sample size.
Compared to \cite{BCCH12} and \cite{Cheng&Liao2015}, we allow both the number of IVs and the number of endogenous regressors to be bigger than the sample size.
\cite{GT11,GT18} construct robust confidence sets based on their Self-Tuned Instrumental Variables (STIV) estimator and confidence bands after bias correction (requiring a type of strong instruments condition).

Much of the present work is devoted to solving the inference problem for parameters of interest in a high-dimensional linear IV model with homoscedasticity by accounting for the prediction error when the first- and second- stage regression models are both high-dimensional.
This contrasts with the methods of \cite{FL14}, who do not account for the need to predict the optimal instruments \cite{A74, A77, IDN03, N90a}.
To our knowledge, such analysis under an~$\ell_1$-regularized estimation procedure is new in the literature.

Recent work by \cite{lin2015regularization, Zhu15} also propose estimation methods for linear IV models when the regressors of both stages are high-dimensional but do not rigorously develop asymptotic methods for inference.
We notice that \cite{neykov15} also provide an inferential method for high-dimensional linear IV models by using a Dantzig selector.
However, their method requires that the number of instruments equals the number of endogenous regressors and they do not account for the first-stage of estimation.
Another related work is concurrently developed by \cite{belloni2018}, who propose regularized estimation of nuisance parameters that appear in carefully constructed empirical orthogonality conditions.
Our one-step update approach can be interpreted as an iteration of the Newton-Raphson method, while \cite{belloni2018} build on the idea of Neyman orthogonality.
Even though the two papers take different approaches, both show that $\sqrt{n}$-consistent estimators for low-dimensional parameters can be constructed in high-dimensional IV models.

\new{In a very recent work \cite{Caner&Kock2018} introduce a de-sparsified, $\ell_1$-penalized two-stage GMM estimator.
They develop estimation error bounds and inferential procedures based on this estimator in the doubly high-dimensional setting and extend our work in that they
(i) allow conditional heteroskedasticity in the second-stage noise elements,
(ii) specify random components in terms of moment conditions rather than sub-Gaussianity, and
(iii) do not require $\ell_0$-sparsity of the first-stage regression vectors.
Also noteworthy is that the authors do not predict the conditional means~$ZA$ for the weighting scheme of the second-stage.
This leads to expressions for the scale factor of the asymptotic pivots in their Theorem 2.(i) that differ those of our Theorem~\ref{wl}.
The scale factors of the latter such pivots is identical to the asymptotic variance of optimal estimators in low-dimensional IV models with homoscedastic structural errors~\cite{chamberlain1987asymptotic, newey1990efficient}.
In \cite{CanerKock18}, the same authors develop asymptotic theory to support both estimation and inference for the conservative Lasso under a wide variety (including heteroscedasticity and non-sub-Gaussianity) of error regimes, thereby laying foundations for future work considering two-stage estimation.

%For high-dimensional inference, we conjecture that our method achieves asymptotic efficiency in the sense of \cite{Jankova16} Theorem 3.
}

%% file: introduction/organization.tex
\subsection{Organization}\label{organization}
The rest of the essay is organized as follows.
We introduce our model and a generic two-stage estimation procedure in Section \ref{s:two-stage}.
In Section \ref{inference}, we propose the one-step update inference procedure and demonstrate conditions under which the update yields asymptotically Gaussian pivots.
In Section \ref{estimation}, we introduce a two-stage Lasso estimator of the parameter $\rb$ and show that it is suitable for use with the inference procedure developed in Section \ref{inference}.
Finally, in Section~\ref{s:numerics}, we present the results of numerical studies that demonstrate the relevance of our theoretical results to finite samples.
All proofs are contained in the Appendix.

%% file: introduction/notation.tex
\subsection{Basic notation and preliminaries}

We adopt the following general notational conventions.
For~$p\in\N$, we let $[p] := \{1,\ldots,p\}$.
We typically use bold and non-bold lowercase letters denote vectors and scalars, respectively.
We use bold uppercase letters to denote matrices.
We typically denote the components of a vector (matrix) by the non-bold (lowercase) counterpart of the letter that denotes the vector (matrix).
If~$\aM \in\R^{n\times p}$, with components~$\amc_{ij}$, we use a superscript to refer to columns~$\am^j = (\amc_{1j},\ldots,\amc_{nj})^\top$ and a subscript to refer to rows~$\am_i = (\amc_{i1},\ldots,\amc_{ip})$.
We let~$\|\cdot\|_q$ and~$\ip{\cdot,\cdot}$ denote the usual~$\ell_q$ norm and inner product over Euclidean spaces, respectively.

For~$\am \in \R^p$, we let~$\supp(\am) := \{j\in[p]: \amc_j \neq 0\}$, $\|\am\|_0 = |\supp(\am)|$, and $\|\am\|_\infty = \max_{j\in[p]} \{\amc_j\}$.
For matrices~$\aM\in\R^{\n\times p}$, we let~$\|\aM\|_\infty=\max_{i,j\in[\n]\times[p]}|\amc_{ij}|$,~$\|\aM\|_{L_1}=\max_{j\in[p]}\|\am^j\|_1$, and~$\|\aM\|_{L_2}=\max_{j\in[p]}\|\am^j\|_2$.
For matrices~$\aM_1$ and $\aM_2\in\R^{\n\times p}$, we write~$\aM_1 \succ \aM_2$ if~$\aM_1 - \aM_2$ is positive-definite.
% We let~$\otimes$ denote the Kronecker product and write~$\am^{\kp2} = \am \kp \am$ for~$\am\in\R^p$;

For quanities~$x$ indexed by~$i\in[n]$, we let~$\En[x_i] = n^{-1}\sum_{i=1}^n x_i$.
If~$X_n$ is a sequence of random variables, we write~$X_n\leadsto X$ if~$X_n$ converges weakly to~$X$.
For~$a,b\in\R$, we let~$a\vee b = \max(a,b)$ and~$a\wedge b = \min(a,b)$.
We write~$a_n \lesssim b_n$ if~$a_n \leq C_nb_n$ for a~$C_n$ that is of constant order.
% we write~$a_n \lesssim_\Prb b_n$ if~$a_n = O_\Prb(b_n)$.
We say that a sequence of events~$\event \equiv \event_n$ occurs with probability approaching one if $\lim_{n\to\infty} \Prb \, \event_n = 1$.

% Many specific scalar quantities we introduce are indexed by a vector or matrix subscript.
% Typically, such scalar quantities refer to a maximum over the components of the subscripting array quantity.
% For instance:
% (i) if the quantity~$\sd$ is a function of the components~$\dvc_j$ of a vector~$\dv$, then we may use~$\maxsdv$ to denote the maximum value of~$\sd$ over such indices;
% (ii) if the quantity~$M$ is a function of the columns~$\ra^j$ of the matrix~$\rA$, then we may use~$\Mra$ to denote the maximum value of~$M$ over such columns.
We recall the following definitions for sub-Gaussian and sub-exponential norms.

\begin{defi}[Sub-Gaussian and sub-exponential norms]\label{def:orlicz}
For~$q\ge1$ and a random variable~$X$, we write
\benn
  \|X\|_{\psi_q} \smallpad{:=} \inf\set{t\in(0,\infty)\,:\, \E[\exp(|X|^q/t^q)-1]\le1}\,.
\eenn
if the infimum exists.
The \emph{sub-Gaussian norm} of a random variable~$X$ is given by~$\|X\|_{\sgn}$;
the \emph{sub-exponential norm} of a random variable~$X$ is given by~$\|X\|_{\sen}$.
The corresponding norms for a random $p$-vector~$\bs{X}$ are given by
\benn
  \|\bs{X}\|_{\psi_q} \smallpad{:=}
  % \sup_{\substack{\bs{x}\in\R^p\\\|\bs{x}\|_2=1}}
  \sup_{\bs{x}\in\R^p\,:\,\|\bs{x}\|_2=1}
  \|\ip{\bs{X},\bs{x}}\|_{\psi_q} \,.
\eenn
\end{defi}

%% file: two-stage/two-stage.tex
\section{Two-stage estimation}\label{s:two-stage}

To contend with endogeneity, the method of instrumental variables isolates variation in the endogenous regressors induced by the instrumental variables.
In Section~\ref{model}, we posit the two-stage linear IV model to describe this relationship.
In Section~\ref{generic}, we discuss a generic two-stage estimation routine that respects the structure of the model.

\subsection{Model}\label{model}

Our model of interest is
\begin{align}
  \dyi \aequals \dx_i^\top\rb + \dui \,,\label{def:yi}\\
  \dxij \aequals \dz_i^\top\ra^j + \dvij \,,\label{def:xij}
\end{align}
where:~$i$ ranges from~1 to~$\n$ (unless stated otherwise);~$j$ ranges from 1 to $\px$ (unless stated otherwise);
the vectors~$\dx_i \in \R^{\px}$ consist of the~\emph{second-stage regressors}~$\dxc_{i1},\ldots,\dxc_{i{\px}}$;
the vector~$\rb\in\R^{\px}$~is the parameter of interest;
the vectors~$\dz_i\in\R^{\pz}$ consist of the~\emph{first-stage regressors}~$\dzc_{i1},\ldots,\dzc_{i{\pz}}$;
% , which satisfy~$\E[\dz_i]=\0$;
the quantities~$\dui$ and~$\dv_i := (\dvc_{i1},\ldots,\dvc_{i\px})^\top$ are random noise elements that satisfy
\be\label{orthogonality}
  \E[\dui|\dz_i] = 0\,,\qquad \E[\dv_i|\dz_i]=\0 \,,
\ee
and the vectors~$\ra^j$ are regression parameters up to which the respective conditional means~$\ddij := \E[\dxij|\dz_i] = \dz_i^\top\ra^j$ are specified.
We call the models of~\eqref{def:xij} and~\eqref{def:yi} the first-stage and second-stage models, respectively.
Note that the setup in (\ref{orthogonality}) is similar to the one in \cite{BCCH12} except for this significant difference: we consider both stages to be high-dimensional while the second stage parameters in \cite{BCCH12} are low-dimensional.
To focus on presenting the main idea and inference steps, we ignore approximation errors in the model but one can add approximation errors in both \eqref{def:xij} and \eqref{def:yi} at the expense of more tedious derivations and notation (see the discussion of Theorem 6.3 in \cite{BvdG11}).

In matrix notation, we write
\benn
  \dy \equals \dX\rb + \du
\eenn
and
\benn
  \dX \equals \dD + \dV \equals \dZ\rA + \dV \,,
\eenn
where the vectors~$\dy,\du\in\R^{\n}$ consist of the responses $\dyi$ and the noise components~$\dui$, respectively;
the matrix~$\dX\in\R^{\n\times\px}$ has columns $\dx^j = (\dxc_{1j},\ldots,\dxc_{\n j})^\top$ and rows~$\dx_i=(\dxc_{i1},\ldots,\dxc_{i\px})$;
the matrix~$\dD=\E[\dX|\dZ]\in\R^{\n\times\px}$ has columns~$\dd^j=(\ddc_{1j},\ldots,\ddc_{\n j})^\top$ and rows~$\dd_i=(\ddc_{i1},\ldots,\ddc_{i\px})$;
the matrix~$\dZ\in\R^{\n\times\pz}$ has columns~$\dz^k=(\dzc_{1k},\ldots,\dzc_{\n k})^\top$ and rows~$\dz_i=(\dzc_{i1},\ldots,\dzc_{i\pz})$;
and the matrix~$\rA \in \R^{\pz\times\px}$ has columns given by $\ra^j$.
We make the following assumption concerning the~$\n$-indexed sequence of regression parameters~$\rA,\rb$.

\begin{assu}[Regularity of~$\rA,\rb$]\label{regreg}
The quantities~$\|\rA\|_{L_1}$ and~$\|\rb\|_1$ are bounded above by universal constants~$\Mra,\Mb<\infty$, respectively.
\end{assu}

\noindent We let~$\rSzh=\dZ^\top\dZ/\n$ denote the empirical Gram matrix of the instrumental variables.

As remarked earlier, the linear IV model has been studied extensively in the low-dimensional setting, where the number $\px$ of endogenous variables $\dx^j$ is fixed.
We are particularly concerned with the high-dimensional regime in which both $\px$ and the number $\pz$ of instrumental variables $\dz^k$ increase with $\n$.
Our results generalize to the low-dimensional case in which $\pz$ and $\px$ are held fixed with respect to $\n$, but we do not treat this case explicitly in the present essay.
Regardless of whether the model is high-dimensional, we require that $\px \le \pz$ in order to maintain identifiability of~$\E[\dx_i|\dz_i]$.
% We note that the results of Sections~\ref{inference} and~\ref{estimation} continue to hold if the~$\dx^j$ are exogenous, though this setting is not our focus.

We study a sub-Gaussian regime for the noise components as well as the instrumental variables, which we treat as random throughout and for which we give marginal results.
This regime encompasses the typical Gaussian model considered in the high-dimensional literature and allows for flexibility in modeling assumptions.

\begin{assu}[Specification of~$\dz_i$]\label{assu:dzi}
The instrumental variables~$\dz_i$ are i.i.d. and sub-Gaussian with sub-Gaussian norm~$\sgnz:=\|\dz_i\|_\sgn$ and satisfy~$\E[\dz_i]=0$ for each~$i\in[\n]$.
Considered as components of an~$\n$-indexed sequence of models, the quanities~$\sgnz$ are bounded away from zero and infinity.
\end{assu}

\begin{rek}[Specification of~$\dz_i$]
We require that the first-stage regressors~$\dz_i$ have mean zero in order to simplify the following exposition and to apply concentration results under more specific distributional assumptions, such as in Lemma~\ref{feassgg}.
This assumption can be relaxed at the expense of brevity and given a sufficient reformulation of the required concentration results.
Similarly, the condition that~$\sgnz=O(1)$ can be relaxed at the expense of introducing more complex growth conditions in later results.
\end{rek}

\begin{assu}[Specification of~$\dv^j$ and~$\du$]\label{assu:subgaussian}
The noise vectors~$\dv^j$ and~$\du$ are sub-Gaussian with sub-Gaussian norms~$\sgnvj:=\|\dv^j\|_{\sgn}$ and~$\sgnu:=\|\du\|_\sgn$.
Considered as components of an~$\n$-indexed sequence of models, the quantities~$\sgnvj$ and~$\sgnu$ are bounded strictly away from zero and infinity.
% The
% \benn
% 	\rSuv \smallpad{:=}
% 	\begin{pmatrix}
% 		\sdu^2 \ \sduv^\top \\
% 		\sduv \ \rSv
% 	\end{pmatrix} \,,
% \eenn
% where~$\sduv:=(\sduvj{1},\ldots,\sduvj{\px})^\top$ consists of the noise covariances~$\sduvj{j}:=\cov(\du,\dv^j)$, and where~$\rSv$ is an unstructured covariance matrix with diagonal entries~$\sdvk{j}^2:=\var(\dv^j)$ for~$j\in[\px]$.
% Further, considered as components of an $\n$-indexed sequence of models, the variances~$\sdu^2\equiv\sdun^2$ and~$\sdvk{j}^2\equiv\sdvkn{j}^2$ are bounded strictly away from zero and infinity for~$j\in[\px]$.
\end{assu}

\noindent
Note that Assumption~\ref{assu:subgaussian} makes no stipulations concerning the joint covariance structure of the~$\dui$ and~$\dv_i$.
The assumption therefore allows for nontrivial covariance between the two stages of noise, which can be used to model endogeneity of the~$\dx_i$.
Furthermore, Assumption~\ref{assu:subgaussian} allows for heteroscedasticity amongst the components of the $\dv^j$.
We require homoscedasticity of the second-stage noise elements~$\dui$ for Theorem~\ref{wl} and Lemma~\ref{wjsehcons}; all other results hold in the presence of heteroscedasticity.

% \begin{rek}[Heterscedastic, non-Gaussian noise]
% \david{What should we do with this remark? I think we could omit it, or reduce it to one sentence.} Heteroscedasticity and non-Gaussianity do not in principle pose problems to the approach we develop in Sections~\ref{estimation} and~\ref{inference}.
% Our results are modular with respect to error regimes: accounting for heteroscedastic, non-Gaussian errors is a matter of combining the appropriate concentration results with the generic bounds we present in Appendix~\ref{proofs}.
% We will pursue this task in a later version of the present essay.
% \end{rek}

\subsection{Generic two-stage estimators}\label{generic}

We formulate our proposed method of inference for the components~$\rbj$ of the second-stage regression parameter~$\rb$ in terms of generic estimators that reflect the structure of the model described above.
We now introduce notation that will be used in Section~\ref{inference}.

For each~$j\in[\px]$, let~$\rah^j \equiv \rah^j(\dx^j, \dZ)$ denote a generic \emph{first-stage estimator} of the first-stage regression vector~$\ra^j$ based on the data~$\dx^j$ and~$\dZ$.
We write~$\rAh := (\rah^1,\ldots,\rah^{\px})$ for the matrix of estimated regression vectors.
From such an estimator~$\rAh$ we may predict the conditional means~$\dd_i = \E[\dx_i|\dz_i]$ for $i\in[\n]$ with~$\ddh_i \smallpad{:=} \dz_i^\top\rAh$;
we write $\dDh$ for the predicted conditional mean matrix whose rows are given by the $\ddh_i$, and we write~$\rSdh:=\dDh^\top\dDh/\n$ and~$\rSd := \E[\rSdh]$.
Our choice of the notation $\dDh$ reflects the fact that this quantity predicts and, under suitable conditions, approaches in probability the conditional mean matrix $\dD$.
We write $\rbh\equiv\rbh(\dy,\dDh)$ for a generic \emph{second-stage estimator} of the second-stage regression parameter $\rb$ based on the response $\dy$ and the predicted conditional means $\dDh$.

%% file: inference/inference.tex
\section{Main proposal}\label{inference}
\myappendix{
\section{Materials required for Section~\ref{inference}}\label{proofs:inference}
}

Our main contribution is to develop a method for statistical inference for the components~$\rbj$ of the second-stage regression vector~$\rb$.
In general, statistical inference for high-dimensional regression parameters is a difficult problem.
Regularized estimators, such as the Lasso and ridge regression, are often used for the purpose of high-dimensional parameter estimation but generally do not have asymptotic distributions suitable for inference \cite{Knight00, Potscher09}.
In studying the model of Section~\ref{model}, we must also account for the dependence of the second-stage estimator on the first-stage estimators. 

The basis for our procedure is to adapt the parametric one-step update to the two-stage estimation procedure described in Section~\ref{generic}.
We first briefly review the use of the one-step estimator in parametric models and its application to high-dimensional inference for the ordinary linear model.
Then we adapt the one-step update to the two-stage estimation procedure described in Section~\ref{generic}.
% We demonstrate that the updated second-stage estimator decomposes into a main term that is linear in the second-stage noise $\du$ and four additional remainder terms.
Section~\ref{normality} discusses high-level conditions under which the scaled updated estimator is asymptotically normal. 

%% file: inference/endogeneity.tex
\subsection{One-step with endogeneity}\label{endogeneity}
\myappendix{
\subsection{Materials required for Section~\ref{endogeneity}}\label{proofs:endogeneity}
}
In this section, we develop a novel adaptation of the one-step update that, under suitable high-level conditions, yields asymptotic pivots for the second-stage components~$\rbj$ of the two-stage model described in Section~\ref{model}.
We note that the present development is valid for any initial second-stage estimator~$\rbh$.
To demonstrate that the high-level conditions are satisfied requires consideration of particular estimators.

In summary, the one-step update is a general method for constructing efficient estimators for parametric and semiparametric models \cite[Sections 2.5, 7.3]{BKRW98}. Recall that the Newton-Raphson method for finding the root in~$\argb$ to a \emph{target system} of~$\px$ equations
\benn
	\score(y_i,\dx_i;\argb) \smallpad\equiv (\scorec_1(y_i,\dx_i;\argb),\ldots,\scorec_{\px}(y_i,\dx_i;\argb))^\top \equals \0
\eenn
is to update an approximation~$\argb^k$ by the rule
\benn
	\argb^{k+1} \equals \argb^k -
	\left[\left.\frac{\bs\partial\score}{\bs\partial\argb}\right|_{\argb=\argb^k}\right]^{-1}\score(y_i,\dx_i;\argb^k)\,,
\eenn
where~$\left.\frac{\bs\partial\score}{\bs\partial\argb}\right|_{\argb=\argb^k}$ is the Jacobian matrix of~$\score$ with respect to~$\argb$ evaluated at~$\argb^k$.
In the ordinary least squares regression model
\be\label{ordinary}
	\dy \equals \dX\rb + \du \,,
\ee
the score function $\score(\dyi, \dx_i; \argb) = -\dx_i(y_i-\dx_i^\top\argb)$
% \benn
% 	\score(\dyi, \dx_i; \argb) \equals -\dx_i(y_i-\dx_i^\top\argb)
% \eenn
satisfies~$\E[\score(\dyi, \dx_i; \rb)] = \0$ given the \emph{orthogonality condition}
\be\label{orthogonality1stage}
	\E[\dx_i\dui] \equals \0\,.
\ee
The \emph{one-step update}~$\rbdb$ to an initial estimate~$\rbh$ of~$\rb$ is given by
\benn
	\rbdb \equals \rbh + \rTh \dX^\top(\dy-\dX\rbh)/n\,,
\eenn
where~$\rTh$ denotes the inverse of
\benn
	\left.\frac{\bs\partial(-\dX^\top(\dy-\dX\argb)/\n)}{\bs\partial\argb}\right|_{\argb=\rbh} \equals \dX^\top\dX/n \equals \rSxh \,.
\eenn
The case when the model is high-dimensional is less well studied. When~$\px>\n$, the empirical covariance matrix~$\rSxh$ is not invertible,
and we have instead
\benn
	\rbdb
	% \equals \rbh + \rTh\dX^\top(\dy-\dX\rbh)/\n
	\equals
	\rb + \rTh\dX^\top\du/\n + \underbrace{(\rTh\rSxh - \Id)(\rb-\rbh)}_{\rem/\sqrt{\n}}\,,
\eenn
where~$\rTh$ denotes an approximate inverse of the Jacobian matrix.
The latter term~$\rem/\sqrt{\n}$ in the above display is the ``remainder'' after incomplete inversion of $\rSxh$.
Thus, in the high-dimensional one-stage linear model, the one-step update satisfies
\be\label{eq:upshot}
	\sqrt{n}(\rbdbj-\rbj)
	\equals
	\frac{1}{\sqrt{\n}}\sum_{i=1}^{\n} \rth_j^\top\dx_i\dui + \remc_j \,,
\ee
where~$\rth_j$ is the~$j^\text{th}$ row of~$\rTh$.
The structure of the main term on the right-hand side above suggests to use~$\sqrt{n}(\rbdbj-\rbj)/\wjseh$, where~$\wjseh$ is an appropriate estimate of~$\wjse = (\E[\ip{\rth_j,\dx_i}^2\dui^2])^{1/2}$, as an asymptotic pivot for~$\rbj$.

When the initial estimator~$\rbh$ is the Lasso, the updated estimator~$\rbdb$ is sometimes called the desparsified~\cite{vdGetal14} or de-biased~\cite{JM14} Lasso, though these authors obtain the form of~$\rbdb$ by means other than the one-step update.
The general upshot of their results is that if~$\|\rem\|_\infty=o_{\Prb}(1)$, and if~$\rth_j$ and~$\dx_i$ are independent of~$\dui$, then the updated Lasso estimator yields asymptotically Gaussian pivots for the parameter components.

In contrast to the ordinary linear regression model, the challenge we face in the case of high-dimensional IV model is that the condition in~\eqref{orthogonality1stage} does not hold.
Instead, the conditional moment restriction~$\E[\dui|\dz_i]=0$ in~\eqref{orthogonality} entails the orthogonality condition $\E[\dd_i\dui] = \0$ for the conditional means~$\dd_i=\E[\dx_i|\dz_i]$.
This suggests that, to develop a one-step update for a generic second-stage estimator~$\rbh$ of~$\rb$ of the present model, we ought to take the empirical analogue
\benn
	\En[\scoret(y_i,\dx_i,\ddh_i; \argb)] \smallpad{:=}  \En[-\ddh_i(\dyi-\dx_i^\top\argb)] \equals -\dDh^\top(\dy-\dX\argb)/\n \equals \0 \,,
\eenn
of~$\E[\dd_i\dui] = \0$ as the target system for which the root is sought via a Newton-Raphson update.
We have elected to base the target system on the moment condition~$\E[\dd_i\dui] = \0$ in accordance with optimal weighting regimes for generalized method of moments~(GMM) estimators; see~\cite{A74,A77,H82,N90a}.
Further, since the~$\dd_i$ are generally unavailable, we instead use the predicted conditional mean matrix~$\dDh$ in the target system above.
The one-step update~$\rbdb$ to a second-stage estimator~$\rbh$ is then given by
\be\label{def:one-step}
	\rbdb \equals \rbh - \rTh\En[\scoret(y_i,\dx_i,\ddh_i; \argb)] \equals \rbh + \rTh\dDh^\top(\dy-\dX\rbh)/\n \,,
\ee
where we continue to let~$\rTh$ denote an (approximate) inverse to the Jacobian matrix in~$\argb$ of the score~$\scoret(y_i,\dx_i,\ddh_i; \argb)$.

If one were to follow strictly the prescription of the Newton-Raphson method for selection of~$\rTh$ for the updated second-stage estimator~$\rbdb$, one would select~$\rTh\approx[\dDh^\top\dX/\n]^{-1}$ to approximate the inverse of the Jacobian of~$\scoret$ evaluated at~$\rbh$.
% ~$\rTh\approx\left[\bs\partial\En[\scoret(\argb)]/\bs\partial\argb(\rbh)\right]^{-1}$
% \begin{align*}
% 	\rTh \smallpad\approx
% 	\left[\frac{\bs\partial\En[\scoret(\argb)]}{\bs\partial\argb}(\rbh)\right]^{-1}
% 	\aequals
% 	\left[\frac{\bs\partial(-\dDh^\top(\dy-\dX\argb)/\n)}{\bs\partial\argb}(\rbh)\right]^{-1} \\
% 	\aequals \big[\dDh^\top\dX/\n\big]^{-1} \,.
% \end{align*}
However, the decomposition obtained in the following lemma suggests that~$\rTh$ ought to control, say, the sup-norm of~$\rTh\rSdh-\Id$, and hence aim to invert~$\rSdh:=\dDh^\top\dDh/\n$ rather than~$\dDh^\top\dX/\n$.
We emphasize that the one-step formulation, insofar as it follows the Newton-Raphson method, is merely a vehicle for producing an updated estimator~$\rbdb$.
In particular, Lemma~\ref{lem:one-step} is valid regardless of what convergence properties an actual Newton-Raphson algorithm incorporating a specific choice of~$\rTh$ may exhibit.
We may choose~$\rTh$ in whatever manner is most appropriate for achieving our goal, which is to obtain a tractable limiting distribution for~$\sqrt{n}(\rbdbj-\rbj)/\wjse$, where~$\wjse$ is an appropriate scale factor.
That said, the two suggestions for how to choose~$\rTh$ may be reconciled somewhat by noting that
% \begin{align*}
% 	\dDh^\top\dX/\n \aequals \dD^\top\dX/\n + (\dDh-\dD)^\top\dX/\n \\
% 	\aequals \dD^\top\dD/\n + \dD^\top\dV/\n + (\dDh-\dD)^\top\dD/\n + (\dDh-\dD)^\top\dV/\n \,,
% \end{align*}
% and
% \begin{align*}
% 	\dDh^\top\dDh/\n \aequals \dD^\top\dD/\n + (\dDh-\dD)^\top\dD/\n \\
% 	&\quad + \dD^\top(\dDh-\dD)/\n + (\dDh-\dD)^\top(\dDh-\dD)/\n \,.
% \end{align*}
both~$\dDh^\top\dZ/\n$ and~$\dDh^\top\dDh/\n$ are equal to the empirical Gram matrix $\rSdh$ modulo additional terms whose sup-norms can be controlled given a rate for~$\|\rAh-\rA\|_{L_1}$ and appropriate concentration results for~$\|\dZ^\top\dv^j/\n\|_\infty$.
In turn, one finds~$\|\rSdh-\rSd\|_\infty = o_{\Prb}(1)$ under appropriate growth restrictions on~$\px$; see Lemma~\ref{concentration}.

For our purposes, we consider the matrix~$\rTh$ primarily as an estimator of the population quantity~$\rT:=\E[\dd_i\dd_i^\top]^{-1}$.
In particular, we require good behavior of~$\rTh$ as such an estimator to derive the asymptotic distribution of~$\sqrt{\n}(\rbdbj-\rbj)$.

The following lemma characterizes a similar decomposition of the updated estimator $\rbdb$ as in the one-stage model.

%%%%%%%%%%%%%%%%%%%%%%%%%%%%%%%%%%%%%%%%%%%%%%%%%%%%%%%%%%%%%%%%%%%%%%%%%%%%%%%
%%%%%%%%%%%%%%%%%%%%%%%%%%%%%%%%%%%%%%%%%%%%%%%%%%%%%%%%%%%%%%%%%%%%%%%%%%%%%%%

\begin{lem}[Decomposition of one-step second-stage estimator]\label{lem:one-step}
Consider the two-stage linear model described in Section~\ref{model}.
Let~$\dDh$ be a prediction of the conditional mean matrix~$\dD$ from an estimate $\rAh$ of the first-stage regression matrix~$\rA$.
Let~$\rbh$ be a second-stage estimator based on the predictions~$\dDh$.
Let~$\rTh$ denote an estimator of~$\rT=\E[\dd_i\dd_i^\top]^{-1}$.
The one-step second-stage estimator
\benn
	\rbdb \equals
  \rbla + \rTh\dDh^\top(\dy-\dX\rbh)/\n
\eenn
satisfies~$\sqrt{\n}(\rbdb-\rb) \equals \rT\dD^\top\du/\sqrt{\n} + \sum_{\ell=1}^4\rem_\ell$, where
% \benn
% 	\rbdb-\rb \equals \rTh\dDh^\top\du/n + \rTh\dDh^\top(\dX-\dDh)(\rb-\rbla)/n + (\rTh\rSdh - \Id)(\rb-\rbla)\,,
% \eenn
\benn
\begin{split}
	\rem_1 \aequals (\rTh-\rT)\dD^\top\du/\sqrt{\n} \,, \\
	\rem_3 \aequals \rTh\dDh^\top(\dX-\dDh)(\rb-\rbh)/\sqrt{\n} \,,
\end{split}\qquad
\begin{split}
	\rem_2 \aequals \rTh(\dDh-\dD)^\top\du/\sqrt{\n} \,, \\
	\rem_4 \aequals \sqrt{\n}(\rTh\rSdh-\Id)(\rb-\rbh) \,.
\end{split}
\eenn
\end{lem}

% \bproof{Proof of Lemma \ref{lem:one-step}}
% \begin{proofatend}{Proof of Lemma \ref{lem:one-step}}
\myappendix{%
\begin{proof}[Proof of Lemma \ref{lem:one-step}]
Note that
\begin{align*}
	\rbdb \aequals \rbh - \rTh\En[\scoret(\rbh)] \\
	\aequals \rbh + \rTh\dDh^\top(\dy-\dX\rbh)/\n \\
	\aequals \rbh + \rTh\dDh^\top(\dX[\rb-\rbh]+\du)/\n \\
	\aequals \rbh + \rTh\dDh^\top(\dDh[\rb-\rbh] + [\dX-\dDh][\rb-\rbh] + \du)/\n \\
	\aequals \rb + \rTh\dDh^\top\du/\n + \underbrace{\rTh\dDh^\top(\dX-\dDh)(\rb-\rbh)/\n}_{\rem_3/\sqrt{\n}} + \underbrace{(\rTh\rSdh-\Id)(\rb-\rbh)}_{\rem_4/\sqrt{\n}} \,.
\end{align*}
Now decompose the second term on the right-hand side above as follows
\begin{align*}
	\rTh\dDh^\top\du/\n
	\aequals
	\rTh\dD^\top\du/\n
		+ \rTh(\dDh-\dD)^\top\du/\n \\
	\aequals
	\rT\dD^\top\du/\n + \underbrace{(\rTh-\rT)\dD^\top\du/\n}_{\rem_{1}/\sqrt{\n}} +
		\underbrace{\rTh(\dDh-\dD)^\top\du/\n}_{\rem_{2}/\sqrt{\n}}
\end{align*}
to complete the proof.
\end{proof}
}%
\noindent\new{As in the case of the main term for the ordinary one-step update discussed in~(\ref{eq:upshot})}, this observation is similar to that of \cite{neykov15}, who derive a similar asymptotic linearization but do not account for prediction of the conditional means $\dd_i$.
Indeed, due to the two stages of estimation, the update incurs four remainder terms as opposed to the single term in~\cite{JM14,vdGetal14}.
In Section~\ref{normality}, we show that the quantity~$\sqrt{\n}(\rbdbj - \rbc_j)/\wjse$, where~$\wjse^2:=\E[\ip{\rt_j,\dd_i}^2\dui^2]$,
% \be\label{def:wjn}
% 	\Wjn \relates{:=}
%   \sqrt{\n}(\rbdbj - \rbc_j)/\wjseh \,,
% \ee
% where~$\wjseh^2$ is an appropriate estimator of
% \be
% 	\wjse^2 \smallpad{:=} \E\big[\En[\ip{\rt_j,\dd_i}^2\dui^2]\big]\,,
% \ee
converges weakly to a~$\Normal(0,1)$ random variable under high level conditions on the remainder terms~$\remc_{\ell,j}$ and that the limit continues to hold if~$\wjse$ is replaced by an appropriate estimator.
From this result one may construct asymptotically valid confidence intervals for the regression components~$\rbj$.

% \todo{Possibly remove the following exposition if we include an analysis of the heteroscedastic SE.} Note that the appropriate form of~$\wjseh$ will depend on the modeling assumptions.
% For instance, in homoscedastic case, the quantity~$\wjse$ reduces to
% % to~$\sdu^2\E[\ip{\rt_j,\dd_i}^2]$, where we recall~$\sdu^2=\E[\dui^2]$.
% \begin{align}\label{homoscedastic}
% 	\E\big[\En[\ip{\rt_j,\dd_i}^2\dui^2]\big]
% 	\aequals
% 	\fon\sum_{i=1}^{\n}\E[\ip{\rt_j,\dd_i}^2\dui^2] \nonumber\\
% 	\aequals
% 	\fon\sum_{i=1}^{\n}\E\big[\ip{\rt_j,\dd_i}^2\E[\dui^2|\dz_i]\big] \nonumber\\
% 	\aequals
% 		\sdu\rt_j^\top\rSd\rt_j \equals \sdu\rTjj \,,
% \end{align}
% which informs the structure of the corresponding estimator.

% We have presented a strategy for how Lemma~\ref{lem:one-step} may be used as the basis for methods of statistical inference for the second-stage regression components $\rbj$.
We have described a strategy for inference for the components~$\rbdbj$.
To implement the strategy for a specific choice of first- and second-stage estimators, one must identify the conditions under which the remainder terms $\rem_\ell$ vanish in probability.
We demonstrate such an implementation in Section~\ref{estimation}.
The conditions in turn depend on the properties of the estimator~$\rTh$.
In the following section, we introduce an estimator suitable for our purposes.

%% file: inference/inverse.tex
\subsection[Estimating the precision matrix]{Estimating $\rT$}\label{ss:inverse}
\myappendix{
\subsection{Materials required for Section~\ref{ss:inverse}}\label{proofs:inverse}
}

%%%%%%%%%%%%%%%%%%%%%%%%%%%%%%%%%%%%%%%%%%%%%%%%%%%%%%%%%%%%%%%%%%%%%%%%%%%%%%%
%%%%%%%%%%%%%%%%%%%%%%%%%%%%%%%%%%%%%%%%%%%%%%%%%%%%%%%%%%%%%%%%%%%%%%%%%%%%%%%

The one-step second-stage estimator~$\rbdb$ depends on an estimator~$\rTh$ of~$\rT=\rSd^{-1}$, the population precision matrix of the conditional means~$\dd_i$.
In general, estimating the population precision matrix incurs two main difficulties in the high-dimensional setting.
First, the empirical covariance matrix~$\rSdh$ is singular when~$\px>\n$ and cannot be inverted to produce an estimator of~$\rT$.
Second, even if an inverse were available, one cannot na\"ively use the continuous mapping theorem to derive asymptotic guarantees if~$\px\to\infty$,
since the sequence of population covariance matrices~$\rSd\equiv\rSdn$ does not itself have a limit if~$\px\to\infty$.
In addition to these general difficulties, we must further contend with the fact that the conditional mean matrix~$\dD$ is unknown.
Hence any estimator of~$\rTh$ will depend on the prediction~$\dDh$, and guarantees for such an estimator must account for such dependence.

We use a slight modification of the CLIME estimator of~\cite{CLL11} to contend with the challenges described above.
The rows~$\rth_j$ of the estimator~$\rTh$ are obtained as solutions to the CLIME program codified below.

\renewcommand{\progarg}{\bs{\theta}}
\begin{prog}[Program~for $\rth_j$]\label{program}
\begin{align*}
	\underset{\progarg\smallpad\in\R^{\px}}{\text{minimize:}}& \quad
		\obj(\progarg) \smallpad{:=} \|\progarg\|_1 \,,
	& \text{subject to:}& \quad
		\|\rSdh\progarg - \basis_j\|_\infty \smallpad\le \tf \,,
\end{align*}
where~$\basis_j$ denotes the~$j^\text{th}$ canonical basis vector in~$\px$ dimensions and~$\tf>0$ is a controlled tolerance.
\end{prog}

The present estimator~$\rTh$ differs in only one respect from that of the CLIME estimator of~\cite{CLL11}.
The latter authors symmetrize the matrix~$\rT$ with rows obtained as solutions to the aforementioned optimization problem, whereas we use the raw solutions.
We omit the symmetrization step for simplicity;
the~$\ell_\infty$ and~$\ell_1$ guarantees that~\cite{CLL11} obtain for the estimation error of the CLIME estimator continue to hold.
We include the requisite guarantees for the unsymmetrized estimator in Section~\ref{proofs:inverse} of the Supplementary Materials.

The present estimator~$\rTh$ also differs in an important respect from that of~\cite{JM14}.
The latter also obtain an inverse Gram matrix approximation as a solution to a convex program with identical constraints as in Program~\ref{program} but with objective function~$\obj(\progarg)=\En[\ip{\progarg,\dx_i}^2]$.
To our knowledge, however, it is currently unknown whether the choice of~$\obj$ in~\cite{JM14} yields guarantees comparable to those of the CLIME estimator.

The~$\ell_1$ bound for~$\rth_j-\rt_j$, which we require for control of the remainder terms~$\rem_\ell$, depends on the following restriction on the class of population precision matrices~$\rT$.
% We adopt the nomenclature and notation of~\cite{CLL11}.

\begin{defi}[Uniformity class]\label{uniformity}
Following~\cite{CLL11}, we define the \emph{uniformity class} of population precision matrices~$\rT = \rSd^{-1}$ relative to the controlled tolerance~$\tuni \in [0,1)$ and the generalized sparsity level~$\suni > 0$ by
% \begin{align}\label{def:uniformity}
% \begin{align*}
% 	&\ \, \uniclass \\
% 	:=&\ \,
% 	\bset{\rT = (\rtc_{jk})_{j,k=1}^{\px} \succ \0 : \|\rT\|_{L_1} \le \Mrt; \maxjpx \sum_{k\in[\px]} |\rtc_{jk}|^q \le \suni } \,.
% \end{align*}
\begin{equation*}
	\uniclass
	\smallpad{:=}\\
	\bset{\rT = (\rtc_{jk})_{j,k=1}^{\px} \succ \0 : \|\rT\|_{L_1} \le \Mrt; \maxjpx \sum_{k\in[\px]} |\rtc_{jk}|^q \le \suni } \,.
\end{equation*}
\end{defi}

\myappendix{
We present the properties of the estimators~$\rth_j$ required for the results of Sections~\ref{normality} and~\ref{ss:remainders}.
Lemma~\ref{lem:eesupprc}, which gives an~$\ell_\infty$ bound for the estimation error~$\rth_j-\rth$, is comparable to~\cite[Theorem 4]{CLL11}; the proofs are similar but depend on different conditions on the covariance estimator~$\rSdh$.
The proof of~Lemma~\ref{lem:eeprc} follows part of that of~\cite[Theorem 6]{CLL11}.
}
% For convenience, we refer to the three results in the following Lemma as ``H\"older's inequality'' throughout the proofs -- in particular that of Lemma~\ref{lem:eesupprc}.
% \newcommand{\mA}{\bs{A}}
% \newcommand{\mB}{\bs{B}}
% \begin{lem}\label{hoelder}
%   Let $\mA\in\R^{n\times m}$,  $\mB\in\R^{m\times p}$. Then,
%   \begin{align*}
%    (i)~~\|\mA\mB\|_\infty &= \|\mA\|_\infty\|\mB\|_1\\
%     (ii)~~\|\mA\mB\|_\infty &= \|\mA^\top\|_1\|\mB\|_\infty\ .
%   \end{align*}
% If further $n=m$ and $\mA$ symmetric,
% \begin{equation*}
%     (iii)~~\|\mA\mB\|_\infty = \|\mA\|_1\|\mB\|_\infty\ .
% \end{equation*}
% \end{lem}
% Compare also to~\cite[Lemma S.1]{Caner&Kock2018}.
%
% \begin{proof}[Proof of Lemma~\ref{hoelder}]
% The first identity can be found in~\cite[Page~44]{van2016estimation}.
% The second identity follows from
% \begin{enumerate}
% 	\item
% 		the fact that the $\ell_\infty$-matrix norm remains unchanged under transposition of the argument and
% 	\item
% 	the first identity
% \end{enumerate}
% as such:
% \begin{align*}
%   \|\mA\mB\|_\infty\ \equals \|(\mB^\top \mA^\top)^\top\|_\infty\ &\equals \|\mB^\top \mA^\top\|_\infty\\
% 	\aequals \|\mB^\top\|_\infty \|\mA^\top\|_1\ \equals \|\mB\|_\infty \|\mA^\top\|_1\ \equals  \|\mA^\top\|_1\|\mB\|_\infty\ \,.
% \end{align*}
% The third identity follows directly from the second one.
% \end{proof}
% }

%%%%%%%%%%%%%%%%%%%%%%%%%%%%%%%%%%%%%%%%%%%%%%%%%%%%%%%%%%%%%%%%%%%%%%%%%%%%%%%
%%%%%%%%%%%%%%%%%%%%%%%%%%%%%%%%%%%%%%%%%%%%%%%%%%%%%%%%%%%%%%%%%%%%%%%%%%%%%%%

\myappendix{
\begin{lem}\label{lem:eesupprc}
Suppose that: (i) the quantity~$\|\rT\|_{L_1}$ is bounded above by a constant~$\Mrt < \infty$; and (ii)~$\rTh$ is an estimate of~$\rT = \rSd^{-1} = \cov(\dd_i)^{-1}$ with rows~$\rth_j$ obtained as solutions to Program~\ref{program}.
Then, on the set~$\Tinv(\tf)$ as defined in \eqref{def:Tinv},
\benn
	\|\rth_j-\rt_j\|_\infty \smallpad\le 2\Mrt\tf
\eenn
for each~$j\in[\px]$.
\end{lem}
}%

%%%%%%%%%%
% PROOF
\myappendix{%
\begin{proof}[Proof of Lemma \ref{lem:eesupprc}]
First, observe that the conditions of the present lemma entail that
\benn
	\|\rT\rSdh-\Id\|_\infty \smallpad\le \tf\,, \qquad\qquad \|\rTh\rSdh-\Id\|_\infty \smallpad\le \tf \,.
\eenn
Now, on the set $\Tinv(\tf)$, each row $\rt_j$ is feasible for the respective Specific Program~\ref{program}.
It then follows from the optimality of $\rth_j$ that $\|\rth_j\|_1 \leq \|\rt_j\|_1$ for each $j\in[\px]$ and hence that $\maxjpx \|\rth_j\|_1 \leq \|\rT\|_{L_1}$.
Next, observe that
\begin{align*}
% 	\rTh-\rT \equals (\rTh\rSd-\Id)\rT \aequals (\rTh\rSdh + \rTh(\rSd-\rSdh) - \Id)\rT \\
% 	\aequals(\rTh\rSdh-\Id)\rT + \rTh(\rSd-\rSdh)\rT \\
% 	\aequals(\rTh\rSdh -\Id)\rT + \rTh(\Id-\rSdh\rT) \,.
	\rT - \rTh \equals \rT(\Id-\rSd\rTh) \aequals \rT(\Id + (\rSdh - \rSd)\rTh - \rSdh\rTh) \\
	\aequals \rT(\Id-\rSdh\rTh) - \rT(\rSd-\rSdh)\rTh \\
	\aequals \underbrace{\rT(\Id-\rSdh\rTh)}_{\I_1} - \underbrace{(\Id - \rT\rSdh)\rTh}_{\I_2} \,.
\end{align*}
From H\"older,
\begin{align*}
\|\I_2\|_\infty \smallpad\le \|\Id - \rT\rSdh\|_\infty\|\rTh\|_{L_1} \smallpad\le \Mrt\tf\,.
\end{align*}
The matrix $\ell_\infty$- and ($L_1$- norm if the argument is symmetric) are invariant under transposition of their arguments.
Use this fact and H\"older to obtain
\begin{align*}
	\|\I_1\|_\infty \aequals \|(\Id-\rSdh\rTh)^\top\rT^\top\|_\infty \\
	\smallpad{&\le} \|(\Id-\rSdh\rTh)^\top\|_\infty\|\rT^\top\|_{L_1} \\
	\smallpad{&\le} \|\Id-\rSdh\rTh\|_\infty\|\rT\|_{L_1} \ \le \ \Mrt\tf \,,
\end{align*}
where the final line follows from the fact that both~$\Id-\rSdh\rTh$ and~$\rT$ are symmetric.
Thus
\begin{equation*}
	\|\rT-\rTh\|_\infty \smallpad\le \|\I_1\|_\infty + \|\I_2\|_\infty \smallpad\le 2\Mrt\tf\,,
\end{equation*}
as required.
\end{proof}
}%

%%%%%%%%%%%%%%%%%%%%%%%%%%%%%%%%%%%%%%%%%%%%%%%%%%%%%%%%%%%%%%%%%%%%%%%%%%%%%%%
%%%%%%%%%%%%%%%%%%%%%%%%%%%%%%%%%%%%%%%%%%%%%%%%%%%%%%%%%%%%%%%%%%%%%%%%%%%%%%%

\myappendix{
\begin{lem}\label{lem:eeprc}
Suppose in addition to the conditions of Lemma \ref{lem:eesupprc} that~$\rT$ belongs to the uniformity class~$\uni(\Mrt, \tuni, \suni)$.
Then,
\benn
	\|\rth_j-\rt_j\|_1 \relates{\le} 2\cprcr(2\Mrt\tf)^{1-\tuni}\suni
\eenn
for each $j\in[\px]$, where $\cprcr := 1+2^{1-\tuni} + 3^{1-\tuni}$.
\end{lem}
}%

%%%%%%%%%%
% PROOF
\myappendix{%
\begin{proof}[Proof of Lemma \ref{lem:eeprc}]
See the proof of line (14) of \cite[Theorem 6]{CLL11}.
\end{proof}
}%

%%%%%%%%%%%%%%%%%%%%%%%%%%%%%%%%%%%%%%%%%%%%%%%%%%%%%%%%%%%%%%%%%%%%%%%%%%%%%%%
%%%%%%%%%%%%%%%%%%%%%%%%%%%%%%%%%%%%%%%%%%%%%%%%%%%%%%%%%%%%%%%%%%%%%%%%%%%%%%%

% As presented in this essay, the $\ell_\infty$ and $L_1$ rates of convergence for estimators $\rTh$ as described above depend on the feasibility of the population precision matrix $\rT$ for Program~\ref{program}.
\noindent In the sequel, we assume as part of high-level regularity conditions that~$\rT\in\uniclass$ and that the model parameters~$\Mrt$ and~$\suni$ are well-behaved as functions of~$\n$.
These parameters appear in the rates for the remainder terms in our analysis of the two-stage Lasso of Section~\ref{estimation}.
For high-level results, we also assume that the probability that the
rows~$\rt_j$ of the population precision matrix are feasible for Program~\ref{program} approaches one.
To express this requirement formally, we define the event
\be\label{def:Tinv}
	\Tinv(\tf)
	\smallpad{:=}
	\bset{\|\rT\rSdh-\Id\|_\infty\le\tf}
\ee
where~$\tf>0$ is the tolerance of Program~\ref{program}, and require that~$\Prb\,\Tinv(\tf)\to1$ as~$\n\to\infty$.
We identify a theoretical choice of~$\tf$ that satisfies the latter requirement in Lemma~\ref{feassgg}; we discuss a practical method for selecting~$\tf$ in Section~\ref{s:numerics}.
Note that, given the event~$\Tinv(\tf)$, the rows~$\rt_j$ of~$\rT_j$ are each feasible for the respective Program~\ref{program}.

Since the quantity~$\tf$ appears in the rates for the remainder terms~$\rem_\ell$, it must be chosen carefully so as to balance the growth of~$\Prb\,\Tinv(\tf)$ with the decay of the~$\|\rem_\ell\|_\infty$.
The appropriate choice of~$\tf$ depends on both the distribution of the~$\dz_i$ as well as the rate for~$\maxerrfs$ --- see Lemma~\ref{feassgg}.

% Given sub-Gaussian~$\dz_i$ and the choice of~$\tf$ described in Lemma~\ref{feassg}, we may ensure~$\limni\Prb\Tinv(\tf)=1$ by satisfying three conditions.
% % First, we require that~$\sqrt{\log(\px)/\n}=o(1)$, which is a standard growth assumption for consistency of high-dimensional estimators.
% First, we must choose~$\gep>0$ so that~$\Prb\set{\maxerrfs>\gep}=o(1)$.
% Such a choice depends on the properties of the specific estimator~$\rAh$; we demonstrate in Lemma~\ref{feassgg} that choosing~$\gep=O(\sqrt{\log(\pz)/\n})$ suffices when the first-stage estimators~$\rah^j$ are obtained via the Lasso program and~
% Third, we must specify~$\ta$ in the tolerance~$\tf$ so that~$\ta^2/(6e^2\subgauss'^2)>2$.
% To do so in practice would require knowledge of quantities such as the sub-Gaussian norm~$\subgauss$, which is not feasible.
% We note that we are not concerned with the feasibility of such quantities in this essay; we explain this position in Section~\ref{bounds}, where we discuss tuning parameter selection.
% In Section~\ref{s:numerics}, we discuss a practical scheme for selecting the quantity~$\tf$ that gives good empirical results.

%% file: inference/normality.tex
\subsection{Asymptotic normality}\label{normality}
\myappendix{
\subsection{Materials required for Section~\ref{normality}}\label{proofs:normality}
}

%%%%%%%%%%%%%%%%%%%%%%%%%%%%%%%%%%%
% New start of section

We saw in Section~\ref{endogeneity} that the updated estimator~$\rbdb$ satisfies~$\sqrt{n}(\rbdbj-\rbj) = \sqrt{\n}\En[\ip{\rt_j,\dd_i}\dui] + \sum_{\ell=1}^4\remj{\ell}$.
If the remainder terms vanish in probability, then~$\sqrt{\n}(\rbdbj-\rbj)$ shares the same weak limit, if it exists, as~$\sqrt{\n}\En[\ip{\rt_j,\dd_i}\dui]$.
If the model were fixed in~$\n$, the Central Limit Theorem would entail that the latter quantity converges weakly to a~$\Normal(0,\wjse^2)$, where~$\wjse^2=\E[\ip{\rt_j,\dd_i}^2\dui^2]$.
In Theorem~\ref{wl}, we provide conditions under which~$\sqrt{n}(\rbdbj-\rbj)/\wjse$ converges weakly to a standard Normal random variable when the model is not fixed in $\n$.
We also show that the limit continues to hold if~$\wjse$ is replaced by an estimator~$\wjseh$ that satisfies~$|\wjseh-\wjse|=o_{\Prb}(1)$.
Note that Theorem~\ref{wl} gives conditions under which the limit holds given homoscedastic Gaussian noise (Condition~\ref{wl:noise1}) as well as conditions under which the limit holds given generic i.i.d. noise (Condition~\ref{wl:noise2}).
Condition~\ref{wl:noise1} is unnecessarily restrictive in practice.
Nonetheless, we include the result under Condition~\ref{wl:noise1} because it serves as a benchmark and requires weaker subsequently weaker assumptions concerning limiting rates.
The latter result includes the case of Assumption~\ref{assu:subgaussian}, as well as any other i.i.d. second-stage noise regime.

\begin{thm}[Weak limits]\label{wl}
Suppose that
\begin{enumerate*}[label=(\arabic{*})]
  \item\label{wl:cond1}
   % there exists a constant sequence ~$c_n=o(1)$ such that for~$\rSdb=\dD^\top\dD/\n$, we have $\Prb\set{\|\rSdb-\rSd\|_\infty > c_n} = o(1)$;
   the quantity~$\|\rSdb-\rSd\|_\infty$ vanishes in probability, where~$\rSdb:=\dD^\top\dD/\n$;
  \item\label{wl:cond2}
    the remainder terms satisfy~$\|\rem_{\ell}\|_\infty=o_{\Prb}(1)$ for each~$1\le\ell\le4$;
  \item\label{wl:cond3}
    $\rT_{jj}>\Tjjmin$ for some universal constant~$\Tjjmin>0$ and each~$j\in[\px]$;
  \item\label{wl:cond4}
    $\maxjpx\|\rt_j\|_1 \le \Mrt$ for some universal constant~$\Mrt<\infty$.
\end{enumerate*}

If either
\begin{enumerate*}[label=(\arabic{*})]
  \setcounter{enumi}{4}
  \item\label{wl:noise1}
    the noise elements~$\dui$ satisfy~$\dui\,|\,\dz_i\sim_{i.i.d.}\Normal(0,\sdu)$, where $\sdu$ is bounded away from zero and infinity uniformly in~$\n$, or
  \item\label{wl:noise2}
    the~$\dz_i$ and~$\dui$ are i.i.d. with~$\E[\dui^{2}|\dz_i]=\sdu^2$, where~$\sdu$ is bounded away from zero and infinity uniformly in~$\n$, and there exist
    $0<\zeta<1/2$ and~$\wlpartwo > 0$ such that
    \begin{enumerate*}
      \item\label{wl:noise2a}
        $\Prb\bset{|\ip{\rt_j,\dd_i}|>\n^{\zeta}} \equals o(1)$ and
      \item\label{wl:noise2b}
        $\E[|\dui|^{2+\wlpartwo}] \smallpad\lesssim \sdu^{2+\wlpartwo}$,
    \end{enumerate*}
\end{enumerate*}
then
\benn
  \sqrt{\n}(\rbdbj-\rbj)/\wjse \smallpad\leadsto \zj\sim\Normal(0,1) \,.
\eenn
Furthermore, the limit continues to hold if~$\wjse$ is replaced by an estimator~$\wjseh$ that satisfies~$|\wjseh-\wjse|=o_{\Prb}(1)$.
\end{thm}

\myappendix{
\begin{proof}[Proof of Theorem~\ref{wl}]
The proof of the first claim consists of two steps.
First, we show that the quantity
\benn
  \zjn \smallpad{:=} \frac{1}{\sqrt{\n}}\sum_{i=1}^{\n} \ip{\rt_j,\dd_i}\dui/\wjse\,.
\eenn
and~$\sqrt{\n}(\rbdbj-\rbj)/\wjse$ share the same weak limit.
Second, we show that~$\zjn\leadsto\Normal(0,1)$.
To establish the first step, we claim that
\be\label{step1g}
  \lim_{\n\to\infty}\Prb\set{\sqrt{\n}(\rbdbj-\rbj)/\wjse \le t}
  \smallpad{\le}
  \lim_{\n\to\infty}\Prb\set{\zjn \le t}\,,
\ee
for all~$t\in\R$.
An analogous lower bound follows by a matching argument, which shows that~$\sqrt{\n}(\rbdbj-\rbj)/\wjse$ and~$\zjn$ share the same weak limit.
To show the claim above, let~$t\in\R$ be given, fix a controlled~$\gep>0$, and note that, by the decomposition of Lemma~\ref{lem:one-step}, we have
\begin{align*}
  \Prb\bset{\sqrt{\n}(\rbdbj-\rbj)/\wjse \le t}
  \smallpad{&\le}
  % \Prb\left\{\frac{1}{\sqrt{\n}}\sumin\ip{\rt_j,\dd_i}\dui/\wjse
  \Prb\left\{\zjn
    + \sum_{\ell=1}^{4}\remj{\ell}/\wjse \leq t + 4\gep\right\} \\
  \smallpad{&\le}
    \Prb\bset{\zjn \le t + \gep} + \sum_{\ell=1}^4\Prb\bset{\remj{\ell}/\wjse > \gep} \,.
\end{align*}
By specification of~$\sdu$ and~$\rTjj$ in Conditions~\ref{wl:noise1} and~\ref{wl:noise2} and Condition~\ref{wl:cond3} of the present theorem, it follows that~$\wjse$ is bounded strictly away from 0 uniformly in~$\n$.
The assumptions of the present theorem then entail that $\Prb\set{\remj{\ell}/\wjse>\gep}=o(1)$ for all~$\gep>0$ and each~$\ell$.
Letting~$\gep$ tend to zero shows the claim of~\eqref{step1g}.
It follows from the analogous lower bound that
\benn
  \lim_{\n\to\infty}\Prb\set{\sqrt{\n}(\rbdbj-\rbj)/\wjse \le t}
  \equals
  \lim_{\n\to\infty}\Prb\set{\zjn \le t}
\eenn
for all~$t\in\R$, thus completing the first step.

Next, we show that, under each of Conditions~\ref{wl:noise1} and~\ref{wl:noise2} in the statement of the present theorem,~$\zjn\leadsto\zj\sim\Normal(0,1)$.
To this end, we define the quantity
\benn
  \wj^2 \smallpad{:=} \rt_j^\top\rSdb\rt_j
  \equals
  \fon \sumin \ip{\rt_j,\dd_i}^2\,,
\eenn
where we recall that~$\rSdb=\dD^\top\dD/\n$.
We claim first that
\benn
  \zjnt \smallpad{:=} \frac{1}{\sqrt{\n}}\sumin\frac{\ip{\rt_j,\dd_i}\dui}{\wj\sdu}
  \smallpad\leadsto\zj\sim\Normal(0,1)
\eenn
under each of Conditions~\ref{wl:noise1} and~\ref{wl:noise2} and second that $\sdu\wj/\wjse \to_{\Prb} 1$.
Note that~$\zjn=\frac{\wj\sdu}{\wjse}\zjnt$, hence the desired limit follows from an application Slutsky's Lemma.

To show the first claim under Condition~\ref{wl:noise1}, note that, by specification of~$\wj$, we have under Assumption~\ref{assu:subgaussian} that
\benn
  \frac{1}{\sqrt{\n}}\sumin\frac{\ip{\rt_j,\dd_i}\dui}{\wj\sdu} \,|\, \dZ
  \smallpad\sim\Normal(0,1)\,.
\eenn
Thus
\begin{equation*}
  \limni \Prb\set{\zjnt\le t}
  \equals
  \limni \E\big[\Prb\set{\zjnt\le t\,|\,\dZ}\big]
  \equals
  \limni \E[\Phi(t)\,|\,\dZ]
  \equals
  \Phi(t)
\end{equation*}
for all~$t\in\R$, where~$\Phi$ denotes the c.d.f. of a standard Normal random variable.
This shows the desired weak limit under Condition~\ref{wl:noise1}.

We use the Lindeberg-Feller Central Limit Theorem to show the limit under Condition~\ref{wl:noise2}.
To begin, write
\begin{equation*}
  \zjnt\equals\frac{1}{\sqrt{n}}\sumin\xii\,,
  \qquad
  \xii\smallpad{:=}\ip{\rt_j,\dd_i}\dui/(\wj\sdu)\,.
\end{equation*}
Note that
\benn
  \E[\xii] \equals \E\big[\ip{\rt_j,\dd_i}/(\wj\sdu)\E[\dui|\dz_i]\big]
  \equals 0
\eenn
and that
\benn
  \sd_{\n}^2 \smallpad{:=}
  \sumin \E[\xii^2]
  \equals
  \E\left[\frac{1}{\wj^2}\sumin\ip{\rt_j,\dd_i}^2\E[(\dui/\sdu)^2|\dz_i]\right]
  \equals \n \,.
\eenn
To demonstrate Lindeberg's condition, let~$\gd>0$ be arbitrary and write
\begin{align*}
  \sd_{\n}^{-2} \sumin \E[\xii^2\Ind{|\xii|>\gd\sd_{\n}}]
  \aequals
  \fon\sumin \E\big[\E[\xii^2\Ind{|\xii|>\gd\sqrt{\n}}\,|\,\dz_i]\big] \\
  \aequals
  \E\left[\fon\frac{1}{\wj^2}\sumin\ip{\rt_j,\dd_i}^2\E[(\dui/\sdu)^2\Ind{|\xii|>\gd\sqrt{\n}}\,|\,\dz_i]\right] \,,
\end{align*}
where we have substituted the definition of $\xii$ and extracted the factor $\ip{\rt_j,\dd_i}^2/\wj^2$ from the conditional expectation to obtain the second line.
Now substitute the definition of $\wj$ below and note
\begin{equation*}
  \frac{1}{\wj^2}\sumin\ip{\rt_j,\dd_i}^2 \equals \frac{n}{\sumin\ip{\rt_j,\dd_i}^2}
    \times \sumin\ip{\rt_j,\dd_i}^2 \equals \n \,.
\end{equation*}
Substitute the above result into the second line of two displays prior and use the law of iterated expectation to obtain
\begin{equation*}
  \sd_{\n}^{-2} \sumin \E[\xii^2\Ind{|\xii|>\gd\sd_{\n}}]
  \equals \E\big[(\duc_1/\sdu)^2\Ind{|\xio|>\gd\sqrt{\n}}\big]\,,
\end{equation*}
where we cite condition \ref{wl:noise2} that the~$\dz_i$ and~$\dui$ are i.i.d. with~$\E[\dui^{2}|\dz_i]=\sdu^2$ to reduce to the case of $\duc_1$.
For brevity, we write~$\duit:=\duc_{1}/\sdu$.
Introduce the set~$\Ti:=\set{|\ip{\rt_j,\dd_1}|\leq\n^{\zeta}}$ and note, since~$|\xio|\le\wj^{-1}|\ip{\rt_j,\dd_i}||\duit|$, that
\benn
  \bset{|\xio|>\gd\sqrt{\n}}\cap\Ti
  \smallpad\subseteq
  \bset{|\duit|>\gd\wj\n^{1/2-\zeta}}\,,
\eenn
and hence that~$\Inds{\{|\xio|>\gd\sqrt{\n}\}\cap\Ti}\le\Ind{|\duit|>\gd\wj\n^{1/2-\zeta}}$.
Combine this inequality and the law of total probability with the result of two displays previous to obtain
\begin{align*}
  \sd_{\n}^{-2} \sumin \E[\xii^2\Ind{|\xio|>\gd\sd_{\n}}]
  \aequals
  \E\big[\duit^2(\Inds{\{|\xio|>\gd\sqrt{\n}\}\cap\Ti}+\Inds{\{|\xio|>\gd\sqrt{\n}\}\cap\Ti^\complement}) \big] \\
  \aequals
  \underbrace{\E[\duit^2 \Ind{|\duit|>\gd\wj\n^{1/2-\zeta}}]}_{\I_1} + \underbrace{\E[\duit^2\Inds{\Ti^\comp}]}_{\I_2} \,,
\end{align*}
where the substitution of indicators in the final line above is permitted since~$\duit^2\geq0$.
To show Lindeberg's condition, it sufficies to show that~$\I_1$ and~$\I_2$ are each~$o(1)$.

To treat~$\I_1$, consider the event~$\set{\wj\le\Tjjmin^{1/2}/\sqrt{2}}$ and write
\begin{align*}
  \Ind{|\duit|>\gd\wj\n^{1/2-\zeta}}
  \aequals
  \Inds{\set{|\duit|>\gd\wj\n^{1/2-\zeta}}\cap\set{\wj\le\Tjjmin^{1/2}/\sqrt{2}}}\\
  &\qquad\qquad +
  \Inds{\set{|\duit|>\gd\wj\n^{1/2-\zeta}}\cap\set{\wj>\Tjjmin^{1/2}/\sqrt{2}}} \\
  \smallpad{&\le}
  \Ind{\wj\le\Tjjmin^{1/2}/\sqrt{2}} + \Ind{|\duit|>\gd\Tjjmin^{1/2}\n^{1/2-\zeta}}
\end{align*}
so that
\begin{equation*}
  \I_1 \smallpad\le
  \underbrace{\E[\duit^2\Ind{\wj\le\Tjjmin^{1/2}/\sqrt{2}}]}_{\I_{1\mathrm{a}}}
  +
  \underbrace{\E[\duit^2\Ind{|\duit|>\gd\Tjjmin^{1/2}\n^{1/2-\zeta}}]}_{\I_{1\mathrm{b}}} \,.
\end{equation*}
Observe that
\begin{align*}
  \I_{1\mathrm{a}} \aequals \E\big[\Ind{\wj\le\Tjjmin^{1/2}/\sqrt{2}}\E[\duit^2\,|\,\dz_i]\big] \\
  \smallpad{&\lesssim}
    \Prb\set{\wj\le\Tjjmin^{1/2}/\sqrt{2}} \\
  \smallpad{&=}
    \Prb\set{\wj^2-\rT_{jj} \le \Tjjmin/2-\rT_{jj}} \\
  \smallpad{&\le}
    \Prb\set{\wj^2-\rT_{jj} \le -\Tjjmin/2} \\
  \smallpad{&\le}
    \Prb\set{|\wj^2-\rT_{jj}| \geq \Tjjmin/2} \,,
\end{align*}
where the inference to the second line of the display above follows from specification of $\E[\duit^2\,|\,\dz_i]$ in condition \ref{wl:noise2} of the present Theorem.
Now note that
\begin{align*}
  |\wj^2-\rT_{jj}| \aequals
  |\rt_j^{\top}\rSdb\rt_j - \rT_{jj}| \\
  \aequals
  |\rt_j^\top(\rSdb-\rSd)\rt_j|
  \smallpad{\le}
  \|\rt_j\|_1^2\|\rSdb-\rSd\|_\infty \,.
\end{align*}
Thus
\begin{align}\label{wl:turn1a}
  \I_{1\mathrm{a}} \smallpad{&\le} \Prb\set{|\wj^2-\rT_{jj}| \geq \Tjjmin/2} \nonumber\\
  \smallpad{&\le}
    \Prb\bset{\|\rSdb-\rSd\|_\infty \geq \Tjjmin/(2\|\rt_j\|_1^2)}
  \equals
    o(1)
\end{align}
by Conditions~\ref{wl:cond1} and~\ref{wl:cond4} of the present theorem.

We now show that
\begin{equation*}
  \I_{1\mathrm{b}} \equals \E[\duit^2\Ind{|\duit|>\gd\Tjjmin^{1/2}\n^{1/2-\zeta}}] \equals o(1)\,.
\end{equation*}
Consider the event $\{ |\duit|>\gd\Tjjmin^{1/2}\n^{1/2-\zeta} \}$,
which is the index of the indicator function above.
Exponentiate both sides of the inequality in that set by $\wlpartwo$ to obtain that
\begin{align*}
  \bset{|\duit| > \gd\Tjjmin^{1/2}\n^{1/2-\zeta}}
  \smallpad{&=}
    \bset{|\duit|^\wlpartwo > (\gd\Tjjmin^{1/2})^{\wlpartwo}\n^{\wlpartwo(1/2-\zeta)}} \,.
  % \aequals
    % \set{\sqrt{2}(\gd\Tjjmin^{1/2})^{-\wlpartwo}\n^{-\wlpartwo(1/2-\zeta)}|\duit|^\wlpartwo > 1},
\end{align*}
Conclude of the respective indicator functions that
\benn
  \Ind{|\duit|>\gd\Tjjmin^{1/2}\n^{1/2-\zeta}}
  \smallpad=
  \Ind{|\duit|^\wlpartwo > (\gd\Tjjmin^{1/2})^{\wlpartwo}\n^{\wlpartwo(1/2-\zeta)}}\,.
  % \Ind{(\gd\Tjjmin^{1/2})^{-\wlpartwo}\n^{-\wlpartwo(1/2-\zeta)}|\duit|^\wlpartwo > 1} \,.
  % \big(\sqrt{2}(\gd\Tjjmin^{1/2})^{-\wlpartwo}\n^{-\wlpartwo(1/2-\zeta)}|\duit|^\wlpartwo\big)\,.
\eenn
Note by direct manipulation that, given the event that indexes the right-hand indicator above, it holds that
\begin{equation*}
  (\gd\Tjjmin^{1/2})^{-\wlpartwo}\n^{-\wlpartwo(1/2-\zeta)}|\duit|^\wlpartwo > 1 \,.
\end{equation*}
Thus, on the support of the right-hand side two displays prior, it holds that
\begin{equation*}
  \Ind{(\gd\Tjjmin^{1/2})^{-\wlpartwo}\n^{-\wlpartwo(1/2-\zeta)}|\duit|^\wlpartwo > 1}
  \smallpad\le
  \big(\gd\Tjjmin^{1/2}\big)^{-\wlpartwo}\n^{-\wlpartwo(1/2-\zeta)}|\duit|^\wlpartwo \,.
\end{equation*}
% In turn, on the support of the indicator~$\Ind{|\duit|>\gd\Tjjmin^{1/2}\n^{1/2-\zeta}}$ in
% % \begin{equation*}
%   $\I_{1\mathrm{b}} \equals \E[\duit^2\Ind{|\duit|>\gd\Tjjmin^{1/2}\n^{1/2-\zeta}}]$,
% % \end{equation*}
% we have
% \begin{equation*}
%   \Ind{|\duit|>\gd\Tjjmin^{1/2}\n^{1/2-\zeta}}
%   \smallpad\le
%   \big(\gd\Tjjmin^{1/2}\big)^{-\wlpartwo}\n^{-\wlpartwo(1/2-\zeta)}|\duit|^\wlpartwo \,.
% \end{equation*}
Hence, by the monotonicity and linearity of expectations,
\begin{align*}
  \I_{1\mathrm{b}} \aequals \E[\duit^2\Ind{|\duit|>\gd\Tjjmin^{1/2}\n^{1/2-\zeta}}] \\
  \smallpad{&\le} \E\big[ \duit^2\big(\gd\Tjjmin^{1/2})^{-\wlpartwo}\n^{-\wlpartwo(1/2-\zeta)}|\duit|^\wlpartwo \big]
  \equals (\gd\Tjjmin^{1/2})^{-\wlpartwo}\n^{-\wlpartwo(1/2-\zeta)}\E[|\duit|^{2+\wlpartwo}] \,.
\end{align*}
Condition~\ref{wl:noise2} of the present theorem stipulates that $\E[|\duit|^{2+\wlpartwo}]$ is of constant order,
hence,
\begin{align*}
  \limni \I_{1\mathrm{b}}
  \smallpad{&\lesssim}
  \limni (\gd\Tjjmin^{1/2})^{-\wlpartwo}\n^{-\wlpartwo(1/2-\zeta)}
  \equals 0
\end{align*}
by the specification that $\wlpartwo > 0$ and $\zeta < 1/2$ in Condition \ref{wl:noise2} of the present theorem.

To show as much for~$\I_2$, observe that
\begin{align*}
  \limni \I_2 \aequals \limni \E\big[\Ind{|\ip{\rt_j,\dd_i}|>\n^{\zeta}}\E[\duit^2\,|\,\dz_i]\big] \\
  \smallpad{&\lesssim}
  \limni \Prb\set{|\ip{\rt_j,\dd_i}|>\n^{\zeta}} \equals 0
\end{align*}
by Condition~\eqref{wl:noise2a} of the present theorem.
This concludes the demonstration of Lindeberg's condition.
It follows that~$\zjnt\leadsto\Normal(0,1)$.
To show as much for~$\zjn$ and hence for~$\sqrt{\n}(\rbdbj-\rbj)/\wjse$, it suffices to show that~$\wj\sdu/\wjse\to_{\Prb}1$.
Note that~$\wjse=\sdu\Tjj$ and hence that~$\wj\sdu/\wjse=\wj/\Tjj$.
Since~$\Tjj$ is bounded strictly away from~zero uniformly in~$\n$, we have~$|\wj/\Tjj-1|=|\wj-\Tjj|/\Tjj$ and hence that it suffices to show that~$|\wj-\Tjj|=o_{\Prb}(1)$.
But, as we established above, we have for arbitrary $\gep>0$
\begin{align*}
  \Prb\set{|\wj-\Tjj| > \gep}
  \smallpad\le
  \Prb\set{\|\rSdb-\rSd\|_\infty > \gep/\Mrt^2} \equals o(1)
\end{align*}
by Condition~\ref{wl:cond1} of the present theorem.
Thus~$\wj\sdu/\wjse\to_{\Prb}1$ hence~$\zjn\leadsto\Normal(0,1)$ under each Condition~\ref{wl:noise1} and~\ref{wl:noise2} of the present theorem.

It remains to show that the limit holds when~$\wjse$ is replaced by an estimator~$\wjseh$ that satisfies~$|\wjseh-\wjse|=o_{\Prb}(1)$.
Suppose that~$\wjseh$ is such an estimator.
We claim that~$\wjseh/\wjse-1=o_{\Prb}(1)$.
To see as much, note first that, from the hypotheses of the present theorem,~$\wjse=\sdu\Tjj$ is bounded strictly away from zero uniformly in~$\n$.
It then follows that~$|\wjseh/\wjse-1|=|\wjseh-\wjse|/\wjse=o_{\Prb}(1)$ by the specification of~$\wjseh$.
By the continuous mapping theorem, it holds that~$\wjse/\wjseh\to_{\Prb}1$ and hence by Slutsky's Lemma that
% that ~$\sqrt{\n}(\rbdbj-\rbj)/\wjseh=(\wjse/\wjseh)\zjn\leadsto\Normal(0,1)$, as claimed.
\begin{equation*}
  \sqrt{\n}(\rbdbj-\rbj)/\wjseh=(\wjse/\wjseh)\zjn\leadsto\Normal(0,1) \,,
\end{equation*}
as claimed.
\end{proof}
}

% \noindent \new{Assuming also $\E[\dui^{2}|\dz_i]=\sdu^2$, \cite{BCCH12} show that their high-dimensional IV estimator is semiparametrically efficient, albeit only for fixed dimensional structural models.
% Theorem~\ref{wl} extends this result:  it shows that by estimating optimal instruments $\dd_i$ with homoscedastic structural errors, asymptotic semiparametric efficiency can be achieved even if both stages are high-dimensional.}

\noindent The proof of Theorem~\ref{wl} can be found in Section~\ref{proofs:normality} of the Supplementary Materials.
We note that, unlike the setting of~\cite[Sections 2.5]{BKRW98}, we do not require~$\sqrt{\n}$-consistency of the initial estimator~$\rbla$ in order for the updates~$\rbdbj$ to be asymptotically Gaussian.
Indeed, the rates of convergence required are dictated by the strategies used to bound the quantities~$\|\rem_\ell\|_\infty$.
See Section~\ref{ss:remainders} for an example of sufficient rates for the two-stage Lasso routine studied in Section~\ref{estimation}.

The main application of Theorem~\ref{wl} is the construction of asymptotically valid confidence intervals under a wide variety of noise regimes.
Given~$j\in[\px]$ and a confidence level~$\tau$, an asymptotic $100(1-\tau)\%$ confidence interval $\hat{\mathcal{I}}_{\tau,j}$ is given by
\be\label{confidence}
	\hat{\mathcal{I}}_{\tau,j} \smallpad{:=} \big[\rbdbj - z_\tau \wjseh, \, \rbdbj + z_\tau \wjseh\big]  \,,
\ee
where $z_\tau = \Phi^{-1}(1-\tau/2)$ and~$\wjseh$ satisfies the conditions of Theorem~\ref{wl}.
We present a simulation study of the finite-sample properties of this procedure for the updated two-stage Lasso estimator in Section~\ref{s:numerics}.

\new{Similarly as for the de-biasing approach in linear regression~\cite[Corollary~2.1]{vdGetal14},
the results in Theorem~\ref{wl} hold uniformly over a class of parameters as long as the assumptions also hold uniformly over that class.
This requires in particular that the remainder terms are $o_{\Prb}(1)$ uniformly.
In the Lasso case considered in Section~\ref{estimation}, such uniformity is known for parameter sets with bounded sparsity --- we refer again to \cite{vdGetal14} for details.
}

\new{Note of the scale factor $\wjse$ that
\begin{equation*}
  \wjse^2 \equals
    \E[\ip{\rt_j,\dd_i}^2\dui^2]
  \equals
    \sdu^2\rt_j^\top \E[\dd_i\dd_i^\top] \rt_j
  \equals
    \sdu^2 \rt_j^\top \rSd \rt_j
  \equals
    \sdu^2 \rT_{jj} \,,
\end{equation*}
which, if the second-stage model were fixed in~$\n$, would equal the asymptotic variance of the optimal linear IV estimator~\cite{A74, chamberlain1987asymptotic, newey1990efficient}.
Indeed, \cite{BCCH12} show that, in such a case, the linear IV estimator still attains the semi-parametric efficiency bound if the first-stage mean model is unknown but well-approximated with a high-dimensional linear model.
While the display above suggests a similar optimality result for our estimator as well, such a property is difficult to state formally in the present framework.
The aforementioned authors study target model parameters assumed to be identical in~$\n$ or that converge to some fixed limit.
In the high-dimensional setting of the present paper, the number of model parameters is allowed to grow infinitely and thus cannot converge to a fixed limit within the parameter space for any given model.
(Note that this is not the same issue as contending with an infinite-dimensional nuisance parameter.)
Thus it is impossible to compare asymptotic covariance matrices directly.
Indeed, we do not prove any results concerning the asymptotic variance of our estimator but rather that of the derived t-statistic in Theorem~\ref{wl}.

One might be tempted to impose regularity conditions that would allow for such comparisons.
For instance, one might require that, for any~$p_0 \in \mathbb{N}$, the~$n$-indexed sequence of parameter subsets~$\beta_{1,n},\ldots, \beta_{p_0,n}$ converge to a well-defined~$p_0$-dimensional vector.
Then one might conjecture along the lines of:
``For each such~$p_0$, the respective~$p_0$-length subset of~$\tilde{\bs\beta}$ has asymptotic variance identical to the optimal variance in a fixed~$p_0$-dimensional linear IV model.''
However, the scope of such asymptotic statements is relative to a \emph{fixed}~$p_0$;
the approach does not describe the asymptotic behavior of~$\rbdb$ under \emph{joint} growth of~$\n$ and~$\px$ when the latter may depend on the former, which is the setting we study.
Hence we do not pursue this approach here.
Instead, we believe that a proper study of optimal inference in high-dimensional linear IV models should follow along the lines of \cite{Jankova16}.
We conjecture that the scale factor in the asymptotic pivot of Theorem~\ref{wl} is optimal in the sense of~\cite[Theorem 3]{Jankova16} and that the updated second-stage estimator achieves a similar efficiency bound.
A more thorough investigation of this matter is required for future work.}

We conclude the present Section with a brief discussion of the feasibility of select conditions of Theorem~\ref{wl}.
Condition~\ref{wl:cond3} can be derived as a consequence of the standard assumption that the minimal and maximal eigenvalues of~$\rSd$ be bounded strictly away from zero and infinity; see Proposition~\ref{rTjjorder}.
The feasibility of Conditions~\ref{wl:cond1} and~\ref{wl:noise2} depends on the distribution of the conditional means~$\dd_i$ and hence of the instrumental variables~$\dz_i$.
The following Lemma shows that both conditions are satisfied if the~$\dz_i$ are sub-Gaussian.

\begin{lem}[Feasibility of Conditions~\ref{wl:cond1} and~\ref{wl:noise2a}]\label{satisfied}
Suppose that (i)~the instrumental variables~$\dz_i$ satisfy Assumption~\ref{assu:dzi}, (ii) we have~$\|\rSz\|_\infty=O(1)$, (iii) we have~$\|\rT\|_{L_1}\le\Mrt$ for~$\Mrt=O(1)$,
% (iv)~$\sqrt{\log\px}/\n^{2\zeta-b}\lesssim1$ for some~$0<b<2\zeta$.
(iv)~$\sqrt{(\log\px)/\n}=o(1)$, and (v)~$\exp(-\n^{2\zeta}/\sqrt{\log\px})=o(1)$ for some~$0<\zeta<1/2$.
Then Conditions~\ref{wl:cond1} and~\ref{wl:noise2a} of Theorem~\ref{wl} are satisfied.
\end{lem}

\myappendix{
\begin{proof}[Proof of Lemma~\ref{satisfied}]\label{proof:satisfied}
We first show that, under the conditions of the present Lemma,
\benn
  \|\rSdb-\rSd\|_\infty = o_{\Prb}(1) \,,
\eenn
thereby demonstrating Condition~\ref{wl:cond1}.
Note that Lemma~\ref{concentration} gives
\benn
  \Prb\bset{\|\rSdb-\rSd\|_\infty > a\sqrt{(\log(\px\vee\n))/\n}}
  \smallpad\le
  2(\px\vee\n)^{2-a^2/(6e^2\sendb^2)}
\eenn
for~$\n$ sufficiently large, where~$a$ is a controlled constant and~$\sendb=\Mra^2(2\sgnz^2+\|\rSz\|_\infty/\log2)$.
Condition~\ref{wl:cond1} then follows given that~$\sqrt{(\log\px)/\n}=o(1)$ and by choosing~$a=O(1)$ large enough so that~$a^2/(6e^2\sendb^2) > 2$.

We now show that, under the conditions of the present Lemma,~$\Prb\set{|\ip{\rt_j,\dd_i}|>\n^{\zeta}} = o(1)$ for~$0<\zeta<1/2$ and hence that Condition~\ref{wl:noise2a} is satisfied.
First, write
\begin{align*}
  \Prb\bset{|\ip{\rt_j,\dd_i}|>\n^{\zeta}}
  \smallpad{\le}
  \Prb\bset{\|\rt_j\|_1\|\dd_i\|_\infty > \n^{\zeta}}
  \smallpad{\le}
  \Prb\bset{\Mrt\|\dd_i\|_\infty > \n^{\zeta}}\,.
\end{align*}
Next, cite~\cite[Chapters~2.1.3 and~2.2, Pages~90--91 and~95--97]{V96} and the proof of Lemma~\ref{concentration} to infer that
\begin{align*}
  \|\maxjpx|\ddij|\|_{\sgn} \smallpad{\le} \Cmaxsgn\sqrt{\log\px+1}\maxjpx\|\ddij\|_{\sgn}
  \smallpad{\le}
  \Cmaxsgn\sgnz\Mra\sqrt{\log\px+1}  \,,
\end{align*}
where~$\Cmaxsgn$ is an absolute constant.
The exponential Markov bound then yields
\begin{align*}
  \Prb\bset{\|\dd_i\|_\infty > t}
  \smallpad\le
  \exp\Big(1-\frac{\Cexpbsgn t^2}{\Cmaxsgn\sgnz\Mra\sqrt{\log\px+1}}\Big)
\end{align*}
for~$t\ge0$, where~$\Cexpbsgn$ is an absolute constant.
Combine the above result with that of three displays previous, choose~$t=\n^{\zeta}/\Mrt$ and cite the growth conditions of the present lemma to conclude that
\benn
  \Prb\bset{|\ip{\rt_j,\dd_i}|>\n^{\zeta}}
  \smallpad\lesssim
  \exp\Big(-\frac{\n^{2\zeta}}{\sqrt{\log\px}}\Big)
  \smallpad{\overset{\n\to\infty}{\to}}0\,,
\eenn
as required for Condition~\ref{wl:noise2a}.
\end{proof}
}

\noindent The proof of Lemma~\ref{satisfied} is found in Section~\ref{proofs:normality} of the Appendix.
The requirements that the quantities~$\subgauss,\|\rSz\|_\infty,$ and~$\Mrt$ in the statement of Lemma~\ref{satisfied} be of constant order can be relaxed at the cost of introducing more complex growth conditions.
% The growth condition~$\exp(-\n^{2\zeta}/\sqrt{\log\px})=o(\n^{-1})$
Condition~(v) of~Lemma~\ref{satisfied}
is satisfied under reasonable constraints on the growth of~$\px$ --- for instance, if~$\sqrt{\log\px}/\n^{2\zeta}\lesssim1$.

%% file: estimation/estimation.tex
\section{Two-Stage Lasso}\label{estimation}
\myappendix{
\section{Materials required for Section~\ref{estimation}}\label{proofs:estimation}
}

Theorem~\ref{wl} depends on high-level assumptions that ensure good behavior of the remainder terms~$\rem_\ell$ and standard error estimate~$\wjseh$.
In this section, we demonstrate how such conditions may be satisfied in the high-dimensional setting.
In particular, we introduce in Section~\ref{ss:two-stage} a two-stage Lasso estimation procedure, for which
we provide theoretical bounds in Section~\ref{bounds}.
The rates for the second-stage estimation error are particularly involved due to the dependence on the predicted conditional means from the first-stage estimation.
In Section~\ref{ss:remainders}, we identify conditions under which the remainder terms~$\rem_\ell$ vanish in probability under the two-stage Lasso procedure.
We also propose estimators of the standard errors and provide conditions under which these estimators converge to the true standard errors. 

%% file: estimation/two-stage.tex
\subsection{Two-stage estimator}\label{ss:two-stage}

For $j\in[\px]$, we let~$\rah^j$ denote the \emph{first-stage Lasso estimator}
\be\label{def:rah}
  \rah^j \relates{\in}
  \argmin_{\bs{a}\in\R^{\pz}} \big\{
    \|\dx^j - \dZ\bs{a}\|_2^2/(2n) + \tunej \|\bs{a}\|_1
  \big\} \,.
\ee
% The quantity $\tunea$ is not a tuning parameter per se: it does not appear in the objective function of any Lasso program of the present essay.
% However, it occupies the same role in the estimation error bounds of the present section as does a Lasso tuning parameter.
% Thus, we preserve the tuning parameter notation for the quantity $\tunea$.
We let~$\ddijh := \dz_i^\top\rah^j$ denote the predicted conditional mean of~$\dxij$ given~$\dz_i$ based on the estimates~$\rah^j$ and write~$\dDh= \dZ\rAh$, where the matrix~$\dDh$ has columns given by~$\ddh^j := (\ddch_{1j},\ldots,\ddch_{\n j})^\top$ and the matrix~$\rAh\in\R^{\pz\times\px}$ has columns given by the~$\rah^j$.

We define the \emph{second-stage Lasso estimator} to be
\be\label{prog:ss}
  \rbla \relates{\in}
  \argmin_{\bs{b}\in\R^{\px}} \big\{
    \|\dy - \dDh\bs{b}\|_2^2/(2n) + \tuneb\|\bs{b}\|_1
  \big\} \,.
\ee
In Section~\ref{second}, we develop sparsity-based results that require the following quantities.
We write~$\Saj=\supp\ra^j$ for the~\emph{active sets} of the first-stage regression parameters, and we write~$\saj:=|\Saj|$ and~$\sa:=\maxjpx\saj$; we write~$\Sb:=\supp\rb$ for the active set of the second-stage regression parameter and~$\sb:=|\Sb|$.
We note that~$\ell_0$ sparsity is not a limitation in principle and that more general regression vectors may be considered at the price of additional complexity \cite[Sections~6.2.3-4]{BvdG11}, \cite{BCCH12}.

\begin{rek}[Model identifiability]
Let~$\ra^j_{\Saj} := \{\rajk \,:\, \rajk \neq 0\}$ denote the restriction of~$\ra^j$ to~$\Saj$,
and let~$S=\cup_j\Saj$.
Write~$\rA_S$ for the~$|S|\times\px$ matrix with columns~$\ra^j_{\Saj}$ and~$\dZ_S$ for the $\n\times|S|$ matrix with columns (re-indexed as necessary) corresponding to the columns of~$\rA_S$.
Note that~$\dD = \dZ_S\rA_S$ and hence that~$\rSd=\cov(\dD)=\rA_S^\top\cov({\dZ_S})\rA_S$.
We require that~$\rSd$ be invertible and hence of full rank
--- for model identifiability, probabilistic gurantees we discuss in the following sections, and because we work directly with the inverse covariance matrix~$\rT$.
Thus we must have ~$\px=\rank\rSd \le \min(\rank\rA_S, \rank\cov(\dZ_S))\le|S|$.
This is a relatively strong assumption, which we implicitly require throughout the sequel.
However, we emphasize that this requirement is an artifact of the $\ell_0$-sparsity-based methods by which we derive the bounds for the first- and second- stage estimation errors.
As noted above, we chose these methods to simplify our exposition and that one can obtain morally similar but more complex bounds even when the~$\ra^j$ are not sparse~\cite[Sections~6.2.3-4]{BvdG11}.
Thus, in general, this restriction is not impractical.
\end{rek}

%% file: estimation/bounds/bounds.tex
\subsection{Estimation error bounds}\label{bounds}
\myappendix{
\subsection{Materials required for Section \ref{bounds}}\label{proofs:bounds}
}

In this section, we present estimation error bounds for the first- and second- stage estimators described in Section~\ref{ss:two-stage}.
Both such bounds depend on the same fundamental strategy for proving finite-sample guarantees for $\ell_1$-regularized estimators.
This strategy consists of two parts.
The first part is the \emph{oracle inequality}, which establishes a deterministic bound for the estimation and prediction performance of the Lasso on a particular set of interest.
The second part is the control of the \emph{empirical process term}, which defines the set of interest.
We include such prerequisites in Section \ref{proofs:bounds} of the Appendix.

\new{Before we present the estimation error bounds for the first- and second-stage Lasso estimators, we first discuss  the compatibility condition, which is required in the proof of the oracle inequality.}

\subsubsection{Compatibility condition}\label{compatibility}

The oracle bounds rely on the good behavior of certain moduli of continuity of the empirical Gram matrices~$\rSzh=\dZ^\top\dZ/n$ and~$\rSdh=\dDh^\top\dDh/n$.
We codify this requirement in the following definition.

\begin{defi}[Compatibility condition]\label{def:cc}
% For a given index set $\idx\subseteq[\px]$ and controlled quantities $\ccc>0$ and $\Mdiff>0$, define the set
% \be\label{cone}
% 	\cone(\idx,\ccc,\Midx) \smallpad{:=}
% 	\bset{\diff\in\R\setminus\0 : \|\diff\|_1\le\Mdiff, \|\diff_{\idx^\comp}\|_1\le\ccc\|\diff_{\idx}\|_1}\,.
% \ee
For a given index set~$\idx\subseteq[p]$,~$p\in\N$, define the double-cone
\be\label{cone}
	\cone(\idx) \smallpad{:=}
	\bset{\diff\in\R^{p}\setminus\0 : \|\diff_{\idx^\comp}\|_1\le3\|\diff_{\idx}\|_1}\,.
\ee
We say that the \emph{compatibility condition} holds for the matrix~$\aM \in \R^{\n \times p}$ relative to the index set~$\idx$ and the constant~$\ccomp^2>0$ defined as
\be\label{eq:ccomp}
	\ccomp^2 \equals
	\inf_{\diff\in\cone(\idx)} \frac{|\idx|\|\aM\diff\|_2^2}{\n\|\diff_{\idx}\|_1^2}
\ee
if the latter is greater than zero.
We call the quantity~$\ccomp^2$ the~\emph{compatibility constant}.
\end{defi}

\noindent \new{The compatibility condition is so named because it interfaces between the~$\ell_1$ norm of the estimation error and the~$\ell_2$ prediction error of the Lasso estimator.
It is a standard tool in the $\ell_1$-regularized estimation literature to ensure identifiability:
it limits the correlations among the predictors such that the estimator can discriminate between the ``relevant'' parameters with index in~$\idx$ and the remaining parameters.
For this purpose, the index set~$\idx$ is taken to be the active set of the target regression parameter~\cite[Chapter 6]{BvdG11}.

The set $\cone(\idx)$ increases with increasing active set~$\idx$;
hence, the larger~$\idx$, the more restrictive the compatibility condition becomes.
This means that in such theories, identifiability and sparsity are closely intertwined.}

The constant~3 is arbitrary.
Alternative choices require adjustment of other constants that appear in the bounds \cite{BRT09}.
A related, slightly stronger condition known as the \emph{restricted eigenvalue condition} is elsewhere used for the same end \cite{BRT09, BvdG11, vdGB09}.
The only known task where such conditions can be avoided is prediction~\cite{Dalalyan17,hebiri2013correlations,lederer2019,vandegeer,doi:10.1002/sta4.186}.

The compatibility and restricted eigenvalue conditions are sometimes defined more generally in terms of the cardinality $s$ of the index set $S$ rather than a specific index set.
For instance,~\cite[Assumption~RE$(s,c_0)$]{BRT09} require for their restricted eigenvalue condition that the quantity
\be\label{re}
	\kappa(s,c_0) \smallpad{:=}
	\min_{\substack{S\subseteq[p]:|S|\le s}}
	\min_{\substack{\diff\neq\0:\\\|\diff_{S^\comp}\|_1\le c_0\|\diff_S\|_1}}
	\|\aM\diff\|_2/(\sqrt{\n}\|\diff_{S}\|_2)
\ee
be bounded away from zero.
The rationale for taking the minimum over all such index sets~$S$ is that the true support of~$\rb$ is unknown.
See also the discussion of~\cite{RZ12}.
We note also that the compatibility and restricted eigenvalue conditions can be replaced by slightly weaker assumptions at the cost of more involved definitions \cite{Dalalyan17}.
% Note also that most expressions of the compatibility or restricted eigenvalue conditions consider the infimum over $\diff\in\R^\px\setminus\0$ that satisfy the ``cone condition''~$\|\diff_{\idx^\comp}\|_1\le\ccc\|\diff_{\idx}\|_1$ but may be unbounded in $\ell_1$ norm.
% We discuss the condition~$\|\diff\|_1\le\Mdiff$ of \eqref{cone} in Section~\ref{second}.

% Suppose that the IV design matrix~$\dZ$ and predicted conditional mean matrix~$\dDh$ satisfy
% (i)~$\dDh$ and (ii)~$\dZ\in\Ccomp(\Saj,\ccompj^2)$.
% Define~$\sa:=\maxjpx\|\Saj\|_0$ and~$\ccompa:=\min\ccompj$.
% For the sake of brevity, we introduce the notation
% \be\label{def:scc}
% \begin{split}
% 	\sccb \smallpad{:=} \sb/\ccompb^2 \,,
% \end{split}
% \quad
% \begin{split}
% 	\scca \smallpad{:=} \sa/\ccompa^2 \,,
% \end{split}
% \ee
% as the right-hand quantities appear frequently in the following results.
% Note that~$\scca \ge \saj/\ccompj$ for each~$j\in[\px]$.

%\subsubsection{Tuning parameters}\label{tuningparameters}

Practical use of the first- and second-stage Lasso estimators requires selection of the tuning parameter~$\tune_j$ and $\tuneb$ for $j\in[\px]$.
A number of proposals for theoretical choices of tuning parameters \cite{BC13, BCH11, BRT09, BvdG11} for the ordinary (one-stage) linear model exist in the regularized regression literature. The bounds of Section~\ref{second} are based on oracle choices of the tuning parameters~$\tune_j,\tuneb$, which depend explicitly on inestimable quantities.
It would be preferable to give results for~\emph{data-adaptive} tuning parameters such as~\cite{Chetelat17,CLW16} for which the respective Lasso problems can feasibly be implemented.
For the present work, we are content to demonstrate that there exist sequences of oracle tuning parameters that tend to zero sufficiently fast to ensure that the remainder terms~$\rem_\ell$ are asymptotically negligible.
In practice, cross-validated choices of Lasso tuning parameters and~$\tf$ chosen according to the scheme described in Section~\ref{ss:inverse} suffice in favorable parameter configurations.
We provide evidence for this claim in Section~\ref{s:numerics}.

%%%%%%%%%%%%%%%%%%%%%%%%%%%%%%%%%%%%%%%%%%%%%%%%%%%%%%%%%%%%%%%%%%%%%%%%%%%%%%%
%%%%%%%%%%%%%%%%%%%%%%%%%%%%%%%%%%%%%%%%%%%%%%%%%%%%%%%%%%%%%%%%%%%%%%%%%%%%%%%
% ORACLE INEQUALITY
\newcommand{\gX}{\bs{W}}
\newcommand{\gx}{\bs{w}}
\newcommand{\gy}{\bs{g}}
\newcommand{\greg}{\bs{\gamma}}
\newcommand{\gregh}{\hat{\greg}}
\newcommand{\gerr}{\bs{h}}
\newcommand{\gS}{S_{\greg}}
\newcommand{\gSc}{\gS^{\comp}}
\newcommand{\gens}{s_{\greg}}
\newcommand{\gcc}{\phi_{\greg}}
\newcommand{\gT}{\mathcal{T}}

\myappendix{%
Our guarantees for estimating $\rA$ and $\rb$ consist of two parts.
The first is the \emph{oracle inequality}, which bounds the $\ell_1$ estimation error of a generic Lasso estimator conditional on the occurrence of a special set $\T$.
The oracle inequality is a fixture of the $\ell_1$ regularized estimation literature; see for instance \cite[Chapter 6]{BvdG11}.
We present it for the sake of completeness.

The oracle inequality itself is specific to neither the first- nor the second- stage estimators of the present work.
Indeed, we require the result to derive bounds for both estimators.
As such, we present the theorem in terms of a generic model that shares notation with neither the first- nor second- stage models described in Section~\ref{model} except for the number of observations $\n$.

\begin{thm}[Oracle inequality]\label{thm:oracle}
Consider the generic linear model
\benn
  \gy \equals \gX\greg + \gerr \,,
\eenn
where~$\gy\in\R^{\n}$ is a vector of univariate responses,~$\gX\in\R^{\n\times p}$ is a design matrix with rows~$\gx_i$, $\greg\in\R^{\n}$ is a noise vector with arbitrary distribution.
Let~$\gregh$ denote the Lasso estimator given by
\benn
  \gregh \smallpad\in
  \argmin_{\bs{a}\in\R^{p}} \big\{
    \|\gy - \gX\bs{a}\|_2^2/(2\n) + \tune \|\bs{a}\|_1
  \big\} \,.
\eenn
Let~$\gS := \supp\greg$, and let $\gens := |\gS|$.
Suppose that~$\gX$ satisfies the compatibility Condition with respect to the index set~$\gS$ and compatibility constant~$\gcc>0$.
Then, on the set
\be\label{soi}
  \gT(\tune) \smallpad{:=} \bset{4\|\gX^\top\gerr/n\|_\infty \le \tune} \,,
\ee
the bound
\benn
  \|\gX(\gregh-\greg)\|_2^2/n + \tune\|\gregh-\greg\|_1
  \smallpad{\leq}
  4\gens\tune^2/\gcc^2
\eenn
holds.
\end{thm}
}%

%%%%%%%%%%
% PROOF
\myappendix{%
\begin{proof}[Proof of Theorem \ref{thm:oracle}]
The proof is algebra. See Theorem 6.1 of \cite{BvdG11}.
\end{proof}
}%

%%%%%%%%%%%%%%%%%%%%%%%%%%%%%%%%%%%%%%%%%%%%%%%%%%%%%%%%%%%%%%%%%%%%%%%%%%%%%%%
% EMPIRICAL PROCESS TERMS

\myappendix{%
\noindent Theorem \ref{thm:oracle} provides a deterministic guarantee for the $\ell_1$ estimation error of a generic Lasso estimator $\gregh$ on the set $\gT(\tune) = \bset{4\|\gX^\top\gerr/n\|_\infty \le \tune}$.
Consequently, it holds that
\begin{align*}
  \Prb\bset{\|\gregh-\greg\|_1 > 4\gens\tune/\gcc^2}
  \smallpad\le
  \Prb\,\gT(\tune)^\comp \,.
\end{align*}
The quantity $\|\gX^\top\gerr/n\|_\infty$ is sometimes called the \emph{empirical process term}; for instance, \cite[Chapter 6]{BvdG11}.
Thus, upper bounds for~$\Prb\,\gT(\tune)^\comp$ yield probabilistic guarantees for the~$\ell_1$ estimation error.
We provide such bounds in Section~\ref{technical} of the Supplementary Materials.
}%

%% file: estimation/bounds/second.tex
\subsubsection{Estimation error bounds}\label{second}
\myappendix{
\subsection{Materials required for Section~\ref{second}}\label{proofs:second}
}

Simultaneous control of the first-stage estimation errors~$\rah^j-\ra^j$ is a straightforward consequence of the standard theoretical results for the Lasso.
We include these bounds in Section~\ref{proofs:second} of the Supplementary Materials.
On the other hand, the bounds for the second-stage estimation error~$\rbla-\rb$, which we study in the present section, are more involved due to the dependence of~$\rbla$ on the predicted conditional means~$\ddh_i$.
Our strategy is to write~$\dy \equals \dDh \rb + \newu$, where~$\newu := \du + [(\dD-\dDh) + \dV]\rb$ and apply concentration results to bound the probability of the event~$\{4\|\dDh^\top\newu\|_\infty\le\tuneb\}$, allowing us to adapt oracle inequality arguments for the Lasso to the present case.
%Below is Appendix B.2
\myappendix{
Recall that the model for~$\dx^j$ is given by
\benn
	\dx^j \equals \dZ\ra^j + \dv^j\,,
\eenn
where~$\dv^j$ has nontrivial covariance with the noise~$\du$.
% Theorem \ref{thm:oracle} entails that, on the set
% \benn
% 	\Tvj(\tunej) \equals \set{\maxzn{\dv^j} \le \tunej} \,,
% \eenn
% it holds that
% \benn
% 	\|\rah^j-\ra^j\|_1 \smallpad\le 4\saj\tunej/\ccompj \,.
% \eenn
It suffices for our purposes to take a na\"ive approach to bounding the quantity~$\maxerrverbfs \equals \maxjpx\|\rah^j-\ra^j\|_1$.
That is, we simultaneously bound the estimation error of each individual task.
One could use a more complex approach such as~\cite{LWZ15} to treat different patterns of joint sparsity amongst the first-stage regression vectors. We make the following assumption.

\noindent The following lemma provides finite-sample guarantees for~$\maxerrfs$ under the choice of tuning parameters in Definition~\ref{tpss}.
}
% \begin{defi}[Generic tuning parameters, $\rAh$]\label{assu:tpgeneric}\
% \begin{enumerate}
% 	\item $\tunea = \maxjpx \tunej$, where $\tunej > 0$ for $j\in[\px]$ are arbitrary;
% \end{enumerate}
% \end{defi}
%
% \begin{defi}[Tuning parameters, $\rAh$, Gaussian noise]\label{assu:tpgaussian}\
% \begin{enumerate}
% 	\item $\tunej = \sdvj c \sqrt{\maxzzn\log\pz/\n}$ for each $j\in[\px]$;
% \end{enumerate}
% \end{defi}

% \begin{rek}[Tuning parameters]
% [Insert remark about possible redundancy between $\tunea$ and $\tuneV$.]
% \end{rek}

%%%%%%%%%%%%%%%%%%%%%%%%%%%%%%%%%%%%%%%%%%%%%%%%%%%%%%%%%%%%%%%%%%%%%%%%%%%%%%%
%%%%%%%%%%%%%%%%%%%%%%%%%%%%%%%%%%%%%%%%%%%%%%%%%%%%%%%%%%%%%%%%%%%%%%%%%%%%%%%

\myappendix{%
The following generic bound for~$\maxerrfs$ can be combined with concentration results for specific distributions of the first-stage noise elements.
We present it separately for the sake of modularity with respect to such assumptions.

\begin{lem}[Generic bound for~$\maxerrfs$]\label{lem:eefs} % Estimation Error, First-Stage
Suppose that Assumption~\ref{assu:ccss} holds.
For each~$j\in[\px]$, let the sets~$\Tvj(\lambda)$ be as defined in Lemma~\ref{lem:epV}.
% \be\label{def:Tvj}
% 	\Tvj(\lambda)
% 	\smallpad{:=}
% 	\big\{4\|\dZ^\top\dv^j/n\|_\infty \leq \lambda\big\}\,.
% \ee
It then holds on the set~$\bigcap_{j\in[\px]} \Tvj(\tunej)$ that
\benn
	\maxerrfs
	\smallpad\le
	4\sa\tunea/\ccompa^2
	\equals
	4\scca\tunea \,,
\eenn
where~$\tunej$ is the tuning parameter for the respective first-stage Lasso problem,~$\tunea$ and~$\sa$ are as defined in Section~\ref{s:two-stage} and~$\ccompa$ is as defined in Assumption~\ref{assu:ccss}.
\end{lem}
}%

%%%%%%%%%%
% PROOF
\myappendix{%
\begin{proof}[Proof of Lemma \ref{lem:eefs}]
For each~$j\in[\px]$, \cite[Theorem 6.1, Lemma 6.2]{BvdG11} entails that
\begin{align*}
	\Tvj(\tunej)
	\equals
	\big\{4\|\dZ^\top\dv^j/n\|_\infty \leq \tunej\big\}
	\smallpad{&\subseteq}
	\big\{\|\rah^j-\ra^j\|_1 \leq 4\saj\tunej/\ccompj^2\big\} \\
	\smallpad{&\subseteq}
		\big\{\|\rah^j-\ra^j\|_1 \leq 4\sa\tunea/\ccompa^2\big\}\,,
\end{align*}
where the latter containment follows by specification of~$\sa, \tunea,$ and~$\ccompa$.
Take intersections over~$j\in[\px]$ on both sides of the above display to conclude.
\end{proof}
}%

%%%%%%%%%%%%%%%%%%%%%%%%%%%%%%%%%%%%%%%%%%%%%%%%%%%%%%%%%%%%%%%%%%%%%%%%%%%%%%%
%%%%%%%%%%%%%%%%%%%%%%%%%%%%%%%%%%%%%%%%%%%%%%%%%%%%%%%%%%%%%%%%%%%%%%%%%%%%%%%

\myappendix{
\begin{lem}[Bound for $\maxerrfs$]\label{lem:eefsg}
Suppose that Assumptions~\ref{assu:ccss} and~\ref{assu:subgaussian} hold.
Set $\tunej$ according to Definition~\ref{tpss}.
Then,
\benn
	\Prb\bset{\maxerrfs > 4\scca\tuneVverb}
	\smallpad\le
	\prbepallvj\,,
\eenn
where~$\tolV$ is as specified in Definition~\ref{tpss} and~$\Choeff$ is as defined in Lemma~\ref{lem:epg}.
% Consequently, if $\sa = o(\sqrt{n/\log\pz})$, then $\maxerrfs = o_P(1)$.
\end{lem}
}

%%%%%%%%%%
% PROOF
\myappendix{
\begin{proof}[Proof of Lemma \ref{lem:eefsg}]
Lemma \ref{lem:eefs} entails that
\benn
	\Prb\{\maxerrfs > 4\sa\tunea/\ccompa^2\}
	\smallpad{\le}
	\Prb \Big\{\bigcap_{j\in[\px]} \Tvj(\tunej)\Big\}^\comp \,.
\eenn
% \begin{align*}
	% \Prb\{\maxerrfs > 4\sa\tunea/\ccompa^2\}
	% \smallpad{&\le}
	% \Prb \Big\{\bigcap_{j\in[\px]} \Tvj(\tunej)\Big\}^\comp \\
	% \smallpad{&\le}
	% 	\sum_{j\in[\px]}\Prb\,\Tvj(\tunej)^\comp \,.
	% \smallpad{&\le}
		% \px\Prb\bset{\sdvj c \sqrt{\maxzzn\log\pz/\n}} \,.
% \end{align*}

Apply the estimate of Lemma~\ref{lem:epV} for the right-hand side to conclude.
\end{proof}

%%%%%%%%%%%%%%%%%%%%%%%%%%%%%%%%%%%%%%%%%%%%%%%%%%%%%%%%%%%%%%%%%%%%%%%%%%%%%%%
%%%%%%%%%%%%%%%%%%%%%%%%%%%%%%%%%%%%%%%%%%%%%%%%%%%%%%%%%%%%%%%%%%%%%%%%%%%%%%%

% \begin{rek}[Choice of $c$]\label{remarkc}
% In order for Lemma~\ref{lem:eefsg} to entail asymptotic negligibility of~$\maxerrfs$, the controlled quantity~$c$ must be chosen so that~$\cepg>1$.
% In the sequel, we assume without further comment that such a choice has been made.
% \end{rek}

\noindent
Note that, in order to obtain a rate of convergence, we must choose each~$\tolj$ so that the quantity~$\min_{j\in[\px]}\tolj^2\Choeff/\sgnvj^2 - 2$ is bounded strictly away from zero.
Such a task may not be feasible in practice.
Our empirical results, which we present in Section~\ref{s:numerics}, suggest that cross-validated choices of the tuning parameters for the first- and second-stage Lasso estimators suffice for good behavior of the resultant updated estimator~$\rbdb$.
For the sake of our theory, we assume that such appropriate choices of~$\tolj$ have been made.
Given such an assumption, Lemma~\ref{lem:eefsg} entails that~$\maxerrfs=O_{\Prb}(\sa\sqrt{\log(\pz)/\n})$, essentially identical to the Lasso rate for single task regression problems.
}

We require the following assumption for the bounds of this section.

\begin{assu}[Compatibility conditions]\label{assu:ccss}

\begin{enumerate*}[label=(\arabic{*})]
	\item
		There exists for each active set~$\Saj$ of the first-stage model a constant~$\ccompj^2 > 0$ such that~$\dZ$ satisfies the compatibility condition with respect to~$\Saj$ and~$\ccompj^2$.
		We assume that the~$\n$-indexed sequences of such constants~$\ccompj^2$ are bounded strictly away from zero uniformly in~$\n$.
		We write~$\ccompa^2:=\maxjpx\ccompj^2$.
	\item
		There exists a constant~$\ccompb$ such that~$\dDh$ satisfies the compatibility condition with respect to~$\Sb$ and~$\ccompb$.
		We assume that the~$\n$-indexed sequence of such constants is bounded strictly away from zero uniformly in~$\n$.
\end{enumerate*}
\end{assu}

\noindent \new{\noindent Assumption \ref{assu:ccss} imposes the compatibility condition on the random object $\widehat{\bm{D}}$.
In Lemma \ref{ccssg} below, we provide sufficient conditions under which Assumption~\ref{assu:ccss} holds with high probability.}

Note that whether~$\dZ$ satisfies the compatibility condition with respect to one active set~$\Sa_{j_1}$ does not bear directly on whether it satisfies the compatibility condition with respect to another active set~$\Sa_{j_2}$ for~$j_1,j_2\in[\px]$.
As such, it is non-trivial to assume that the compatibility condition as specified in Definition~\ref{def:cc} holds for each active set~$\Saj$ for~$j\in[\px]$ when~$\px$ tends to infinity.
However, the condition that~$\dZ$ satisfies the compatibility condition with respect to each active set~$\Saj$ is entailed by requiring that~$\kappa(s, c_0)$ of~\eqref{re} for~$s=\max_{j\in[\px]\saj}$ and~$c_0=3$ be bounded away from 0.
Thus, the need to accommodate multiple active sets does not thereby significantly alter the treatment of the compatibility condition.

\myappendix{
The following bound is required for Lemma~\ref{lem:eessg}

\begin{lem}[Control of $\|\dDh^\top\newu/\n\|_\infty$]\label{lem:ep}
Let $\newu=\du + [(\dD-\dDh) + \dV]\rb$.
Then,
\begin{align*}
	\|\dDh^\top\newu/\n\|_\infty
	\smallpad{&\le}
	\maxerrfs\maxzzn(\maxerrfs\|\rb\|_1+\|\rA\|_{L_1})\\
	&\quad
		+ (\maxerrfs+\|\rA\|_{L_1})(\maxzn{\dV}\|\rb\|_1+\maxzn{\du})\,.
\end{align*}
\end{lem}

%%%%%%%%%
% PROOF
\begin{proof}[Proof of Lemma \ref{lem:ep}]
Write $\dDh^\top = (\dDh - \dD)^\top + \dD^\top$
% and
% \benn
% 	\newu \equals (\dD - \dDh)\rb + \dV\rb  + \du
% \eenn
to find that
\begin{align}\label{ep:pledge}
	\|\dDh^\top\newu/n\|_\infty
	\aequals
	\big\| \big[(\dDh - \dD)^\top + \dD^\top\big] \big[(\dD - \dDh)\rb + (\dV\rb
		+ \du)\big]/n\big\|_\infty \nonumber\\
	\smallpad{&\le}
		\|(\dDh - \dD)^\top(\dDh - \dD)\rb/n\|_\infty + \| (\dDh-\dD)^\top(\dV\rb
			+ \du)/n\|_\infty  \nonumber\\
	& \qquad
		+ \|\dD^\top (\dDh - \dD)\rb/n \|_\infty + \| \dD^\top (\dV\rb
			+ \du)/n \|_\infty \nonumber\\
	\smallpad{&:=}
		\I_1 + \I_2 + \I_3 + \I_4
\end{align}
We treat each quantity in the right-hand side above in turn.

For $\I_1$, write
\begin{align*}
	\I_1 \smallpad{=} \|(\dDh-\dD)^\top(\dDh-\dD)\rb/n\|_\infty
		\smallpad{&\le} \|(\dDh-\dD)^\top(\dDh-\dD)/n\|_\infty\|\rb\|_1 \,.
\end{align*}
Recall that $\dDh-\dD=\dZ(\rAh-\rA)$ and write
\begin{align*}
	\|(\dDh-\dD)^\top(\dDh-\dD)/n\|_\infty
	\smallpad{&=} \|(\rAh-\rA)^\top\rSzh(\rAh-\rA)\|_\infty \\
	\smallpad{&\le} \maxerrverbfs^2\|\rSzh\|_\infty
	\smallpad= \maxerrfs^2\maxzzn  \,,
\end{align*}
where the second line follows from repeated application of H\"older's inequality.
Combine the two previous displays to conclude that
\be\label{turn1}
	\I_1 \smallpad\le \maxerrfs^2\maxzzn\|\rb\|_1\,.
\ee

\providecommand{\It}[1]{\I_{2,\mathrm{#1}}}
For $\I_2$, write
\begin{align*}
	\I_2 \smallpad{&=} \|(\dDh-\dD)^\top(\dV\rb + \du)/n\|_\infty\\
	\smallpad{&=} \|(\rAh-\rA)^\top\dZ^\top(\dV\rb + \du)/n\|_\infty\\
	\smallpad{&\le}
		\underbrace{\|(\rAh-\rA)^\top\dZ^\top\dV\rb/n\|_\infty}_{\It{a}}
	+ \underbrace{\|(\rAh-\rA)^\top\dZ^\top\du/n\|_\infty}_{\It{b}}\,.
\end{align*}
Applications of H\"older's inequality yield
\begin{align*}
	\It{a} \smallpad{\le} \maxerrverbfs \|\dZ^\top\dV/n\|_\infty \|\rb\|_1
\end{align*}
and
\begin{align*}
	\It{b} \smallpad{\le} \maxerrverbfs\|\dZ^\top\du/n\|_\infty
	\smallpad{=} \maxerrfs\maxzn{\du}\,.
\end{align*}
Combine the previous three displays to conclude that
\be\label{turn2}
	\I_2 \smallpad\le \maxerrfs(\maxzn{\dV}\|\rb\|_1 + \maxzn{\du}) \,,
\ee

For the quantity $\I_3$, write
\benn
	\I_3 \equals \|\dD^\top(\dDh - \dD)\rb/n\|_\infty
	\smallpad{\le} \|\dD^\top(\dDh-\dD)/n\|_\infty\|\rb\|_1
\eenn
and observe that
\begin{align*}
	\dD^\top(\dDh-\dD)/n
	\equals \rA^\top\dZ^\top\dZ(\rAh-\rA)/n
	\equals \rA^\top\rSzh(\rAh-\rA) \,,
\end{align*}
which yields
\begin{align*}
		\|\dD^\top(\dDh-\dD)/n\|_\infty
		\aequals
			\|\rA^\top\rSzh(\rAh-\rA)\|_\infty\\
		\smallpad{&\le}
			\|\rA\|_{L_1}\|\rSzh\|_\infty\|(\rAh-\rA)\|_{L_1}
\end{align*}
after repeated application of H\"older's inequality.
Conclude from the previous three displays that
\be\label{turn3}
	\I_3
	\smallpad{\le} \maxerrfs\maxzzn\|\rA\|_{L_1} \,.
\ee

For the quantity $\I_4$, write
\begin{align}\label{turn4}
	\I_4
	\aequals
	\|\dD^\top(\dV\rb + \du)/n \|_\infty \nonumber\\
	\aequals
		\|\rA^\top\dZ^\top(\dV\rb + \du)/n \|_\infty \nonumber\\
	\smallpad{&\le}
		\|\rA\|_{L_1}
			\big(
				\|\dZ^\top\dV\rb/n\|_\infty
				+
				\|\dZ^\top\du/n\|_\infty
			\big) \nonumber\\
	\smallpad{&\le}
		\big(
			\maxzn{\dV}\|\rb\|_1 + \maxzn{\du}
		\big)
		\|\rA\|_{L_1} \,.
\end{align}
The original claim follows from line \eqref{ep:pledge} and lines \eqref{turn1}-\eqref{turn4}.
\end{proof}
}%

We will refer to the following choices of tuning parameters throughout the sequel.

\begin{defi}[Tuning parameters]\label{tpss}
(1) For the first-stage Lasso estimator, set~$\tunej:=\tunejverb$,
% \benn
% 	\tunej\smallpad{:=}\tunejverb \,,
% \eenn
where~$\tolj>0$ is a controlled quantity.
We let~$\tunefs := (\tune_1,\ldots,\tune_{\px})$ denote the tuple of first-stage tuning parameters, and we write~$\tunea \smallpad{:=} \maxjpx \tunej =\tuneVverb$,
% \benn
% 	\tunea \smallpad{:=} \maxjpx \tunej \equals \tuneVverb \,,
% \eenn
where~$\tolV=\maxjpx\tolj$.
(2) For the second-stage Lasso estimator, set
\benn
  \tuneb \smallpad{=}
  16\scca\tunea
		\maxzzn \big(4\Mb\scca\tunea+\Mra\big)
		+ \big(4\scca\tunea+\Mra\big)
		\big(\Mb\tuneV+\tuneu\big)\,,
\eenn
where~$\Mra,\Mb$ are as defined in Assumption~\ref{regreg}
% ~$\tuneV := \tolV(\maxzzn \cdot (\log\pz)/\n)^{1/2}$, $\tuneu = \tuneuverb$, and~$\tolu>0$ is a controlled quantity.
and
\benn
\begin{split}
  \tuneV \smallpad{:=} \tuneVverb \,,
\end{split}
\qquad
\begin{split}
  \tuneu \smallpad{:=} \tuneuverb \,,
\end{split}
\eenn
where~$\tolu>0$ is a controlled quantity.
\end{defi}

\noindent Note that, unlike as in much of the related literature, the tuning parameters we identify above are random, in particular due to the term $\maxzzn$.
Our theory handles the consequences of this allowance in Condition~\ref{gr:zzn} of Assumption~\ref{gr}, for which we provide subsequent justification.
As we note above, our practical choice of tuning parameters is guided by cross-validation.

We now present probabilistic bounds for the~$\ell_1$ estimation error for the second-stage Lasso estimator~$\rbla$.

\myappendix{%
\noindent The following generic bound for~$\|\rbla-\rb\|_1$ can be combined with concentration results for specific distributions of the first- and second-stage noise elements.
We present it separately for the sake of modularity with respect to such assumptions.

\begin{lem}[Generic bound for $\|\rbla-\rb\|_1$]\label{lem:eess}
Suppose that Assumption~\ref{assu:ccss} holds.
Let $\tuneV, \tuneu > 0$ be arbitrary.
Set $\tuneb$ according to Definition~\ref{tpss}.
Then, on the set $\TV(\tunefs)\cap\TV(\tuneV)\cap\Tu(\tuneu)$, where~$\TV$ and~$\Tu$ are defined as in Lemmas~\ref{lem:epu} and~\ref{lem:epV}, respectively, and~$\tunefs$ is the tuple of first-stage tuning parameters, we have
\benn
	\|\rbla-\rb\|_1
	\smallpad\le
	4\sccb\tuneb \,.
\eenn
\end{lem}
}%

%%%%%%%%%%
% PROOF
\myappendix{%
\begin{proof}[Proof of Lemma \ref{lem:eess}]
By Theorem \ref{thm:oracle}, we have
\benn
	\Tnewu = \big\{4\|\dDh^\top\newu/n\|_\infty \leq \tune\big\}
	\smallpad\subseteq
	\big\{\|\rbla-\rb\|_1\le4\sb\tune/\ccomp^2\big\} \,.
\eenn
It therefore suffices to show that $\TV(\tunefs)\cap\Tu(\tuneu)\subseteq\Tnewu(\tuneb)$ for the present choice of $\tuneb$.
Lemma \ref{lem:ep} gives the bound
% \begin{align*}
% 	\|\dDh^\top\newu/\n\|_\infty
% 	\smallpad{&\le}
% 	\maxerrfs\maxzzn(\maxerrfs\|\rb\|_1+\|\rA\|_{L_1})\\
% 	&\quad
% 		+ (\maxerrfs+\|\rA\|_{L_1})(\maxzn{\dV}\|\rb\|_1+\maxzn{\du})\,.
% \end{align*}
\begin{multline*}
	\|\dDh^\top\newu/\n\|_\infty
	\smallpad{\le}
	\underbrace{\maxerrfs\maxzzn(\Mb\maxerrfs+\Mra)}_{\I_1} \\
		+ \underbrace{(\maxerrfs+\Mra)(\maxzn{\dV}\|\rb\|_1+\maxzn{\du})}_{\I_2}\,.
\end{multline*}
Cite Lemma \ref{lem:eefs} to conclude that, on the set $\TV(\tunefs)$,
\begin{align*}
	\I_1
	\smallpad\le
	4\scca\tunea\maxzzn(4\Mb\scca\tunea+\Mra)\,.
\end{align*}
Note that, on the set $\TV(\tunefs)\cap\TV(\tuneV)\cap\Tu(\tuneu)$,
\benn
	\I_2
	\smallpad\le
	\foq\big(4\scca\tunea+\Mra\big)\big(\Mb\tuneV+\tuneu\big)
\eenn
by specification.
Multiply the two previous displays by 4 and combine with the third-previous display to conclude that~$\TV(\tunefs)\cap\TV(\tuneV)\cap\Tu(\tuneu)\subseteq\Tnewu(\tuneb)$ for the present choice of~$\tuneb$, as required.
\end{proof}
}%

%%%%%%%%%%%%%%%%%%%%%%%%%%%%%%%%%%%%%%%%%%%%%%%%%%%%%%%%%%%%%%%%%%%%%%%%%%%%%%%

\newcommand{\compatibilityprob}{t_n}
\newcommand{\compatibilityset}{\mathcal{T_{\operatorname{cc}}}}
\begin{lem}[Bound for $\|\rbla-\rb\|_1$]\label{lem:eessg}
Suppose that Assumption~\ref{assu:subgaussian} holds and that the compatibility conditions~\ref{assu:ccss} are satisfied with probability at least~$\compatibilityprob=o(1)$.
For each $j\in[\px]$, set $\tunej$ according to Definition~\ref{tpss};
set $\tuneb$ according to Definition~\ref{tpss}.
Then,
\begin{align*}
	\Prb \Big\{& \|\rbla-\rb\|_1 >
		 4\sccb\Big(
			4\scca\tolV\big(\Mb\tolV[16\scca\maxzzn+1]+\tolu\big)\maxzzn(\log\pz)/\n \\
			&\qquad~~~~~~~~~ +
			\Mra\big(16\scca\maxzzn\tolV+\tolu[\Mb+1]\big)(\maxzzn(\log\pz)/\n)^{1/2}
		\Big)\Big\} \\
		& \smallpad\le
		\prbepu + \prbepallvj+\compatibilityprob\,,
\end{align*}
where $\Choeff$ is as specified in Lemma \ref{lem:epg}.
% Then,
% \benn
% 	\Prb \big\{\|\rbla-\rb\|_1 > 4\sccb\tuneb \big\}
% 	\smallpad\le
% 	2\pz^{1-\cepg} + 2\pz^{-\cepg}\,.
% \eenn
\end{lem}

%%%%%%%%%%
% PROOF
\myappendix{%
\begin{proof}[Proof of Lemma \ref{lem:eessg}]
Note first that
\benn
	\TV(\tuneV)
	\equals
	\bigcap_{j=1}^{\px}\bset{4\maxzn{\dv^j}>\tunej}
	\subseteq
	\TV(\tunefs)\,.
\eenn
Let~$\compatibilityset$ be the set on which Assumption~\ref{assu:ccss} holds.
 Lemma~\ref{lem:eess} then entails that
\benn
	\big\{\|\rbla-\rb\|_1 > 4\sccb\tuneb \big\}\cap \compatibilityset
	\smallpad\subseteq
	(\TV(\tunefs)\cap\Tu(\tuneu))^\comp
	\equals
	\TV(\tunefs)^\comp \cup \Tu(\tuneu)^\comp\,.
\eenn
Thus,
\benn
	\Prb\{\|\rbla-\rb\|_1 > 4\sccb\tuneb\big\}
	\smallpad\le
	\Prb\,\TV(\tunefs)^\comp + \Prb\,\Tu(\tuneu)^\comp+\compatibilityprob\,.
\eenn
Now substitute the present choices of tuning parameters and cite the estimates of Lemmas~\ref{lem:epu} and~\ref{lem:epV}.
\end{proof}
}%

\noindent Lemma~\ref{lem:eessg} entails that~$\|\rbla-\rb\|_1 = O_{\Prb}\big(\sb\sa^2(\log\pz)/\n + \sb\sa\sqrt{(\log\pz)/\n}\big)$.
Thus, we see that the convergence rate of the second-stage Lasso estimator is slower than the typical rate of~$\sb\sqrt{\log(\px)/\n}$ in the ordinary (sub-)Gaussian linear model.
In an interesting paper, \cite{Zhu15} provides $L_{1}$ convergence rate of the second stage estimator for unconditional IV models under different sets of assumptions.
If some stronger assumptions such as the ones in Theorem 2.1 of \cite{Zhu15} are satisfied, one can show that the convergence rate in Lemma ~\ref{lem:eessg} can be improved.
Whether it can be improved under the assumptions of the present paper is a direction for future work.

%%%%%%%%%%%%%%%%%%%%%%%%%%%%%%%%%%%%%%%%%%%%%%%%%%%%%%%%%%%%%%%%%%%%%%%%%%%%%%%
%%%%%%%%%%%%%%%%%%%%%%%%%%%%%%%%%%%%%%%%%%%%%%%%%%%%%%%%%%%%%%%%%%%%%%%%%%%%%%%
% SECOND STAGE COMPATIBILITY CONDITION

% \subsubsection{Second-stage compatibility condition}\label{sscc}
% \myappendix{
% \subsection{Materials required for Section~\ref{sscc}}\label{proofs:sscc}
% }

Since Lemma~\ref{lem:eessg} requires that Assumption~\ref{assu:ccss} holds, we need to demonstrate the latter's feasibility.
The following Lemma provides such a guarantee.
For other approaches to studying the empirical compatibility constants and restricted eigenvalues of random matrices, see~\cite{RZ12, vdGM14}.
Unlike the extant literature on the compatibility condition, however, we must account for the prediction error of~$\dDh$.
% The following result is therefore a novel extension of the literature to a situation in which the matrix of interest has structure other than that of a sample covariance matrix.

\myappendix{
\begin{lem}[Second-stage compatibility constant]\label{lem:cc}
% For square matrices~$\aM\in\R^{p\times p}$, let~$\eigmin(\aM)$ denote the minimal eigenvalue of~$\aM$.
Let~$\idx\subseteq[p]$ be an arbitrary index set with~$\sidx=|\idx|$.
For a given matrix $\aM\in\R^{\n\times p}$, define the quantity
\benn
  \ccmin(\aM,\idx)
  \equals
  \inf_{\gd\in\cone(\idx)} \frac{\sidx\|\aM\diff\|_2^2}{\n\|\diff_{\idx}\|_1^2} \,.
\eenn
% Suppose that~$\diff\in\R^{\px}\setminus\0$ satisfies~$\|\diff_{\idx^{\comp}}\|_1 \leq \ccc \|\diff_{\idx}\|_1$, and let $\gep_1,\gep_2>0$ be arbitrary.
Let $\gep_1,\gep_2>0$ be arbitrary. Then,
\begin{align*}
 & \Prb\left\{\ccmin(\dDh,\idx)<\eigmin(\rSd)-\gep_2-\gep_1\right\} \\
  &\qquad\qquad\qquad\le\ \,
    \Prb\bset{16\sidx(2\Mra\maxerrfs + \maxerrfs^2)\maxzzn > \gep_1} \\
		&\qquad\qquad\qquad\qquad + \Prb\bset{16\sidx\|\rSdb-\rSd\|_\infty > \gep_2} \nonumber\,,
\end{align*}
% \begin{align*}
%   &\ \, \Prb\bset{\ccmin(\rSdh, \idx,\ccc)>(1-\gep_1)\big(\eigmin(\rSd)-(1+\ccc)^2\gep_2\big)} \\
%   \le&\ \, \Prb\bset{(2\Mra\maxerrfs + \maxerrfs^2)\maxzzn\|\diff\|_1^2 > \gep_1} + \Prb\bset{\sidx\|\rSdb-\rSd\|_\infty > \gep_2} \,,
% \end{align*}
where~$\rSdb=\dD^\top\dD/\n$.
\end{lem}

%%%%%%%%%%
% PROOF
\begin{proof}[Proof of Lemma \ref{lem:cc}]
Let~$\idx,\sidx$ be as in the statement of Lemma~\ref{lem:cc}, and let~$\diff \in \R^{\px} \setminus \{\bs{0}\}$ satisfying~$\|\diff_{\idx^{\comp}}\|_1 \leq 3 \|\diff_{\idx}\|_1$ be arbitrary.
Write~$\dDh = \dD + (\dDh - \dD)$, so that
\begin{align*}
	\|\dDh\diff\|_2^2 \aequals \|[\dD + (\dDh-\dD)]\diff\|_2^2 \\
	\aequals \big\langle [\dD+(\dDh-\dD)]\diff,[\dD+(\dDh-\dD)]\diff\big\rangle \\
	\aequals \|\dD\diff\|_2^2 + 2\ip{\dD^\top(\dDh-\dD)\diff,\,\diff} + \ip{(\dDh-\dD)^\top(\dDh-\dD)\diff,\diff} \,.
\end{align*}
Thus,
\begin{align}\label{cc:pledge}
	\|\dDh\diff\|_2^2/n
  \smallpad{&\ge}
    \|\dD\diff\|_2^2/n - \underbrace{2|\ip{\dD^\top(\dDh-\dD)\diff/n,\diff}|}_{\I_1} \nonumber\\
  &\quad
    - \underbrace{|\ip{(\dDh-\dD)^\top(\dDh-\dD)\diff/n,\diff}|}_{\I_2} \, .
\end{align}
We now obtain bounds for the quantities~$\I_1,\I_2$.
% Thus,
From repeated applications of H\"older's inequality, write
\begin{align*}
  \I_1 \smallpad{\lesssim}
  |\ip{\dD^\top(\dDh-\dD)\diff/\n,\diff}|
  \smallpad{&\le}
  \|\dD^\top(\dDh-\dD)/\n\|_\infty \|\diff\|_1^2 \\
  \aequals
  \|\rA^\top\dZ^\top\dZ(\rAh-\rA)/\n\|_\infty \|\diff\|_1^2 \\
  \smallpad{&\le}
  \|\rA\|_{L_1}\|\dZ^\top\dZ/\n\|_\infty\|\rAh-\rA\|_{L_1}\|\diff\|_1^2 \\
  \smallpad{&\le}
  \Mra\maxzzn\maxerrfs\|\diff\|_1^2
\end{align*}
and
\begin{align*}
  \I_2 \aequals
  |\ip{(\dDh-\dD)^\top(\dDh-\dD)\diff/\n,\diff}| \\
  \smallpad{&\le}
  \|(\dDh-\dD)^\top(\dDh-\dD)/\n\|_\infty \|\diff\|_1^2 \\
  \aequals
  \|(\rAh-\rA)^\top\dZ^\top\dZ(\rAh-\rA)/\n\|_\infty \|\diff\|_1^2 \\
  \smallpad{&\le}
  \|\rAh-\rA\|_{L_1}\|\dZ^\top\dZ/\n\|_\infty\|\rAh-\rA\|_{L_1}\|\diff\|_1^2 \\
  \smallpad{&\le}
  \maxzzn\maxerrfs^2\|\diff\|_1^2 \,.
\end{align*}
Combine the previous two displays with~\eqref{cc:pledge} to find that
\benn
  \|\dDh\diff\|_2^2/\n \smallpad\ge \|\dD\diff\|_2^2/\n - (2\Mra\maxerrfs + \maxerrfs^2)\maxzzn\|\diff\|_1^2 \,.
\eenn
By assumption, we have~$\|\diff\|_1\le4\|\diff_{\idx}\|_1$.
Substitute this expression in the right-hand side above and multiply through by~$\sidx/\|\diff_{\idx}\|_1^2$ to obtain
\begin{align*}
	\frac{\sidx\|\dDh\diff\|_2^2}{\n\|\diff_{\idx}\|_1^2}
	\smallpad\ge
	\frac{\sidx\|\dD\diff\|_2^2}{\n\|\diff_{\idx}\|_1^2}
		- 16\sidx\big(2\Mra\maxerrfs+\maxerrfs^2\big)\maxzzn \,.
\end{align*}
Thus, on the set~$\bset{16\sidx\big(2\Mra\maxerrfs+\maxerrfs^2\big)\maxzzn\le\gep_1}$,
we have
\begin{align*}
	\frac{\sidx\|\dDh\diff\|_2^2}{\n\|\diff_{\idx}\|_1^2}
	\smallpad{&\ge}
	\frac{\sidx\|\dD\diff\|_2^2}{\n\|\diff_{\idx}\|_1^2}
		- \gep_1 \\
	\aequals
		\left(\frac{\sidx\diff^\top\rSd\diff}{\|\diff_{\idx}\|_1^2}-\frac{\sidx\diff^\top(\rSdb-\rSd)\diff}{\|\diff_{\idx}\|_1^2}\right) - \gep_1 \\
	\smallpad{&\ge}
		\frac{\sidx\diff^\top\rSd\diff}{\|\diff_{\idx}\|_1^2}-\frac{\sidx\|\rSdb-\rSd\|_\infty\|\diff\|_1^2}{\|\diff_{\idx}\|_1^2} - \gep_1 \\
	\smallpad{&\ge}
		\frac{\sidx\diff^\top\rSd\diff}{\|\diff_{\idx}\|_1^2}-16\sidx\|\rSdb-\rSd\|_\infty - \gep_1\,.
\end{align*}
From Cauchy-Schwartz we have~$\|\diff_{\idx}\|_1\le\sqrt{\sidx}\|\diff_{\idx}\|_2$ and hence that~$\|\diff_{\idx}\|_1^2\le\sidx\|\diff\|_2^2$.
Substitute this bound into the first term on the right-hand side above to obtain
\begin{align*}
	\frac{\sidx\|\dDh\diff\|_2^2}{\n\|\diff_{\idx}\|_1^2}
	\smallpad{&\ge}
	\frac{\diff^\top\rSd\diff}{\|\diff\|_2^2}-16\sidx\|\rSdb-\rSd\|_\infty - \gep_1 \\
	\smallpad{&\ge}
	\eigmin(\rSd) -16\sidx\|\rSdb-\rSd\|_\infty - \gep_1\,,
\end{align*}
where~$\eigmin(\rSd)$ denotes the minimal eigenvalue of~$\rSd$.
The right-hand side above does not depend on~$\diff$, so we may take the infimum of the left-hand side above over~$\diff\in\cone(\idx)$ to write, for any $\gep_1,\gep_2>0$ as in the statement of the present Lemma,
\begin{align*}
 & \Prb\left\{\ccmin(\dDh,\idx)<\eigmin(\rSd)-\gep_2-\gep_1\right\} \\
  &\qquad\qquad\qquad\le\ \,
    \Prb\bset{16\sidx(2\Mra\maxerrfs + \maxerrfs^2)\maxzzn > \gep_1} \\
		&\qquad\qquad\qquad\qquad + \Prb\bset{16\sidx\|\rSdb-\rSd\|_\infty > \gep_2} \nonumber\,,
\end{align*}
% \begin{multline*}
% 	\Prb\left\{\ccmin(\dDh,\idx,\ccc)<\eigmin(\rSd)-(1+\ccc)^2\gep_2-\gep_1\right\} \smallpad\le
% 	\Prb\bset{\sidx\|\rSdb-\rSd\|_\infty > \gep_2}
% \end{multline*}
as claimed.
\end{proof}

\noindent Lemma~\ref{lem:cc} may be combined with results for the maximum first-stage estimation error~$\maxerrfs$ and the maximum entry-wise difference~$\|\rSdb-\rSd\|_\infty$ to obtain specific bounds for~$\ccmin(\dDh,\Sb)$ under different error and design matrix regimes,
such as in Lemma~\ref{ccssg}.
}%

\begin{lem}[Second-stage compatibility constant]\label{ccssg}
Suppose that the~$\dz_i$ and~$\dv^j$ satisfy Assumptions~\ref{assu:dzi} and~\ref{assu:subgaussian}, respectively.
Set~$\tunej$ according to Definition~\ref{tpss} for each~$j\in[\px]$; set~$\tunea=\maxjpx\tunej$.
Let~$\sqrt{(\log\pz)/\n}=o(1)$.
Then, for~$\n$ sufficiently large,
\begin{multline*}
	\Prb\Big\{\ccmin(\dDh,\Sb) < \eigmin(\rSd) - \big(a + 384\Mra\sb\scca\tolV \maxzzn^{3/2}\big)\sqrt{(\log\pz)/\n}\Big\} \\
	\smallpad\le
	\prbepallvj + 2\px^{2-a^2/(6e^2\sendb^2)}\,,
\end{multline*}
where~$a>0$ is a controlled quantity,~$\eigmin(\rSd)$ denotes the minimal eigenvalue of~$\rSd$, and~$\sendb=\Mra^2(2\sgnz^2+\|\rSz\|_\infty/\log2)$.
\end{lem}

\noindent
Lemma~\ref{ccssg} entails that, under mild growth conditions~$\ccmin(\dDh,\Sb)$ is bounded below by a sequence of quantities approaching~$\eigmin(\rSd)$, and hence that~$\dDh$ satisfies the compatibility condition with probability approaching one.
Note that we must also choose the controlled quantity~$a$ above so that the exponent~${2-a^2/(6e^2\sendb^2)}$ is negative.
The sub-Gaussian regime of the present paper entails that we may choose such an~$a$ of constant order.
Throughout the present essay we use the symbol~$a$ in various bounds to denote a controlled quantity that plays this role as above.

\myappendix{
\begin{proof}[Proof of Lemma~\ref{ccssg}]
Set~$\gep_1$ in the statement of Lemma~\ref{lem:cc} as
\benn
  \gep_1 \equals 128\sb(\Mra\scca\tunea+2\frac{\sa^2}{\ccompa^4}\tunea^2)\maxzzn\,,
\eenn
so that
\begin{equation*}
  \Prb\bset{16\sb(2\Mra\maxerrfs + \maxerrfs^2)\maxzzn > \gep_1} \\
  \smallpad{\le}
  \Prb\bset{\maxerrfs > 4\scca\tunea} \,.
\end{equation*}
The present choices of tuning parameters, along with the estimates of Lemma~\ref{lem:eefsg} and Lemma~\ref{lem:cc}, entail that
\begin{multline*}
	\Prb\left\{\ccmin(\dDh,\Sb)<\eigmin(\rSd)-\gep_2-128\sb(\Mra\scca\tunea+2\frac{\sa^2}{\ccompa^4}\tunea^2)\maxzzn\right\} \\
	\smallpad{\le}
	\prbepallvj + \Prb\bset{16\sidx\|\rSdb-\rSd\|_\infty > \gep_2} \,.
\end{multline*}
Cite the growth assumptions of the present lemma to observe that, for~$\n$ sufficiently large, we have~$\Mra\scca\tunea+2\big(\scca\big)^2\tunea^2\le3\Mra\scca\tunea$\,,
from which it follows that, for such~$\n$,
\begin{multline*}
	\Prb\Big\{\ccmin(\dDh,\Sb) < \eigmin(\rSd) - \gep_2 - 384\Mra\sb\scca\tolV \maxzzn^{3/2}\sqrt{(\log\pz)/\n}\Big\} \\
	\smallpad\le
	\prbepallvj + \Prb\bset{16\sb\|\rSdb-\rSd\|_\infty > \gep_2} \,,
\end{multline*}
Cite a slight modification of Lemma~\ref{concentration} to conclude that
\begin{multline*}
	\Prb\Big\{\ccmin(\dDh,\Sb) < \eigmin(\rSd) - \big(a + 384\Mra\sb\scca\tolV \maxzzn^{3/2}\big)\sqrt{(\log\pz)/\n}\Big\} \\
	\smallpad\le
	\prbepallvj + 2\px^{2-a^2/(6e^2\sendb^2)}\,,
\end{multline*}
where~$a>0$ is a controlled quantity, as claimed.
\end{proof}
}

%% file: inference/remainders.tex
\subsection{Remainder terms and scale factors}\label{ss:remainders}
\myappendix{
\subsection{Materials required for Section~\ref{ss:remainders}}\label{proofs:remainders}
}

The asymptotic results of Section~\ref{normality} depend on the high-level assumption that the remainder terms~$\rem_\ell$ and satisfy~$\|\rem_\ell\|_\infty=o_{\Prb}(1)$.
In this Section, we identify the specific conditions under which this assumption is satisfied for the two-stage Lasso.
The primary goal of these conditions, which we present in Assumption~\ref{gr} below, is to ensure the~$\ell_1$ consistency of the first- and second-stage estimators and of the estimator~$\rTh$ specified in Section~\ref{ss:inverse}.
We implicitly refer to~$\n$-indexed sequences of all quantities mentioned below.

\begin{assu}[Model regularity for inference of $\beta_{j}$]\label{gr}
\begin{enumerate*}[label=(\arabic*)]
	\item\label{gr:cc} Assumption~\ref{assu:ccss} holds;
	\item\label{gr:satunea}
		The growth condition~$\maxjpx \saj\tunej = o(1)$ holds;
	\item\label{gr:uni}
		The sequence of population quantities~$\rT$ satisfies~$\rT \in \uni(\Mrt, \tuni, \suni)$ for a universal constant $\Mrt$,~$\suni>0$ and controlled~$\tuni\in[0,1)$;
	\item\label{gr:feas} The condition~$\Prb\,\Tinv(\tf)^\comp = o(1)$ holds;
	\item\label{gr:zzn} The quantity~$\maxzzn=\|\En[\dz_i\dz_i^\top]\|_\infty$ satisfies~$\lim_{n\to\infty}\Prb\set{\maxzzn>\maxsdz}=0$ for a universal constant~$\maxsdz$;
	\end{enumerate*}
	and
	\begin{enumerate}[label=(\arabic{*})]
	\setcounter{enumi}{5}
	\item\label{gr:growthg} The following growth conditions hold:
	\begin{enumerate}
		\item\label{gg1} $\tf^{1-\tuni}\suni\sqrt{\log\pz} = o(1)$;
		\item\label{gg2} $\sa^3\sb(\log\pz)^{3/2}/\n + \sa^2\sb(\log\pz)/\sqrt{\n} = o(1)$
		\item\label{gg3} $\mu\sb\big(\sa^2\log\pz/\sqrt{\n} + \sa\sqrt{\log\pz}\big) = o(1)$.
	\end{enumerate}
	\end{enumerate}
\end{assu}

\noindent Condition~\ref{gr:cc} is a prerequisite for the bounds on the first- and second-stage estimation errors; we discuss the feasibility of these assumptions in Section~\ref{second}.
Condition~\ref{gr:satunea} is required for asymptotic control of the remainder terms and is comparable to typical growth rates required for Lasso consistency.
Condition~\ref{gr:uni} is required to control~$\rth_j-\rt_j$ under~$\ell_\infty$ and~$\ell_1$ norms as discussed in Section~\ref{ss:inverse}.
Condition~\ref{gr:feas} is a high-level requirement for asymptotic negligibility of the remainder terms~$\rem_\ell$; it can be obtained as a consequence of specific model assumptions as in Lemma~\ref{feassgg}.
Condition~\ref{gr:zzn} is similarly a high-level condition required for asymptotic negligibility of the remainder terms: it ensures that the empirical quantity~$\maxzzn=\|\dZ^\top\dZ/\n\|_\infty$ behaves in probability as of constant order.
It can be derived as a consequence of the standard requirement that the minimal and maximal eigenvalues of~$\rSz$ be bounded away from zero and infinity uniformly in $\n$ if~$\|\rSzh-\rSz\|_\infty=o_{\Prb}(1)$; the latter condition can in turn be derived from distributional assumptions on the~$\dz_i$ using the tools of, say, \cite{V12}.
Condition~\ref{gr:growthg} lists the model parameter growth conditions required for asymptotic negligibility of the remainder terms under the sub-Gaussian noise regime of Assumption~\ref{assu:subgaussian}.
\new{We can compare these conditions with the requirement~$s\log p / \sqrt{n}=o(1)$ in~\cite{JM14, vdGetal14} for negligibility of the single remainder term that occurs under the ordinary linear model.
The conditions on the sparsity here are generally similar to those of the ordinary linear model, yet slightly more strict when comparing the powers at which the sparsity factors enter.
It is not clear if our conditions can be relaxed further, or whether there are more fundamental reasons for the differences.}

Note that while the quantity~$\maxsdz$ of Condition~\ref{gr:zzn} appears in the bounds of Lemmas \ref{lem:rem1g}-\ref{lem:rem4g} of the Appendix, which give the rates for the remainder terms, we do not include it in the growth conditions of Condition~\ref{gr:growthg}.
This is because, under the presently studied regime,~$\maxsdz$ is assumed of constant order.
One could consider more general scenarios where the maximum entry of~$\rSzh$ is not bounded in probability and include the quantity~$\maxsdz$ in the aforementioned growth conditions.
Doing so would in turn affect the rate at which~$\sa,\sb$, and~$\pz$ may be allowed to grow with~$\n$ while maintaining asymptotic negligibility of the remainder terms~$\rem_\ell$.

In addition to the model regularity conditions of Assumption~\ref{gr}, we require appropriate choices of the first- and second-stage Lasso tuning parameters and the estimator~$\rTh$.
We codify such choices in the following Assumption and then conclude the asymptotic negligibility of the remainder terms.

\begin{assu}[Specification of estimators]\label{assum:estimators}
Let~$\rAh$ and~$\rbla$ be the first- and second-stage Lasso estimators, respectively. The tuning parameters under the sub-Gaussian noise regime of Assumption~\ref{assu:subgaussian} are chosen according to
\begin{enumerate*}[label=(\roman*)]
	\item Definition~\ref{tpss} for the first-stage tuning parameters~$\tunefs = (\tunej)_{j=1}^{\px}$ and the quantity~$\tunea$;
	\item Definition~\ref{tpss} for the second-stage tuning parameter~$\tuneb$ and quantities~$\tuneu$ and~$\tuneV$;
\end{enumerate*}
(iii) let ~$\rTh$ be an estimator of~$\rT$ with rows~$\rth_j$ given by solutions to Program~\ref{program}.
\end{assu}

%%%%%%%%%%%%%%%%%%%%%%%%%%%%%%%%%%%%%%%%%%%%%%%%%%%%%%%%%%%%%%%%%%%%%%%%%%%%%%%

\myappendix{%

Lemmas~\ref{lem:rem1}, \ref{lem:rem2}, \ref{lem:rem3}, and \ref{lem:rem4} provide finite-sample bounds for the quantities~$\|\rem_\ell\|_\infty$ for~$\ell\in[4]$ that are generic over various noise regimes.
We present them separately for the sake of modularity with respect to such assumptions.
Lemmas~\ref{lem:rem1g}, \ref{lem:rem2g}, \ref{lem:rem3g}, and~\ref{lem:rem4g} in turn provide specific rates for the~$\|\rem_\ell\|_\infty$ under the sub-Gaussian noise regime of Assumption~\ref{assu:subgaussian}.

\begin{lem}[Control of $\rem_1$]\label{lem:rem1}
Suppose that Assumption~\ref{assu:ccss} and Conditions~\ref{gr:satunea}, \ref{gr:uni} of Assumption~\ref{gr} hold and that~$\rTh$ is chosen according to Assumption~\ref{assum:estimators}.
Then, on the set $\Tu(\tuneu)\cap\Tinv(\tf)$, where $\tuneu >0$ is arbitrary, the remainder term
\benn
	\rem_1 \equals (\rTh-\rT)^\top\dD^\top\du/\sqrt{n}
\eenn
satisfies
\benn
	\|\rem_1\|_\infty \relates{\le} 2^{-\tuni}\sqrt{n}\Mra\cprcr(\Mrt\tf)^{1-\tuni}\suni\tuneu \,,
\eenn
where $\cprcr$ is as in Lemma~\ref{lem:eeprc}.
\end{lem}
}%

%%%%%%%%%%
% PROOF
\myappendix{
\begin{proof}[Proof of Lemma~\ref{lem:rem1}]
Lemma~\ref{lem:eeprc} entails that, on the set $\Tinv(\tf)$,
\benn
	\maxjpx\|\rth_j-\rt_j\|_1 \relates{\le} 2\cprcr (2\Mrt\tf)^{1-\tuni}\suni \,.
\eenn
We therefore find that
\begin{align*}
	\|(\rTh-\rT)^\top\dD^\top\du/\sqrt{n}\|_\infty
	\smallpad{&\le}
	\sqrt{n}\|\rTh-\rT\|_{L_1}\|\dD^\top\du/n\|_\infty \\
	\smallpad{&\le}
		\sqrt{n}\|\rTh-\rT\|_{L_1}\|\rA\|_{L_1}\|\dZ^\top\du/n\|_\infty \\
 	\smallpad{&\le}
		\sqrt{n}\Mra\|\rTh-\rT\|_{L_1}\maxzn{\du} \\
	\smallpad{&\le}
	 	2\sqrt{n}\Mra\cprcr(2\Mrt\tf)^{1-\tuni}\suni\maxzn{\du} \,.
\end{align*}
On the set~$\Tu(\tuneu)$ we have~$\maxzn{\du} \le \tuneu/4$.
From this bound and the previous display we conclude that, on the set~$\Tu(\tuneu)\cap\Tinv(\tf)$,
\benn
	\|(\rT-\rTh)^\top\dD^\top\du/\sqrt{n}\|_\infty \relates{\le} 2^{-\tuni}\sqrt{n}\Mra\cprcr(\Mrt\tf)^{1-\tuni}\suni\tuneu \,,
\eenn
as claimed.
\end{proof}
}%

%%%%%%%%%%%%%%%%%%%%%%%%%%%%%%%%%%%%%%%%%%%%%%%%%%%%%%%%%%%%%%%%%%%%%%%%%%%%%%%
%%%%%%%%%%%%%%%%%%%%%%%%%%%%%%%%%%%%%%%%%%%%%%%%%%%%%%%%%%%%%%%%%%%%%%%%%%%%%%%

\myappendix{
\begin{lem}[Control of $\rem_1$, sub-Gaussian noise]\label{lem:rem1g}
Suppose that (i) Assumption~\ref{assu:ccss} and Conditions~\ref{gr:satunea}, \ref{gr:uni} of Assumption~\ref{gr} hold and (ii)~Assumption~\ref{assu:subgaussian} holds.
Choose~$\rTh$ according to Assumption~\ref{assum:estimators}.
Then,
\begin{equation*}
	\Prb\left\{\|\rem_1\|_\infty > 2^{-\tuni}\Mra\cprcr\tolu (\Mrt\tf)^{1-\tuni}\suni\sqrt{\maxsdz\log\pz} \right\}
	\smallpad{\le} e\pz^{1-\tolu^2\Choeff/\sgnu^2}+\Prb\,\Tinv(\tf)^\comp \,,
%\Prb\big\{\maxzzn > \maxsdz\big\} + \Prb\,\Tinv(\tf)^\comp \,,
\end{equation*}
% \begin{align*}
% 	\smallpad{&\ }
% 	\Prb\left\{\|\rem_1\|_\infty > 2^{-\tuni}\Mra\cprcr\sdu c(\Mrt\tf)^{1-\tuni}\suni\sqrt{\maxsdz\log\pz} \right\} \\
% 	\smallpad{\le&\ }
% 	2\pz^{-\cepg} + \Prb\big\{\maxzzn > \maxsdz\big\} + \Prb\,\Tinv(\tf)^\comp \,,
	% \smallpad{&\ }
	% \Prb\left\{\|\rem_1\|_\infty > 2^{-\tuni}\Mra\cprcr\sdu c(\Mrt\tf)^{1-\tuni}\suni\sqrt{\maxzzn\log\pz} \right\} \\
	% \smallpad{\le&\ }
	% 2\pz^{-\cepg} + \Prb\,\Tinv(\tf)^\comp \,,
% \end{align*}
where~$\tolu$ is as in Definition~\ref{tpss} and~$\Choeff$ is as defined in Lemma~\ref{lem:epg}.
If Conditions~\ref{gr:cc}, \ref{gr:feas}, \ref{gr:zzn}, and~\ref{gg1} of Assumption~\ref{gr} also hold, then~$\|\rem_1\|_\infty = o_{\Prb}(1)$.
\end{lem}
}%

%%%%%%%%%%
% PROOF
\myappendix{%
\begin{proof}[Proof of Lemma~\ref{lem:rem1g}]
Lemma~\ref{lem:rem1} entails that
\begin{align*}
	\Prb\left\{\|\rem_1\|_\infty > 2^{-\tuni}\sqrt{n}\Mra\cprcr(\Mrt\tf)^{1-\tuni}\suni\tuneu \right\} \
	\leq \ \Prb\,\Tu^{\comp} + \Prb\,\Tinv(\tf)^\comp \,.
\end{align*}
Substitute~$\tuneu$ chosen according to Definition~\ref{tpss} into the display above and cite the estimate of Lemma~\ref{lem:epu} to deduce the original claim.
\end{proof}
}%

%%%%%%%%%%%%%%%%%%%%%%%%%%%%%%%%%%%%%%%%%%%%%%%%%%%%%%%%%%%%%%%%%%%%%%%%%%%%%%%
%%%%%%%%%%%%%%%%%%%%%%%%%%%%%%%%%%%%%%%%%%%%%%%%%%%%%%%%%%%%%%%%%%%%%%%%%%%%%%%

\myappendix{%
\begin{lem}[Control of $\rem_2$]\label{lem:rem2}
Suppose that Assumption~\ref{assu:ccss} and Condition~\ref{gr:uni} of Assumption~\ref{gr} hold.
Choose~$\rTh$ according to Assumption~\ref{assum:estimators}, set~$\tunea$ according to Definition~\ref{tpss}, and let~$\tuneu > 0$ be arbitrary.
Then, on the set $\TV(\tunefs) \cap \Tu(\tuneu) \cap \Tinv(\tf)$, the remainder term
\benn
	\rem_2 \equals \rTh^\top(\dDh-\dD)^\top\du/\sqrt{n}
\eenn
satisfies
\benn
	\|\rem_2\|_\infty \smallpad{\le} \sqrt{n}\Mrt\scca\tunea\tuneu
\eenn
for $\n$ sufficiently large.
\end{lem}
}%
%%%%%%%%%%
% PROOF
\myappendix{
\begin{proof}[Proof of Lemma~\ref{lem:rem2}]
Observe that
\begin{align*}
	\rTh^\top(\dDh-\dD)^\top\du/\sqrt{n} \equals \sqrt{n}\rTh^\top(\rAh-\rA)^\top(\dZ^\top\du/n) \,.
\end{align*}
On the set $\Tinv(\tf)$, each row $\rt$ is feasible for Program~\ref{program}.
Then, $\|\rth_j\|_1 \le \|\rt_j\|_1$ for each $j\in[\px]$ by specification.
Lemmas \ref{lem:eefs} and \ref{lem:epu} then entail that, on the set $\TA(\tunefs) \cap \Tu(\tuneu) \cap \Tinv(\tf)$,
\begin{align*}
	\|\rTh^\top(\dDh-\dD)^\top\du/\sqrt{n}\|_\infty \ & \leq \ \maxjkpx \sqrt{n}\|\rth_j\|_1\|\rah^k-\ra^k\|_1\|\dZ^\top\du/n\|_\infty \\
	\smallpad{&\le}
		\sqrt{n}\Mrt\scca\tunea\tuneu \,,
\end{align*}
as claimed.
\end{proof}
}%

%%%%%%%%%%%%%%%%%%%%%%%%%%%%%%%%%%%%%%%%%%%%%%%%%%%%%%%%%%%%%%%%%%%%%%%%%%%%%%%
%%%%%%%%%%%%%%%%%%%%%%%%%%%%%%%%%%%%%%%%%%%%%%%%%%%%%%%%%%%%%%%%%%%%%%%%%%%%%%%

\myappendix{
\begin{lem}[Control of~$\rem_2$, sub-Gaussian noise]\label{lem:rem2g}
Suppose that (i) Assumption~\ref{assu:ccss} and Condition~\ref{gr:uni} of Assumption~\ref{gr} hold and (ii) Assumption~\ref{assu:subgaussian} holds.
Choose~$\rTh$,~$\tunefs = (\tune_1,\ldots,\tune_{\px})$, and~$\tunea$ according to Assumption~\ref{assum:estimators}.
Then,
\begin{multline*}
	\Prb\Big\{ \|\rem_2\|_\infty > \Mrt \tolV\tolu \maxsdz\scca\log\pz/\sqrt{\n}\Big\} \\
	\smallpad{\le}
	\prbepV + \prbepu + \Prb\big\{\maxzzn > \maxsdz\big\} + \Prb\,\Tinv(\tf)^\comp \,,
\end{multline*}
where~$\tolV$ is as in Definition~\ref{tpss},~$\tolu$ is as in Definition~\ref{tpss}, and~$\Choeff$ is as defined in Lemma~\ref{lem:epg}.
Consequently, if Conditions~\ref{gr:cc}, \ref{gr:feas}, \ref{gr:zzn}, and~\ref{gg2} of Assumption~\ref{gr} also hold, then $\|\rem_2\|_\infty = o_{\Prb}(1)$.
\end{lem}
}%

%%%%%%%%%%
% PROOF
\myappendix{
\begin{proof}[Proof of Lemma~\ref{lem:rem2g}]
Lemma~\ref{lem:rem2} entails that
\benn
	\Prb\left\{\|\rem_2\|_\infty > \foq\sqrt{n}\Mrt\scca\tunea\tuneu\right\} \relates{\le} \Prb\,\TA(\tunefs)^\comp + \Prb\,\Tu(\tuneu)^\comp + \Prb\,\Tinv(\tf)^\comp \,.
\eenn
Substitute the present choices of~$\tunea$ and~$\tuneu$ into the previous display and cite the estimates of Lemmas~\ref{lem:epu} and~\ref{lem:epV} to deduce the original claim.
\end{proof}
}%

%%%%%%%%%%%%%%%%%%%%%%%%%%%%%%%%%%%%%%%%%%%%%%%%%%%%%%%%%%%%%%%%%%%%%%%%%%%%%%%
%%%%%%%%%%%%%%%%%%%%%%%%%%%%%%%%%%%%%%%%%%%%%%%%%%%%%%%%%%%%%%%%%%%%%%%%%%%%%%%

\myappendix{%
\begin{lem}[Control of $\rem_3$]\label{lem:rem3}
Suppose that Assumptions~\ref{assu:ccss} and \ref{assu:ccss} and Conditions~\ref{gr:satunea} and~\ref{gr:uni} of Assumption~\ref{gr} hold.
Choose~$\rTh$ according to Assumption~\ref{assum:estimators}; let $\tunefs = (\tune_1,\ldots,\tune_{\px}) > \0, \tuneu > 0$, and~$\tuneV>0$ be arbitrary.
Set
\benn
\tuneb \smallpad{=}
16\scca\tunea
	\maxzzn \big(4\Mb\scca\tunea+\Mra\big)
	+ \big(4\scca\tunea+\Mra\big)
	\big(\Mb\tuneV+\tuneu\big)\,.
\eenn
Then, on the set
\benn
	 \TV(\tunefs)\cap\Tu(\tuneu)\cap\Tinv(\tf)\,,
\eenn
the remainder term
\benn
	\rem_3 \equals \rTh\dDh^\top(\dX-\dDh)(\rb-\rbla)/\sqrt{n}
\eenn
satisfies
\begin{align*}
	\|\rem_3\|_\infty
	\smallpad{\le}
	8\Mrt\Mra\sqrt{\n}(4\scca\tunea+\tuneu)\sccb\tuneb \,.
\end{align*}
\end{lem}
}%

%%%%%%%%%%
% PROOF
\myappendix{
\begin{proof}[Proof of Lemma~\ref{lem:rem3}]
We first observe that
\begin{align*}
	\|\rTh\dDh^\top(\dX-\dDh)(\rb-\rbla)/\sqrt{n}\|_\infty \ & \leq \ \sqrt{n} \|\rTh\|_{L_1} \|\dDh^\top(\dX-\dDh)/n\|_\infty\|\rbla-\rb\|_1 \,.
\end{align*}
Now,
\begin{align*}
	\dDh^\top(\dX-\dDh)/n \ & = \ \dDh^\top(\dD+\du-\dDh)/n \\
	& = \ \dDh^\top(\dD-\dDh)/n + \dDh^\top\du/n \\
	& = \ \rAh^\top(\dZ^\top\dZ/n)(\rA-\rAh) + \dDh^\top\du/n \,.
\end{align*}
For the first term on the right-hand side above, write
\begin{align*}
	\|\rAh^\top(\dZ^\top\dZ/n)(\rA-\rAh)\|_\infty
	\smallpad{&\le}
	\|\rA^\top(\dZ^\top\dZ/n)(\rA-\rAh)\|_\infty \\
	&\quad
	+ \|(\rAh-\rA)^\top(\dZ^\top\dZ/n)(\rA-\rAh)\|_\infty \\
	% \smallpad{&\le}
		% \maxjpx\maxkpx \|\ra^j\|_1\|\rah^k-\ra^k\|_1 + \left(\maxjpx \|\rah^j-\ra^j\|_1\right)^2  \\
	\smallpad{&\le}
		\maxzzn[\Mra\maxerrfs + \maxerrfs^2] \,.
\end{align*}
On the set $\TV(\tunefs)$, the right-hand side above is less than or equal to~$2\Mra \cdot \maxzzn\maxerrfs$ for $n$ sufficiently large by Lemma~\ref{lem:eefs} and the hypotheses of the present lemma.
For the second term on the right-hand side of two displays previous, write
\begin{align*}
	\|\dDh^\top\du/n\|_\infty
	\smallpad{&=}
	\|\rAh^\top\dZ^\top\du/n\|_\infty \\
	\smallpad{&=}
		\|(\rA+[\rAh-\rA])^\top(\dZ^\top\du/n)\|_\infty \\
	\smallpad{&\le}
		2\Mra\maxzn{\du} \,,
\end{align*}
where the final line holds on the set~$\TV(\tunefs)$ for such~$n$ by Lemma~\ref{lem:eefs}.
Thus, on the set~$\TV(\tunefs)\cap\Tu(\tuneu)$, we have
\begin{align*}
	\|\dDh^\top(\dX-\dDh)/n\|_\infty
	\smallpad{&\le}
	2\Mra(\maxzzn\maxerrfs+\maxzn{\du}) \\
	\smallpad{&\le}
	2\Mra(4\maxzzn\scca\tunea+\tuneu) \,,
\end{align*}
where the latter substitutions are justified by Lemma~\ref{lem:eefs} and the definition of $\Tu(\tuneu)$.

On the set $\Tinv(\tf)$, each row $\rt$ is feasible for Program~\ref{program}.
Then, $\|\rth_j\|_1 \le \|\rt_j\|_1$ for each $j\in[\px]$ by specification.

Finally, Lemma~\ref{lem:eess} entails that, for the present choice of $\tuneb$,
\benn
	\|\rbla-\rb\|_1
	\smallpad{\le}
	4\sccb\tuneb
\eenn
on the set $\TV(\tuneV)\cap\Tu(\tuneu)$.

Combining the foregoing results, we see that, on the set
\benn
	\TA(\tunefs)\cap\Tu(\tuneu)\cap\TV(\tuneV)\cap\Tinv(\tf)\,,
\eenn
it holds that
\begin{align*}
	\|\rTh\dDh^\top(\dX-\dDh)(\rb-\rbla)/\sqrt{n}\|_\infty
	\smallpad{\le}
	8\Mrt\Mra\sqrt{\n}(4\maxzzn\scca\tunea+\tuneu)\sccb\tuneb \,.
\end{align*}
Now note that, under the present choices of~$\tunefs$ and~$\tuneV$, the set~$\TV(\tuneV)$ is contained in~$\TV(\tunefs)$.
\end{proof}
}%

%%%%%%%%%%%%%%%%%%%%%%%%%%%%%%%%%%%%%%%%%%%%%%%%%%%%%%%%%%%%%%%%%%%%%%%%%%%%%%%
%%%%%%%%%%%%%%%%%%%%%%%%%%%%%%%%%%%%%%%%%%%%%%%%%%%%%%%%%%%%%%%%%%%%%%%%%%%%%%%

\myappendix{
\begin{lem}[Control of $\rem_3$, sub-Gaussian noise]\label{lem:rem3g}
Suppose that Conditions~\ref{gr:cc}, \ref{gr:satunea}, and~\ref{gr:uni} of Assumption~\ref{gr} hold and (ii) Assumption~\ref{assu:subgaussian} holds.
Choose~$\rTh$,~$\rAh$ and~$\rbla$ according to Assumption~\ref{assum:estimators}.
Then,
\begin{align*}
	& \Prb\big\{\|\rem_3\|_\infty >
		8\Mrt\Mra\sqrt{\n}(4\maxsdz\scca\tunea+\tuneu)\sccb\tuneb\big\} \\
	&\qquad\qquad\qquad\qquad\qquad\qquad \smallpad{\le}
	\prbepallvj + \prbepu
	+ \Prb\,\Tinv(\tf)^\comp \,,
\end{align*}
where~$\tolV$ is as in Definition~\ref{tpss},~$\tolu$ is as in Definition~\ref{tpss}, and~$\Choeff$ is as defined in Lemma~\ref{lem:epg}.
% \begin{multline*}
% 	\Prb\big\{\|\rem_3\|_\infty >
% 		8\Mrt\Mra\sqrt{\n}(4\maxsdz\scca\tunea+\tuneu)\sccb\tuneb\big\} \\
% 	\smallpad{\le}
% 	\prbepallvj + \prbepu + \Prb\big\{\maxzzn > \maxsdz\big\} + \Prb\,\Tinv(\tf)^\comp \,.
% 	% \smallpad{&\ } \Prb\big\{\|\rem_3\|_\infty >
% 	% 	8\Mrt\Mra\sqrt{\n}(4\maxzzn\scca\tunea+\tuneu)\sccb\tuneb\big\} \\
% 	% \smallpad{\le&\ }
% 	% 2\pz^{1-\cepg} + 2\pz^{-\cepg} + \Prb\,\Tinv(\tf)^\comp \,.
% \end{multline*}
Consequently, if Conditions~\ref{gr:feas}, \ref{gr:zzn}, and~\ref{gg2} of Assumption~\ref{gr} also hold, then~$\|\rem_3\|_1 = o_{\Prb}(1)$.
\end{lem}
}%

%%%%%%%%%%
% PROOF
\myappendix{
\begin{proof}[Proof of Lemma~\ref{lem:rem3g}]
Lemma~\ref{lem:rem3} entails that
\begin{multline*}
	\Prb\big\{\|\rem_3\|_\infty >
		8\Mrt\Mra\sqrt{\n}(4\maxsdz\scca\tunea+\tuneu)\sccb\tuneb\big\}\\
	% & \leq \ \Prb\TV^\comp + \Prb\setcomp{\Tu} \,.
	\smallpad{\le}
	\Prb\,\TV(\tuneV)^\comp + \Prb\,\Tu(\tuneu)^\comp + \Prb\,\Tinv(\tf)^\comp \,.
\end{multline*}
Substitute the present choices of tuning parameters into the display above and cite the estimates of Lemmas~\ref{lem:epu} and~\ref{lem:epV}, to deduce the first claim.
Expand the the present choices of tuning parameters to find
\begin{align*}
	&8\Mrt\Mra\sqrt{\n}(4\maxsdz\scca\tunea+\tuneu)\sccb\tuneb \\
	&\qquad\qquad\qquad\qquad \lesssim\,
		\sa^3\sb(\log\pz)^{3/2}/\n + \sa^2\sb(\log\pz)/\sqrt{\n} \\
		&\qquad\qquad\qquad\qquad\qquad+ \sa^2\sb(\log\pz)^{3/2}/\n + \sa\sb\log\pz/\sqrt{n} \,,
\end{align*}
from which the latter claim follows.
\end{proof}
}%

%%%%%%%%%%%%%%%%%%%%%%%%%%%%%%%%%%%%%%%%%%%%%%%%%%%%%%%%%%%%%%%%%%%%%%%%%%%%%%%
%%%%%%%%%%%%%%%%%%%%%%%%%%%%%%%%%%%%%%%%%%%%%%%%%%%%%%%%%%%%%%%%%%%%%%%%%%%%%%%

\myappendix{%
\begin{lem}[Control of $\rem_4$]\label{lem:rem4}
Suppose that Assumption~\ref{assu:ccss} and Condition~\ref{gr:uni} of Assumption~\ref{gr} hold.
Choose~$\rTh$ according to Assumption~\ref{assum:estimators}; let~$\tuneV, \tuneu > 0$ be arbitrary; set
\benn
\tuneb \smallpad{=}
16\scca\tunea
	\maxzzn \big(4\Mb\scca\tunea+\Mra\big)
	+ \big(4\scca\tunea+\Mra\big)
	\big(\Mb\tuneV+\tuneu\big)\,.
\eenn
Then, on the set~$\TV(\tuneV)\cap\Tu(\tuneu)\cap\Tinv(\tf)$, the remainder term
\benn
	\rem_4 \equals \sqrt{n}(\rTh\rSdh - \Id)(\rb-\rbla)
\eenn
satisfies
\benn
	\|\rem_4\|_\infty
	\smallpad{\le}
	4\sqrt{n}\tf\sccb\tuneb \,.
\eenn
\end{lem}
}%

%%%%%%%%%%
% PROOF
\myappendix{%
\begin{proof}[Proof of Lemma~\ref{lem:rem4}]
Note first that
\begin{align*}
	\|\sqrt{n}(\rTh\rSdh - \Id)(\rb-\rbla)\|_\infty
	\smallpad{&\le}
	\sqrt{n}\|\rTh\rSdh - \Id\|_\infty \|\rb-\rbla\|_1 \\
	\smallpad{&\le}
		\sqrt{n}\tf \|\rb-\rbla\|_1 \,,
\end{align*}
where the latter inequality follows from the specification of $\rTh$ and the fact that Program~\ref{program} is feasible given $\ta$.
By Lemma~\ref{lem:eess}, on the set~$\TV(\tuneV)\cap\Tu(\tuneu)\cap\Tinv(\tf)$,
\benn
	\|\rb-\rbla\|_1 \smallpad\le 4\sccb\tuneb \,.
\eenn
Combine the two previous displays to deduce the original claim.
\end{proof}
}%

%%%%%%%%%%%%%%%%%%%%%%%%%%%%%%%%%%%%%%%%%%%%%%%%%%%%%%%%%%%%%%%%%%%%%%%%%%%%%%%
%%%%%%%%%%%%%%%%%%%%%%%%%%%%%%%%%%%%%%%%%%%%%%%%%%%%%%%%%%%%%%%%%%%%%%%%%%%%%%%

\myappendix{
\begin{lem}[Control of $\rem_4$, Gaussian noise]\label{lem:rem4g}
Suppose that (i) Assumption~\ref{assu:ccss} and Condition~\ref{gr:uni} of Assumption~\ref{gr} hold and (ii) Assumption~\ref{assu:subgaussian} holds.
Choose~$\rTh$,~$\rAh$ and~$\rbla$ according to Assumption~\ref{assum:estimators}.
Then,
\begin{align*}
	\Prb\big\{\|\rem_4\|_\infty > 4\sqrt{n}\tf\sccb\tuneb\big\}
	\smallpad{&\le}
	\prbepV + \prbepu + \Prb\,\Tinv(\tf)^\comp \,,
\end{align*}
where~$\tolV$ is as in Definition~\ref{tpss},~$\tolu$ is as in Definition~\ref{tpss}, and~$\Choeff$ is as defined in Lemma~\ref{lem:epg}.
Consequently, if Conditions~\ref{gr:feas}, \ref{gr:zzn}, and~\ref{gg3} of Assumption~\ref{gr} also hold, then~$\|\rem_4\|_1 = o_{\Prb}(1)$.
\end{lem}
}%

%%%%%%%%%%
% PROOF
\myappendix{%
\begin{proof}[Proof of Lemma~\ref{lem:rem4g}]
Lemma~\ref{lem:rem4} entails that
\benn
	\Prb\big\{\|\rem_4\|_\infty > 4\sqrt{n}\tf\sccb\tuneb\big\}
	\smallpad{\le}
	\Prb\,\TV(\tuneV)^\comp + \Prb\,\Tu(\tuneu)^\comp + \Prb\,\Tinv(\tf) \,.
\eenn
Substitute the present choices of tuning parameters into the display above and cite the estimates of Lemmas~\ref{lem:epu} and~\ref{lem:epV} to deduce the first claim.
Expand the the present choices of tuning parameters to find
\begin{align*}
	\sqrt{n}\tf\sccb\tuneb
	\smallpad{&\lesssim}
	\mu\sb\big(\sa^2\log\pz/\sqrt{\n} + \sa\sqrt{\log\pz}\big)\,,
\end{align*}
from which the second claim follows.
\end{proof}
}%

%%%%%%%%%%%%%%%%%%%%%%%%%%%%%%%%%%%%%%%%%%%%%%%%%%%%%%%%%%%%%%%%%%%%%%%%%%%%%%%
%%%%%%%%%%%%%%%%%%%%%%%%%%%%%%%%%%%%%%%%%%%%%%%%%%%%%%%%%%%%%%%%%%%%%%%%%%%%%%%

% \noindent We can now conclude the asymptotic negligibility of the remainder terms~$\rem_\ell$.

\begin{lem}[Negligibility of remainders for Theorem \ref{wl}]\label{negligibility}
Suppose that Assumption~\ref{assu:subgaussian} and Conditions~\ref{gr:cc}-\ref{gr:growthg} of Assumption~\ref{gr} hold and that
the estimators~$\rAh,\rbla,$ and~$\rTh$ are chosen according to Assumption \ref{assum:estimators} (i)-(iii).
Then,~$\|\rem_\ell\|_\infty = o_{\Prb}(1)$ for~$\ell\in[4]$.
\end{lem}

\myappendix{
\begin{proof}[Proof of Lemma~\ref{negligibility}] %Proof of Lemma 4.8
The result follows from Lemmas~\ref{lem:rem1g}, \ref{lem:rem2g}, \ref{lem:rem3g}, and~\ref{lem:rem4g}.
\end{proof}
}

\noindent The primary use of Lemma~\ref{negligibility} is to verify Condition~\ref{wl:cond2} of Theorem~\ref{wl}.
Indeed, the result justifies the use of the one-step update to the second-stage Lasso estimator to construct asymptotically valid confidence intervals for the components~$\rbj$ according to~\eqref{confidence}.

We note that the quantity~$\tf$, which we recall is the tolerance parameter for Program~\ref{program}, must be given careful consideration.
Conditions~\ref{gg1} and~\ref{gg2} of Assumption~\ref{gr} require $\tf$ to be of small order~$(\suni\sqrt{\log\pz})^{\frac{1}{\tuni-1}}$ and~$(\sb\sa\log\pz)^{-1}$, respectively.
However,~$\tf$ must not tend to zero so fast that the probability that~$\rT$ is feasible for Program~\ref{program}, which we recall is formally denoted by~$\Prb\,\Tinv(\tf)$, becomes bounded away from zero.
% The growth of~$\tf$ must therefore balance two competing objectives, which is non-trivial.
The following Lemma identifies a choice of~$\tf$ that satisfies these competing objectives for sub-Gaussian~$\dz_i$ and first-stage noise elements.

\myappendix{
\noindent The following two Lemmas are required for Lemma~\ref{feassgg}.

\begin{lem}\label{feas} %LemmaB.15
For~$\gep<\Mra$, it holds that
\begin{multline*}
	\Prb\bset{\|\rT\rSdh-\Id\|_\infty > t + 3\Mrt\Mra\|\rSzh\|_\infty\gep} \\
	\smallpad\le \Prb\bset{\|\rT\rSdb-\Id\|_\infty > t} +
	\Prb\bset{\maxerrfs\leq\gep}\,,
\end{multline*}
where $\rSdb = \En[\dd_i\dd_i^\top]$.
\end{lem}
}%

%%%%%%%%%%
% PROOF
\myappendix{
\begin{proof}[Proof of Lemma \ref{feas}]%Proof of Lemma B.15
Note that
\begin{align*}
	\rSdh \equals \dDh^\top\dDh/n \aequals \dD^\top\dD/n + (\rAh-\rA)^\top\dZ^\top\dZ\rA/n + \rAh^\top\dZ^\top\dZ(\rAh-\rA)/n \\
	\aequals\rSdb + (\rAh-\rA)^\top\dZ^\top\dZ\rA/n + \rA^\top\dZ^\top\dZ(\rAh-\rA)/n \\
	& \qquad + (\rAh-\rA)^\top\dZ^\top\dZ(\rAh-\rA)/n \,,
\end{align*}
so that
\begin{align*}
	\|\rT\rSdh - \Id\|_\infty
	\smallpad{&\le}
	\|\rT(\dD^\top\dD/n) - \Id\|_\infty + \|\rT(\rAh-\rA)^\top\dZ^\top\dZ\rA/n\|_\infty \\
	& \quad
		+ \|\rT \rA^\top\dZ^\top\dZ(\rAh-\rA)/n\|_\infty \\
	&\quad	+ \|\rT(\rAh-\rA)^\top\dZ^\top\dZ(\rAh-\rA)/n\|_\infty \\
	\smallpad{&:=}
		\|\rT(\dD^\top\dD/n) - \Id\|_\infty + \I_1 + \I_2 + \I_3 \,.
\end{align*}
Note that
\begin{align*}
	\I_1
	\equals
	\|\rT(\rAh-\rA)^\top\dZ^\top\dZ\rA/n\|_\infty
	\smallpad{&\le}
		\|\rT\|_{L_1}\|(\rAh-\rA)^\top\dZ^\top\dZ\rA/n\|_\infty \\
	\smallpad{&\le}
		\Mrt \|\rSzh\|_\infty \|\rAh-\rA\|_{L_1}\|\rA\|_{L_1} \\
	\aequals
		\Mrt\Mra\|\rSzh\|_\infty\maxerrfs \,.
\end{align*}
The same bound holds for $\I_2$ by symmetry of the $\ell_\infty$ norm under transposition of its argument.

For the term $\I_3$, similar reasoning yields
\begin{align*}
	\I_3
	\aequals
	\|\rT(\rAh-\rA)^\top\dZ^\top\dZ(\rAh-\rA)/n\|_\infty \\
	\smallpad{&\le}
	\Mrt\|\rSzh\|_\infty\|\rAh-\rA\|_{L_1}^2
	\equals
	\Mrt\|\rSzh\|_\infty\maxerrfs^2 \,.
\end{align*}
If $\gep<\Mra$, then, on the set~$\set{\maxerrfs\leq\gep}$, it holds that~$\I_3 \le \Mrt\Mra \cdot \|\rSzh\|_\infty\maxerrfs$.
Conclude that
\begin{align*}
	\Prb\bset{\I_1+\I_2+\I_3 > 3\Mrt\Mra\|\rSzh\|_\infty\gep}
	\smallpad\le
	\Prb\set{\maxerrfs>\gep}
\end{align*}
and hence that
\begin{multline*}
	\Prb\bset{\|\rT\rSdh-\Id\|_\infty > t + 3\Mrt\Mra\|\rSzh\|_\infty\gep} \\
	\smallpad\le
	\Prb\bset{\|\rT\rSdb-\Id\|_\infty > t} + \Prb\set{\maxerrfs>\gep} \,,
\end{multline*}
as claimed.
\end{proof}
}%

%%%%%%%%%%%%%%%%%%%%%%%%%%%%%%%%%%%%%%%%%%%%%%%%%%%%%%%%%%%%%%%%%%%%%%%%%%%%%%%
%%%%%%%%%%%%%%%%%%%%%%%%%%%%%%%%%%%%%%%%%%%%%%%%%%%%%%%%%%%%%%%%%%%%%%%%%%%%%%%

\myappendix{
\begin{lem}\label{feassg} %Lemma B.16
Suppose that the~$\dz_i$ satisfy Assumption~\ref{assu:dzi}.
Set
\benn
	\tf \equals \Mrt a\sqrt{\log(\px)/\n} + 3\Mrt\Mra\|\rSzh\|_\infty\gep
\eenn
where~$a>0$ and~$\gep>0$ are controlled quantities.
Then
\benn
	\Prb\Tinv(\tf)
	\smallpad\le
	2\px^{2-a^2/(6e^2\sendb^2)} + \Prb\set{\maxerrfs>\gep} \,,
\eenn
where~$\sendb=\Mra^2(2\sgnz^2+\|\rSz\|_\infty/\log2)$.
\end{lem}
}

\myappendix{
\begin{proof}[Proof of Lemma~\ref{feassg}]\label{proof:feassg} %Proof of Lemma B.16
Lemma~\ref{feas} entails that
\begin{multline*}
	\Prb\bset{\|\rT\rSdh-\Id\|_\infty > t + 3\Mrt\Mra\|\rSzh\|_\infty\gep} \\
	\smallpad\le
	\underbrace{\Prb\bset{\|\rT\rSdb-\Id\|_\infty > t}}_{\I_1(t)} + \Prb\set{\maxerrfs>\gep} \,.
\end{multline*}
For~$t>0$.
Now write
\begin{align*}
	\rT\rSdb-\Id \equals \rT\rSd+\rT(\rSdb-\rSd) - \Id
	\equals \rT(\rSdb-\rSd)
\end{align*}
to infer that
\benn
	\|\rT\rSdb-\Id\|_\infty \equals
	\|\rT(\rSdb-\rSd)\|_\infty
	\smallpad\le
	\|\rT\|_{L_1}\|\rSdb-\rSd\|_\infty
	\smallpad\le
	\Mrt\|\rSdb-\rSd\|_\infty
\eenn
and hence that
\begin{align*}
	% \Prb\bset{\|\rT\rSdb-\Id\|_\infty > t}
	\I_1(t)
	\smallpad{&\le}
	\Prb\bset{\Mrt\|\rSdb-\rSd\|_\infty > t}
\end{align*}
Choose~$t=\Mrt \ta\sqrt{\log(\px)/\n}$ for a controlled quantity~$\ta>0$ and cite a slight modification of Lemma~\ref{concentration} to find
\begin{align*}
	\I_1\big(\Mrt \ta\sqrt{\log(\px)/\n}\big)
	\smallpad{&\le}
	\Prb\bset{\|\rSdb-\rSd\|_\infty > \ta\sqrt{\log(\px)/\n}} \\
	\smallpad{&\le}
	2\px^{2-\ta^2/(6e^2\sendb^2)}\,,
\end{align*}
where~$\sendb=\Mra^2(2\sgnz^2+\|\rSz\|_\infty/\log2)$.
Substitute the above bound into the first display of the present proof to conclude that
\begin{multline*}
	\Prb\bset{\|\rT\rSdh-\Id\|_\infty > \Mrt \ta\sqrt{\log(\px)/\n} + 3\Mrt\Mra\|\rSzh\|_\infty\gep} \\
	\smallpad\le
	2\px^{2-\ta^2/(6e^2\sendb^2)} + \Prb\set{\maxerrfs>\gep} \,,
\end{multline*}
as claimed.
\end{proof}
}

\begin{lem}[Probability of~$\Tinv(\tf)$]\label{feassgg} %Lemma 4.19
Suppose that (i)~the~$\dz_i$ and~$\dv^j$ satisfy Assumptions~\ref{assu:dzi} and~\ref{assu:subgaussian}, respectively; (ii)~$\rAh$ consists of first-stage Lasso estimates of the~$\ra^j$ with tuning parameters~$\tunej$ chosen according to Definition~\ref{tpss}.
Set
\benn
	\tf \equals \frac{\Mrt}{\sqrt{\n}}\Big(\ta\sqrt{\log\px} + 12\Mra\tolV \maxzzn^{3/2}\scca\sqrt{\log\pz}\Big) \,,
\eenn
where~$\tolV$ is as in Definition~\ref{tpss}.
Then, for~$\n$ sufficiently large,
\benn
	\Prb\bset{\|\rT\rSdh-\Id\|_\infty > \tf}
	\smallpad\le
	2\px^{2-\ta^2/(6e^2\eta^2)} + \prbepallvj\,,
\eenn
where~$\eta=\Mra^2(2\sgnz^2+\|\rSz\|_\infty/\log2)$ and~$\Choeff$ is as defined in Lemma~\ref{lem:epg}.
\end{lem}

\myappendix{
\begin{proof}[Proof of Lemma~\ref{feassgg}]
%Proof of Lemma 4.19
% From Lemma~\ref{feassg} we write
% \begin{multline*}
% 	\Prb\bset{\|\rT\rSdh-\Id\|_\infty > \Mrt \ta\sqrt{\log(\px)/\n} + 3\Mrt\Mra\|\rSzh\|_\infty\gep} \\
% 	\smallpad\le
% 	2\px^{2-\ta^2/(6e^2\sendb^2)} + \Prb\set{\maxerrfs>\gep} \,,
% \end{multline*}
% where~$\sendb=\Mra^2(2\subgauss^2+\|\rSz\|_\infty/\log2)$.
Set~$\gep\equals4\scca\tolV(\maxzzn(\log\pz)/\n)^{1/2}$, plug this choice into the result of Lemma~\ref{feassg} and cite the estimate of Lemma~\ref{lem:eefsg} to find
\begin{multline*}
	\Prb\Big\{\|\rT\rSdh-\Id\|_\infty > \frac{\Mrt}{\sqrt{\n}}\Big(\ta\sqrt{\log\px} + 12\Mra\scca\tolV \maxzzn^{3/2}\sqrt{\log\pz}\Big)\Big\} \\
	\smallpad\le
	2\px^{2-\ta^2/(6e^2\sendb^2)} + \prbepallvj
\end{multline*}
for~$\n$ sufficiently large, as claimed.
\end{proof}
}

\noindent If there exists~$\maxsdz=O(1)$ that satisfies~$\Prb\set{\maxzzn>\maxsdz}=o(1)$,
we may substitute the former quantity into the specification of~$\tf$ in Lemma~\ref{feassgg} to obtain
\begin{align*}
	\Prb\bset{\|\rT\rSdh-\Id\|_\infty>\tf}
	\smallpad\le
	2\px^{2-\ta^2/(6e^2\eta^2)} + 2\pz^{1-\cepg} + \Prb\set{\maxzzn>\maxsdz} \,,
	% \smallpad{\overset{\n\to\infty}{\to}} 0
\end{align*}
which tends to zero under appropriate specification of the controlled quantities~$\ta$ and~$c$.
Note that~$\tf\lesssim\sa\sqrt{\log(\pz)/\n}$;
Condition~\ref{gg1} of Assumption~\ref{gr} then becomes
\be
	\suni\sa^{1-\tuni}(\log\pz)^{1-\frac{\tuni}{2}}/\n^{\frac{1-\tuni}{2}}
	\equals
	o(1) \,,
\ee
and Condition~\ref{gg3} becomes
\be
	\sb\sa^3(\log\pz)^{3/2}/\n +\sb\sa^2\log(\pz)/\sqrt{\n} \equals o(1) \,,
\ee
which is identical to Condition~\ref{gg2}.

Recall that Theorem~\ref{wl} specifies the conditions under which~$\sqrt{\n}(\rbdbj-\rbj)/\wjse$ converges weakly to a~$\Normal(0,1)$ random variable and demonstrates that the limit continues to hold if~$\wjse$ is replaced with an estimator~$\wjseh$ that satisfies~$|\wjseh-\wjse|=o_{\Prb}(1)$.
In practice,~$\wjse$ is not available, hence it is crucial to demonstrate the existence of an estimator~$\wjseh$ that satisfies the foregoing condition.
The following Lemma identifies such an estimator~$\wjseh$ and the conditions under which it is appropriate for use with Theorem~\ref{wl}.

\begin{lem}[Existence of appropriate~$\wjseh$]\label{wjsehcons}
Suppose that~(i) (a) the second-stage noise elements~$\dui$ satisfy~$\E[\dui^2\,|\,\dz_i]=\sdu^2$ and~$\|\du\|_2^2/\n-\sdu^2 = o_{\Prb}(1)$,
(b)~the~$\dz_i$ and~$\dv^j$ satisfy Assumptions~\ref{assu:dzi} and~\ref{assu:subgaussian},
(c)~we have that~$\sb\sa^2\log(\pz)/\n + \sb\sa\sqrt{\log(\pz)/\n}=o(1)$ and~$\maxjpx\maxin\E[\dxij^2]=O(1)$.
Let~$\rAh,\rbla$, and~$\rTh$ be as specified in Assumption~\ref{assum:estimators}.
Define the estimator~$\sduh$ of the second-stage noise level~$\sdu$ by
\be\label{def:sduhhomo}
	\sduh^2 \smallpad{:=} \|\dy-\dX\rbla\|_2^2/\n\,.
\ee
Then~$\sduh-\sdu=o_{\Prb}(1)$.
If, in addition to the conditions in (i), we have (ii) (a)~Condition~\ref{gr:feas} of Assumption~\ref{gr} holds, (b) the sequence of minimal and maximal eigenvalues of~$\rSd$, denoted respectively by~$\eigmin(\rSd)$ and~$\eigmax(\rSd)$, are bounded away from zero and infinity uniformly in $\n$; and (c)~the tolerance~$\tf$ satisfies~$\tf=o(1)$, then~$\wjseh$ defined by
\be\label{def:wjsehhomo}
	\wjseh^2 \smallpad{:=} \sduh^2\rTjjh
\ee
satisfies~$\wjseh-\wjse=o_{\Prb}(1)$.
\end{lem}

\myappendix{
\begin{proof}[Proof of Lemma~\ref{wjsehcons}]
(i) We first show that~$\sduh-\sdu=o_{\Prb}(1)$.
To begin, write
\begin{align*}
	\sduh^2 \equals \|\dy-\dX\rbla\|_2^2/\n
	\aequals \|\du + \dX(\rbla-\rb)\|_2^2/\n \\
	\aequals \|\du\|_2^2/\n + 2\underbrace{\ip{\du,\dX(\rbla-\rb)}/\n}_{\I_1} + \underbrace{\|\dX(\rbla-\rb)\|_2^2/\n}_{\I_2}\,.
\end{align*}
It follows that
\begin{align*}
	(\|\du\|_2^2/\n-\sdu^2) - 2|\I_1| + \I_2
	\smallpad{&\le} \sduh^2-\sdu^2 \\
	\smallpad{&\le} |\sduh^2-\sdu^2|
	\smallpad\le \big|\|\du\|_2^2/\n-\sdu^2\big| + 2|\I_1| + \I_2 \,.
\end{align*}
We claim that $\I_1$ and~$\I_2$ are each~$o_{\Prb}(1)$. It then follows that~$\sduh^2-\sdu^2=o_{\Prb}(1)$,
since~$\|\du\|_2^2/\n-\sdu^2=o_{\Prb}(1)$ by assumption.
To show the claim, note first that
\benn
	\I_1^2 \smallpad\le \big(\|\du\|_2^2/\n\big)\big(\|\dX(\rbla-\rb)\|_2^2/\n\big)
	\equals \big(\|\du\|_2^2/\n\big)\I_2 \,.
\eenn
From the assumption that~$\|\du\|_2^2/\n-\sdu^2=o_{\Prb}(1)$ we infer that~$\|\du\|_2^2/\n=O_{\Prb}(1)$.
Thus, it suffices to show that~$\I_2=o_{\Prb}(1)$.
To this end, note that
\benn
	\|\dX(\rbla-\rb)\|_2\smallpad\le\|\dX\|_{L_2}\|\rbla-\rb\|_1
\eenn
hence~$\I_2 \smallpad\le \|\dX\|_{L_2}^2/\n \big(\|\rbla-\rb\|_1^2\big)$.
Now cite Lemma~\ref{concentration2} to obtain that
\benn
	\|\dX\|_{L_2}^2/\n\equals O_{\Prb}\big(\maxjpx\E[\dxij^2]+a\sqrt{(\log\px)/\n}\big)
\eenn
for a suitably chosen~$a$ of constant order.
It then follows from the growth conditions of the present lemma and the estimate of Lemma~\ref{lem:eessg} for~$\|\rbla-\rb\|_1$ that~$\I_2^2 = O_{\Prb}(1)o_{\Prb}(1) = o_{\Prb}(1)$, as required.
Thus,~$\sduh^2-\sdu^2=o_{\Prb}(1)$.
To show as much for~$\sduh-\sdu$, it suffices to show that~$\Prb\set{|\sduh-\sdu|>\gep}\to0$ for each~$\gep>0$.
We will show that~$\Prb\set{\sduh-\sdu>\gep}\to0$; the matching limit follows from an analogous argument.
Fix~$\gep>0$ and note that
\begin{align*}
	\Prb\set{\sduh-\sdu > \gep} \aequals \Prb\set{\sduh > \sdu+\gep} \\
	\aequals \Prb\set{\sduh^2 > \sdu^2 + 2\sdu\gep + \gep^2} \\
	\aequals \Prb\set{\sduh^2-\sdu^2 > 2\sdu\gep + \gep^2} \smallpad\to 0
\end{align*}
as~$\n\to\infty$ since, by Assumption~\ref{assu:subgaussian},~$\sdu$ is bounded strictly away from zero uniformly in~$\n$.
The previous display and the matching limit for~$\Prb\set{\sdu-\sduh>\gep}$ entail that~$\sduh-\sdu = o_{\Prb}(1)$, as claimed.

(ii) We now show that~$\wjseh^2=\sduh^2\rTjjh$ satisfies~$\wjseh-\wjse=o_{\Prb}(1)$.
We first show as much for~$\wjseh^2-\wjse^2$; the original claim then follows from reasoning analogous to that above.
To this end, note that, since the noise components~$\dui$ are homoscedastic, we have
\benn
	\wjse^2 \equals \E[\ip{\rt_j,\dd_i}^2\dui^2]
	\equals \E\big[\ip{\rt_j,\dd_i}^2\E[[\dui^2\,|\,\dz_i]\big]
	\equals \sdu^2\rTjj \,.
\eenn
Now write
\begin{align*}
  \wjseh^2-\wjse^2
  \aequals
  \sduh^2\rTjjh-\sdu^2\rTjj \\
  \aequals
  \sdu^2(\rTjjh-\rTjj)+(\sduh^2-\sdu^2)(\rTjjh-\rTjj)+(\sdu^2-\sduh^2)\rTjj\,.
\end{align*}
Next, recall that, on the set $\Tinv(\tf)$, it holds due to Lemma~\ref{lem:eesupprc} that~$\|\rth_j-\rt_j\|_\infty \le 2\Mrt\tf$.
Letting $\gep>0$ be arbitrary, it then follows from the previous display that
\begin{align*}
  \Prb\bset{|\wjseh^2-\wjse^2|>2\Mrt\tf(\sdu^2+\gep)+\gep\rTjj}
  \smallpad{&\le}
  \Prb\,\Tinv(\tf)^\comp + \Prb\set{|\sduh^2-\sdu^2|>\gep}\,.
\end{align*}
Now note that, due to the first claim of the present lemma and condition~\ref{gr:feas} of Assumption~\ref{gr}, we have
\begin{align*}
  & \lim_{\n\to\infty}
  \Prb\bset{|\wjseh-\wjse|>2\Mrt\tf(\sdu+\gep)+\gep\rTjj} \\
	& \qquad\smallpad\le
  \lim_{\n\to\infty}
  \Prb\,\Tinv(\tf)^\comp + \Prb\set{|\sduh-\sdu|>\gep} \equals 0 \,.
\end{align*}
Finally, cite Proposition~\ref{rTjjorder} and the present assumption that~$\tf=o(1)$ to find that
\benn
  \limni \Mrt\tf(\sdu+\gep)+\gep\rTjj \equals 0\,.
\eenn
Conclude that~$\wjse^2-\wjseh^2 = o_{\Prb}(1)$ and hence that~$\wjseh-\wjse=o_{\Prb}(1)$, as claimed.
\end{proof}
}

%% file: experiments/experiments.tex
\section{Numerical experiments}\label{s:numerics} \rm %checked. 

In this section, we present a Monte Carlo simulation study of the finite-sample properties of the inferential procedure developed in Section~\ref{inference} using the two-stage Lasso studied in Section~\ref{estimation}.
Our objective is to test this method under a variety of parameter configurations chosen to reflect settings of practical interest.
In Section~\ref{ss:design}, we describe the general scheme according to which the data for each trial are generated and the metrics gathered for each configuration.
In Section~\ref{results}, we enumerate the specific parameter configurations studied and discuss the results.

%% file: experiments/design.tex
\subsection{General experimental design}\label{ss:design}

Each trial contains a data-generation step and an estimation step.
We specify the regression parameters~$\rb$ and~$\rA$ for the data-generation step as follows.
For each configuration, we set the second-stage regression parameter~$\rb$ according to~$\rbc_j = 1$ for~$j \in \Sb$ and~$\rbc_j = 0$ otherwise, where~$\Sb \subset [\px]$ is a random set of~$\sb$ generated by uniformly random draws from~$[\px]$ without replacement.
% and~$b$ is a model parameter that can vary over configurations.
Similarly, we set the first-stage regression parameters~$\ra^j$ for~$j\in[\px]$ according to~$\rac^j_k = 1$ for~$k\in\Saj$ and~$\alpha_{k}=0$ otherwise, where~$\Saj\subset[\pz]$ is a random set of~$\sa$ generated by uniformly random draws from~$[\pz]$ without replacement.
% and~$a$ is a model parameter that can vary over configurations.
We let~$\sb, \sa$ vary over configurations.

Having specified the regression parameters, we then draw $n$ i.i.d. observations~$(\dyi, \dx_i, \dz_i)$ according to
% \begin{align*}
% 	\dz_i \ & \sim \ \Normal_{\pz}(\zeros, \rSz) \,, \\
% 	(u_i, \dv_i) \,|\,\dz_i \ & \sim \ \Normal_{1+\px}(\zeros, \rSuv) \,, \\
% 	x_{ij} \ & = \ \ip{\dz_i, \ra^j} + v_{ij} \,,\\
% 	\dyi \ & = \ \ip{\dx_i, \rb} + u_i \,,
% \end{align*}
\benn
\begin{split}
\dz_i \ & \sim \ \Normal_{\pz}(\zeros, \rSz) \,, \\
(u_i, \dv_i) \,|\,\dz_i \ & \sim \ \Normal_{1+\px}(\zeros, \rSuv)
\end{split}
\qquad\qquad
\begin{split}
	x_{ij} \ & = \ \ip{\dz_i, \ra^j} + v_{ij} \,,\\
	\dyi \ & = \ \ip{\dx_i, \rb} + u_i \,,
\end{split}
\eenn
where~$\n, \px, \pz, \rb, \{\ra^j\}_{j\in[\px]},$ and the structure of~$\rSz$ vary amongst configurations.
For all configurations, we set
\benn
	\rSuv \equals \begin{pmatrix} \sdu & \sduv^\top \\ \sduv & \sdv \Id \end{pmatrix} \,,
\eenn
where~$\sdu, \sdv = \sqrt{.7}$ are held fixed and~$\sduv\equals(\sd_{uv^{1}}, \ldots, \sd_{uv^{\px}})$ is given as follows.
For each configuration, we set one~$\sd_{uv^j}$ chosen at random equal to~.5, nine~$\sd_{uv^j}$ chosen at random equal to~.25, and the remaining~$\sd_{uv^j}$ equal to~.05.
The present covariance structure for the noise reflects the constraint that~$\rSuv$ be positive-definite;
our choices of~$\sd_{uv^j}$ are an attempt to balance this requirement with the goal of studying non-trivial correlations between the first- and second-stage noise elements.

We consider two forms for the covariance matrix~$\rSz$.
The first form is a Toeplitz~(TZ) structure given by
\benn
	\left.\rSztz\right|_{jk} \equals\rho^{|j-k|}\,, \qquad j,k \in[\pz] \,, \qquad \rho = 0.8 \,.
\eenn
The second is a circulant-symmetric~(CS) structure given for~$j \leq k$ by
\benn
	\threecase{\left.\rSzcs\right|_{jk} \equals}{1}{$k = j$\,,}{0.1}{$k\in\{j+1,\ldots,j+5\} \cup \{j+\pz-5,\ldots,j+\pz-1\}$\,,}{0}{otherwise\,.}
\eenn
Within a configuration study, the random quantities~$\dz_i, u_i,$ and~$\dv_i$ are re-drawn for each trial; the quantities~$\rb, \{\ra^j\}_{j=1}^{\px}, \rSz, \sduv, \n, \px, \pz$ are held fixed.

For the estimation step of each trial, we compute the first- and second-stage Lasso as defined in Section~\ref{ss:two-stage} using the \texttt{glmnet} package \cite{glmnet}.
Tuning parameters~$\tune$ for the Lasso estimators are selected by 10-fold cross-validation over a grid~$\{\tune_\ell\}_{\ell=1}^L$, where~$L=100$, $\tune_L = .01\tune_1$, and~$\tune_1$ is the least quantity for which the respective Lasso estimator is identically 0.
The tuning parameter $r_{\bm{\beta}}$ is chosen by similar cross-validation procedure.

The rows~$\rth_j$ of~$\rTh$ are obtained as solutions to the respective Program~\ref{program} with tuning parameter~$\mu_j$ chosen as follows.
For each~$j\in[\px]$, we set~$\tf_j := \kappa\times\inf_{\rt\in\R^{\px}}\|\rSdh\rt-\basis_j\|_\infty$,
where~$\basis_j$ denotes the $j^{\text{th}}$ canonical basis vector in~$\px$ dimensions and $\kappa > 1$ is chosen at our discretion.
Note that, under this choice of~$\mu_j$, the respective Program~\ref{program} is guaranteed feasible.
The factor~$\kappa$ is chosen to balance the performance of~$\rTh$ as a surrogate inverse for~$\rSdh$, for which a smaller~$\kappa$ is desirable, with the size of the objective function~$\|\rt\|_1$, for which a larger~$\kappa$ is desirable.
The following results were obtained under~$\kappa=1.2$.
To obtain the infimum, we cast~$\mathrm{minimize}_{\rt\in\R^{\px}}\|\rSdh\rt-\basis_j\|_\infty$ as a linear programming problem, which we solve using MOSEK optimization software \cite{mosek}.

In a given trial $t,t=1,\dots,T$, we set~$\tau = 0.05$ and compute the respective~$100(1-\tau)\%$ confidence interval
\benn
	\hat{\mathcal{I}}_{\tau,t,j} \equals\big[\tilde{\beta}_{t,j} - z_{\tau} \SEh(\tilde{\beta}_{t,j}), \, \rbdbj + z_{\tau}\SEh(\tilde{\beta}_{t,j})\big]  \,,
\eenn
for each component~$\tilde{\beta}_{t,j}$ of~$\tilde{\bm\beta}_{t}$, where~$z_{\tau} = \Phi^{-1}(1-\tau/2)$ and
\begin{equation*}
	\SEh(\rbdbj)^2 \equals \En[(\dyi - \rbla\dX)^2\ip{\rth_j, \ddh_i}^2] \,.
\end{equation*}
% To assess the finite-sample performance of the procedure developed in Section~\ref{inference},
For each configuration of~$\n,\, \px,\, \pz,\, \sb,\, \sa,\, \rSz$, we generate $T = 100$ trials and calculate the average coverage~$\cvgh$ for the~95\% confidence intervals $\hat{\mathcal{I}}_{\tau,j}$ about components of $\rbdb$ and the average interval length $\lenh$ given by
\benn
\begin{split}
	\cvgh \equals\frac{1}{p_{\bm{x}}}\sum_{j=1}^{p_{\bm{x}}}\frac{1}{T}\sum_{t=1}^{T}1\left\{ \beta_{t,j}\in\hat{\mathcal{I}}_{\tau,t,j}\right\},
\end{split}\qquad
\begin{split}
	\lenh \equals\frac{1}{p_{\bm{x}}}\sum_{j=1}^{p_{\bm{x}}}\frac{1}{T}\sum_{t=1}^{T}\mbox{len}(\hat{\mathcal{I}}_{\tau,t,j}).
\end{split}
\eenn
For each configuration, we also provide the average mean squared error of the second-stage Lasso estimator given by
$\mseh=\frac{1}{T} \sum_{t=1}^T \frac{1}{\px}\sum_{j=1}^{\px} (\tilde{\beta}_{t,j}-\beta_{t,j})^2 \,$.
We present the results for the study described above in Table 1 in Section~\ref{results}.

%% file: experiments/results.tex
\subsection{Specifications and results}\label{results}
We conduct simulations according the design described in Section~\ref{ss:design} for all configurations belonging to
\benn
\underbrace{\begin{pmatrix*}[r]
		(100,125,150) \\
		(200, 250,275) \\
		(300,400,500) \\
		(500,600,700)
\end{pmatrix*}}_{(n\,,\px\,,\pz)} \times
\underbrace{\begin{pmatrix*}[r]
		(3, 5) \\
		(5,10) \\
		(10, 15)
\end{pmatrix*}}_{(\sb\,,\sa)} \times
\underbrace{\begin{pmatrix*}[r]
	\rSzcs \\[.4em] \rSztz
\end{pmatrix*}}_{\rSz} \,.
\eenn
The results, which are presented in Table~\ref{table I}, show that our estimator achieves close to nominal coverage under a variety of configurations.
% We see that the coverage over the parameter~$\rb$ as a whole and over the inactive components~$\rb_{\sb^{\comp}}$ tends to be closer to the nominal level than the coverage over the active components~$\rb_{\sb}$.
We also see that arguably the greatest determinant of coverage is the relative magnitude of~$\px$ and~$\pz$ to the size of the active set~$\sb$. As the latter grows, coverage diverges from the nominal level.
This phenomenon is expected and has been found in ordinary linear regression models as well \cite{vdGetal14, JM14}, since the bounds for the estimation error of the Lasso is proportional to the size of the active set. Nevertheless, the performance improves significantly when we increase the sample size.
Finally, we observe that the covariance structure of the instrumental variables $\dz_i$ has a strong influence on coverage: the Toeplitz structure features greater correlation among the instrumental variables in general, and this is reflected in coverage that tends to be farther from the nominal level than in the case of the circulant-symmetric covariance structure.
%On the whole, the finite-sample performance of our estimator is comparable to those of~\cite{BvdG11,JM14}.
These results suggest that our proposed method of inference for the low-dimensional components of a high-dimensional regression vector is relevant to practical scenarios that may exhibit non-trivial degrees of correlation between the noise components and nontrivial correlation among the instrumental variables.

\begin{table}[ht]
\centering
\caption{Simulation results}\label{table I}
% \makebox[\textwidth][c]{%
\begin{tabular}{llrrrrrr}
	\toprule
& & \multicolumn{3}{c}{$\rSzcs$} & \multicolumn{3}{c}{$\rSztz$}	\\
\cmidrule(lr){3-5}\cmidrule(lr){6-8}
$(\n,\, \px,\, \pz)$ & $(\sb,\saj)$ & $\cvgh$ & $\lenh$ & $\mathrm{MSE}(\rbh)$
& $\cvgh$ & $\lenh$ & $\mathrm{MSE}(\rbh)$\\
  \midrule
 (100, 125, 150) & (3, 5) & 0.942 & 0.225 & 0.004 & 0.895 & 0.201 & 0.005 \\
 & (5,10) & 0.941 &0.211 & 0.004 & 0.672 & 0.212 & 0.014\\
  & (10, 15) & 0.930 & 0.190 & 0.003 & 0.545 & 0.219 & 0.030 \\

 \addlinespace[.5em]
 (200, 250, 275) & (3, 5) & 0.947 & 0.157 & 0.002 & 0.942 & 0.140 & $\le 0.001$ \\
 & (5,10) & 0.941 & 0.171 & 0.002 & 0.673 & 0.192 & 0.011\\
 & (10, 15) & 0.930 & 0.190 & 0.003 & 0.545 & 0.219 & 0.030 \\

  \addlinespace[.5em]
 (300, 400, 500) & (3, 5) & 0.947 & 0.094 & $\le 0.001$ & 0.952 & 0.092 & $\le 0.001$ \\
 & (5,10) & 0.955 & 0.085 &$\le 0.001$ &0.945 & 0.082 & $\le 0.001$\\
  & (10, 15) & 0.961 & 0.067 & $\le 0.001$ & 0.927 & 0.064 & $\le 0.001$ \\

 \addlinespace[.5em]
 (500, 600, 700) & (3, 5) & 0.947 & 0.094 & $\le 0.001$ & 0.952 & 0.092 & $\le 0.001$ \\
 & (5,10) & 0.951 & 0.082 & $\le 0.001$ & 0.950 & 0.088 & $\le 0.001$\\
 & (10, 15) & 0.961 & 0.067 & $\le 0.001$ & 0.927 & 0.064 & $\le 0.001$ \\

   \bottomrule
\end{tabular}
% }
\end{table}

%
% \begin{table}[ht]
% \centering
% \caption{Simulation results for $\bm{\kappa}^{\top}\bm{\beta}$}\label{table II}
% % \makebox[\textwidth][c]{%
% \begin{tabular}{llrrrrrr}
% 	\toprule
% & & \multicolumn{3}{c}{$\rSzcs$} & \multicolumn{3}{c}{$\rSztz$}	\\
% \cmidrule(lr){3-5}\cmidrule(lr){6-8}
% $(\n,\, \px,\, \pz)$ & $(\sb,\saj)$ & $\cvgh$ & $\lenh$ & $\mathrm{MSE}(\bm{\beta})$
% & $\cvgh$ & $\lenh$ & $\mathrm{MSE}(\bm{\beta})$\\
%   \midrule
%  (100, 125, 150) & (3, 5) & 0.935 & 0.231 &0.005 & 0.946 & 0.031 & 0.005\\
%  & (5,10) & 0.931 & 0.211 &0.005 & 0.672 & 0.212 & 0.007\\
%   & (10, 15) & 0.921 &0.321 &0.005 &0.673 &0.232 & 0.030\\
%
%  \addlinespace[.5em]
%  (200, 250, 275)  & (3,5) & 0.946 &0.323 &0.005 &0.824 &0.463 & $\le 0.001$\\
%  & (5,10) & 0.946 &0.454 &0.004 &0.843 &0.452 & 0.010\\
%  & (10, 15)  & 0.942 &0.432 &0.004 &0.799 &0.523 & 0.023\\
%
%   \addlinespace[.5em]
%  (300, 400, 500) & (3, 5) & 0.947 & 0.094 & 0.002 & 0.942 & 0.092 &$\le 0.001$\\
%  & (5,10) & 0.953 &0.097 & 0.002  & 0.948 & 0.092 &$\le 0.001$\\
%   & (10, 15) & 0.948 & 0.098 & 0.002 & 0.951 & 0.094 &$\le 0.001$\\
%
%  \addlinespace[.5em]
%  (500, 600, 700) & (3, 5) & 0.951 & 0.081& 0.002  & 0.950 & 0.087 &$\le 0.001$\\
%  & (5,10) & 0.951 & 0.082 & 0.002 & 0.951 & 0.088 &$\le 0.001$\\
%  & (10, 15) & 0.950 & 0.082 & $\le 0.001$  & 0.951 &0.089 &$\le 0.001$\\
%
%    \bottomrule
% \end{tabular}
% % }
% \end{table}

%% file: conclusion/conclusion.tex
\section{Conclusion}\label{sec:conclusion}
In this paper, we propose inference methods for the components of a high-dimensional instrumental variables regression parameter despite possible endogeneity of each regressor.
We allow both the number of instruments and the number of regressors to be greater than the sample size.
We construct asymptotically valid confidence intervals for the components of the second-stage regression coefficients.
Though our estimator is not a nonlinear generalized method of moments (GMM) estimator \cite{H82}, we expect that our results can be extended to that more general setting.

Our Sections~\ref{s:two-stage} and~\ref{inference} comprise a general pipeline for estimation and inference, while the remainder is then an exemplification with a Lasso approach.
Therefore, it would be interesting to use our pipeline with other regularized estimators as well (one could also use different estimators for the two different stages).
Candidates include, for example, SCAD~\cite{Fan&Li2001} and MCP~\cite{zhang2010nearly}.

We finally refer to our software on \href{https://github.com/LedererLab/HDIV}{github.com/LedererLab/HDIV}.

%% file: proofs/proofs.tex
% \section{Supplementary materials}\label{proofs}

\proofs

\section{Technical lemmas}\label{technical}

The probabilistic guarantees for Lasso estimation performance require that the tuning parameter dominate the empirical process term.
In our consideration of both the first- and second- stage Lasso estimators, we encounter a number of such terms, each of the form~$\maxzn{\gerr}$ for various noise vectors~$\gerr$.
As such, we formulate the following lemma, which is used throughout our consideration of the sub-Gaussian error regime, in terms of a generic sub-Gaussian vector~$\gauss$ with i.i.d. components.
The lemma itself is a standard application of basic concentration results for sub-Gaussian random variables.
In the subsequent corollaries, we derive bounds for various empirical process terms by taking~$\gauss$ to be, for instance,~$\du$ and~$\dv^j$ for $j\in[\px]$.
Such bounds are key ingredients of the results presented in Sections~\ref{second}.

\begin{lem}[Control of $\maxzn{\gauss}$]\label{lem:epg}
Let~$\gauss\,|\,\dZ$ be sub-Gaussian and let~$\|\gauss\|_\sgn=\subgauss$.
Set~$\tp=\tol(\maxzzn(\log\pz)/\n)^{1/2}$, where~$\tol>0$ is a controlled quantity.
Then
\benn
  \Prb\bset{4\maxzn{\gauss} > \tp} \smallpad{\le} e\pz^{1-\tol^2\Choeff/\subgauss^2}\,,
\eenn
where~$\Choeff$ is an absolute constant.
\end{lem}

%%%%%%%%%%
% PROOF
\begin{proof}[Proof of Lemma \ref{lem:epg}]
For a given~$j\in[\px]$, cite \cite[Proposition~5.10]{V12} to observe that for all~$t>0$ and all~$\veca=(a_1,\ldots,a_\n)\in\R^\n$ it holds that
\begin{align*}
  \Prb\left\{ \left|\sumin a_i\gaussi\right| > t \,|\, \dZ  \right\}
  \smallpad{&\le}
  e\cdot\exp\left(-\frac{Ct^2}{\max_{i\in[\n]} \|\gaussi\|_\sgn^2\|\veca\|_2^2}\right) \\
  \smallpad{&\le}
  e\cdot\exp\left(-\frac{Ct^2}{\|\gauss\|_\sgn^2\|\veca\|_2^2}\right)\,,
\end{align*}
where~$C$ is an absolute constant.
Take~$a_i=\dzij/\n$, and observe that~$\|\veca\|_2^2= \|\dz^j/\sqrt{\n}\|_2^2/\n$. The display above then yields
\begin{align*}
  \Prb\left\{ \left|\fon\sumin\dzij\gaussi\right| > t \,|\, \dZ\right\}
  \smallpad{\le}
  e\cdot\exp\left(-\frac{Ct^2\n}{\subgauss^2\|\dz^j/\sqrt{\n}\|_2^2} \right) \,.
\end{align*}
Now choose~$t=\foq\tol(\maxzzn(\log\pz)/\n)^{1/2}=\foq\tp$ to find that
\begin{align*}
  \Prb\left\{ 4\left|\fon\sumin\dzij\gaussi\right| > \tp \,|\, \dZ \right\}
  \smallpad{&\le}
  e\cdot\exp\left(-\frac{\tol^2C\maxzzn\log\pz}{16\subgauss^2\|\dz^j/\sqrt{\n}\|_2^2} \right) \\
  \smallpad{&\le}
  e\pz^{-\tol^2C\maxzzn/(16\subgauss^2\|\dz^j/\sqrt{\n}\|_2^2)} \\
  \smallpad{&\le}
  e\pz^{-\tol^2C/(16\subgauss^2)}
  \smallpad{=:}
  e\pz^{-\tol^2\Choeff/\subgauss^2}\,.
\end{align*}
Take the union bound over~$j\in[\px]$ and use iterated expectations to conclude that
\begin{align*}
  \Prb\bset{4\maxzn{\gauss}>\tp}
  \aequals \E\big[\Prb\bset{4\maxzn{\gauss}>\tp\,|\,\dZ}\big] \\
  \smallpad{&\le} \E[e\pz^{1-\tol^2\Choeff/\subgauss^2}]
  \equals
  e\pz^{1-\tol^2\Choeff/\subgauss^2}\,,
\end{align*}
% \begin{multline*}
%   \Prb\bset{4\maxzn{\gauss}>(\maxzzn(\log\pz)/\n)^{1/2}} \\
%   \equals \E\big[\Prb\bset{4\maxzn{\gauss}>(\maxzzn(\log\pz)/\n)^{1/2}\,|\,\dZ}\big] \smallpad{\le} \E[e\pz^{1-\Choeff/\subgauss}]
%   \equals
%   e\pz^{1-\Choeff/\subgauss}\,,
% \end{multline*}
as claimed.

\end{proof}

%%%%%%%%%%%%%%%%%%%%%%%%%%%%%%%%%%%%%%%%%%%%%%%%%%%%%%%%%%%%%%%%%%%%%%%%%%%%%%%
%%%%%%%%%%%%%%%%%%%%%%%%%%%%%%%%%%%%%%%%%%%%%%%%%%%%%%%%%%%%%%%%%%%%%%%%%%%%%%%

\noindent Lemmas~\ref{lem:epu} and~\ref{lem:epV} follow from Lemma \ref{lem:epg} and are required throughout the results for the first- and second-stage estimation errors.

\begin{lem}[Control of $\maxzn{\du}$]\label{lem:epu}
Suppose that $\du$ satisfies Assumption~\ref{assu:subgaussian}.
For~$\lambda>0$, define the set
\be\label{def:Tu}
	\Tu(\lambda)
	\smallpad{:=}
	\{4\maxzn{\du} \leq \lambda\} \,.
\ee
Set~$\tuneu=\tolu(\maxzzn(\log\pz)/\n)^{1/2}$, where~$\tolu>0$ is a controlled quantity.
Then,
\benn
  \Prb\,\Tu(\tuneu)^{\comp}
  \smallpad\le
  e\pz^{1-\tolu^2\Choeff/\sgnu^2} \,,
\eenn
where~$\Choeff$ is as in Lemma~\ref{lem:epg} and~$\sgnu=\|\du\|_\sgn$.
\end{lem}

%%%%%%%%%%
% PROOF

\begin{proof}[Proof of Lemma~\ref{lem:epu}]
This follows immediately from Lemma~\ref{lem:epg}.
\end{proof}

%%%%%%%%%%%%%%%%%%%%%%%%%%%%%%%%%%%%%%%%%%%%%%%%%%%%%%%%%%%%%%%%%%%%%%%%%%%%%%
%%%%%%%%%%%%%%%%%%%%%%%%%%%%%%%%%%%%%%%%%%%%%%%%%%%%%%%%%%%%%%%%%%%%%%%%%%%%%%

\begin{lem}[Simultaneous control of $\maxzn{\dv^j}$]\label{lem:epV}
Suppose that $\dV$ satisfies Assumption~\ref{assu:subgaussian}.
For each $j\in[\px]$ and $\lambda_j>0$, define the set
\be\label{def:Tvj}
	\Tvj(\lambda_j)
	\smallpad{:=}
	\bset{4\|\dZ^\top\dv^j/n\|_\infty \le \lambda_j} \,.
	% \smallpad{=}
	% \big\{\|\dZ^\top\dV/n\|_\infty \le \maxjpx\tunej\big\}\,.
\ee
Let~$\tunevec:=(\lambda_1,\ldots,\lambda_{\px})$.
Define the set
\be\label{def:TV}
  \TV(\tunevec)
  \smallpad{:=}
  \bigcap_{j=1}^{\px}\Tvj(\lambda_j) \,.
\ee
If~$\lambda_j=\lambda$ for each~$j\in[\px]$ and a given~$\lambda>0$, we abuse notation and write
\benn
  \TV(\lambda)
  \equals
  \bigcap_{j=1}^{\px}\Tvj(\lambda)
  \equals
  \bset{\maxzn{\dV}\le\lambda}
\eenn
For each~$j\in[\px]$, set~$\lambda_j=\tunejverb$, where each~$\tolj>0$ is a controlled quantity.
Set~$\tuneV=\tuneVverb$, where~$\tolV>0$ is a controlled quantity.
Then,
\benn
  \Prb\TV(\tunevec)^{\comp}
  \smallpad\le
  \prbepallvj
\eenn
and
\benn
  \Prb\TV(\tuneV)^{\comp}
  \smallpad\le
  \prbepV \,,
\eenn
where~$\Choeff$ is as in Lemma~\ref{lem:epg} and~$\sgnV=\maxjpx\|\dv^j\|_\sgn$.
\end{lem}

\begin{proof}[Proof of Lemma~\ref{lem:epV}]
We show the first claim.
Note that, for~$\tunevec$ as specified in the statement of the present lemma,
\begin{align*}
	\TV(\tunevec)^\comp
	\equals
	\bigcup_{j\in[\px]} \bset{4\|\dZ^\top\dv^j/\n\|_\infty > \lambda_j} \,.
\end{align*}
Hence
\begin{align*}
  \Prb\,\TV(\tunevec)^{\comp}
  \smallpad{&\le}
  \sum_{j=1}^{\px}\Prb\bset{4\maxzn{\dv^j} > \lambda_j} \\
  \smallpad{&\le}
  \sum_{j=1}^{\px} e\pz^{1-\tolj^2\Choeff/\|\dv^j\|_\sgn^2} \\
  \smallpad{&\le}
  e\px\pz^{1-\Choeff\min_{j\in[\px]}\{\tolj^2/\|\dv^j\|^2_\sgn\}} \\
  \smallpad{&\le}
  e\pz^{2-\Choeff\min_{j\in[\px]}\{\tolj^2/\|\dv^j\|^2_\sgn\}} \,,
\end{align*}
where the inference to the second line above follows from Lemma~\ref{lem:epg} and the inference to the last line follows from the assumption under the presently studied regime that~$\px\le\pz$.

To show the second claim, set~$\tolj=\tolV>0$ for each~$j\in[\px]$ and then note that
\begin{equation*}
  \min_{j\in[\px]}\{\tolj^2/\|\dv^j\|^2_\sgn\} \equals \tolV/\maxjpx\|\dv^j\|_\sgn^2 \,.
\end{equation*}
\end{proof}

\begin{prop}[Regularity of~$\rTjj$]\label{rTjjorder}
Suppose that the minimal and maximal eigenvalues~$\eigmin(\rSd)$ and~$\eigmin(\rSd)$ of~$\rSd$ satisfy
\benn
  0 \smallpad\le \Cmin \smallpad\le \eigmin(\rSd) \smallpad\le \eigmax(\rSd) \smallpad< \Cmax < \infty
\eenn
for universal constants~$\Cmin, \Cmax$.
Then
\benn
  \Cmin/\Cmax \smallpad\le \rTjj \smallpad\le \Cmax/\Cmin\,.
\eenn
\end{prop}

\begin{proof}
First write
\benn
  \rTjj \equals \rt_j^\top\rSd\rt_j \smallpad\ge \eigmin(\rSd)\|\rt_j\|_2^2
\eenn
and note that
\benn
  \|\rt_j\|_2^2 \equals \|\rT\basis_j\|_2^2 \smallpad\ge \eigmin(\rT) \equals 1/\eigmax(\rSd) \,,
\eenn
where~$\basis_j$ denotes the~$j^{\text{th}}$ canonical basis vector in~$\px$ dimensions.
It follows that~$\rTjj\geq\eigmin(\rSd)/\eigmax(\rSd)$, as claimed.
The upper bound follows by analogous reasoning.
\end{proof}

\noindent The following lemma is required for Lemmas~\ref{satisfied},~\ref{ccssg}, and~\ref{feassgg}.
The proof is similar to that of \cite[Theorem 2.4.(b)]{JM14}; both use the Bernstein-type inequality of \cite[Proposition 5.16]{V12} and union bounds to derive concentration results for a sup-norm of interest.

\begin{lem}[Concentration of~$\|\rSdb-\rSd\|_\infty$]\label{concentration}
Suppose that the~$\dz_i$ satisfy Assumption~\ref{assu:dzi}.
Then, for~$\n$ sufficiently large,
\benn
  \Prb\bset{\|\rSdb-\rSd\|_\infty > a\sqrt{(\log(\px\vee\n))/\n}}
  \smallpad\le
  2(\px\vee\n)^{2-a^2/(6e^2\sendb^2)}\,,
\eenn
where~$a$ is a controlled quantity, $\sendb:=\Mra^2(2\sgnz^2+\|\rSz\|_\infty/\log2)$, and~$\sgnz=\|\dz_i\|_\sgn$.
\end{lem}

\begin{proof}[Proof of Lemma~\ref{concentration}]\label{proof:concentration}
Note that~$\rSdb-\rSd=\fon\sumin\dd_i\dd_i^\top-\rSd$ and hence that~$[\rSdb-\rSd]_{jk}=\fon\sumin\ddij\ddik-\sdjk$, where $\sdjk:=\rSdjk$.
Now, for any two random variables~$X$ and~$Y$, it holds that~$\|XY\|_\sen\le2\|X\|_\sgn\|Y\|_\sgn$.
Further, if~$\mu\in\R$ is a constant, then~$\|\mu\|_\sen=|\mu|/\log2$.
Thus,
\benn
  \|\ddij\ddik-\sdjk\|_\sen
  \smallpad\le
  2\|\ddij\|_\sgn\|\ddik\|_\sgn+|\sdjk|/\log2 \,.
  % \smallpad\le\|\dx_i\|_\sgn+|\sdjk|/\sqrt{2}\,.
\eenn
Note that~$\ddij=\dz_i^\top\ra^j=\sum_{k=1}^{\pz}\rajk\dzik$ and that
\begin{align*}
  % \|\dxij\|_\sgn \equals \|\dz_i^\top\ra^j\|_\sgn
  % \equals
  \|\sum_{k=1}^{\pz}\rajk\dzik\|_\sgn
  \smallpad{\le} \sum_{k=1}^{\pz}\|\rajk\dzik\|_\sgn
  \smallpad{\le} \|\dz_i\|_\sgn\sum_{k=1}^{\pz} |\rajk|
  \smallpad{\le} \sgnz\Mra
\end{align*}
and similarly for~$\|\ddik\|_\sgn$.
Thus
\benn
  \|\ddij\ddik-\sdjk\|_\sen
  \smallpad\le
  2\sgnz^2\Mra^2+|\sdjk|/\log2
  \smallpad\le
  \Mra^2(2\sgnz^2+\|\rSz\|_\infty/\log2)
  \smallpad{=:} \sendb \,.
  % \smallpad\le\|\dx_i\|_\sgn+|\sdjk|/\sqrt{2}\,.
\eenn
Now apply the Bernstein-type inequality of \cite[Proposition 5.16]{V12}
to conclude that
\benn
  \Prb\Big\{\fon\Big|\sumin\ddij\ddik-\sdjk\Big|>t\Big\}
  \smallpad{\le}
  2\exp\Big(-\frac{\n}{6}\min\Big\{\big(\frac{t}{e\sendb}\big)^2,\frac{t}{e\sendb}\Big\}\Big) \,.
\eenn
Choose~$t=a\sqrt{(\log(\px\vee\n))/\n}$.
If~$\n\ge(a/(e\sendb))^2\log(\px\vee\n)$, then
\benn
  \Prb\Big\{\fon\Big|\sumin\ddij\ddik-\sdjk\Big|>a\sqrt{(\log(\px\vee\n))/\n}\Big\}
  \smallpad{\le}
  2(\px\vee\n)^{-a^2/(6e^2\sendb^2)} \,,
\eenn
and hence
\benn
  \Prb\bset{\|\rSdb-\rSd\|_\infty > a\sqrt{(\log(\px\vee\n))/\n}}
  \smallpad\le
  2(\px\vee\n)^{2-a^2/(6e^2\sendb^2)}\,,
\eenn
which follows from taking the union bound over~$j,k\in[\px]$.
\end{proof}

\noindent The following lemma is required for Lemma~\ref{wjsehcons}.

\begin{lem}[Concentration of~$\|\dX\|_{L_2}$]\label{concentration2}
Suppose that the~$\dz_i$ and~$\dv^j$ satisfy Assumptions~\ref{assu:dzi} and~\ref{assu:subgaussian}, respectively.
Then, for
\benn
  \n \smallpad> \maxjpx \maxin a^2\log\px/\|\dxij^2-\E[\dxij^2]\|_\sen^2
\eenn
where~$a>0$ is a controlled quantity, it holds that
\benn
  \Prb\bset{\|\dX\|_{L_2}^2/\n > \maxjpx\maxin\E[\dxij^2] + a\sqrt{(\log\px)/\n}} \\
  \smallpad{\le}
  2\px^{1-\Cbernstein a^2/\senxscb^2} \,,
\eenn
where~$\senxscb=\maxjpx \big\{(\Mra\sgnz + \sgnvj)^2 + \maxin\E[\dxij^2]/\log2\big\}$.
\end{lem}

\begin{proof}[Proof of Lemma~\ref{concentration2}]
Fix~$j\in[\px]$ and write
\begin{align*}
  \Prb\bset{\|\dx^j\|_2^2/\n > \maxkpx\E[\dxik^2] + t}
  \smallpad{&\le}
  \Prb\left\{\fon\sumin\dxij^2 > \E[\dxij^2] + t\right\} \\
  \smallpad{&\le}
  \Prb\left\{\left|\fon\sumin \dxij^2 - \E[\dxij^2]\right| > t \right\} \\
  \smallpad{&\le}
  2\exp\left(-\Cbernstein\n\min\left(\frac{t^2}{\subgauss_{j}^2}, \frac{t}{\subgauss_{j}}\right)\right) \,,
\end{align*}
where~$\Cbernstein$ is an absolute constant,~$\subgauss_{j}=\maxin\|\dxij^2-\E[\dxij^2]\|_{\sen}$ and the final line follows from~\cite[Corollary 5.17]{V12}.
Set~$t=a\sqrt{(\log\px)/\n}$ for a controlled quantity~$a>0$ and note that if~$\n>a^2\log\px/\subgauss_{j}^2$ then~$t/\subgauss_{j}<1$ and hence
\benn
  \Prb\bset{\|\dx^j\|_2^2/\n > \maxkpx\E[\dxik^2] + t}
  \smallpad\le
  2\px^{-\Cbernstein a^2/\subgauss_{j}^2}
  \smallpad\le
  2\px^{-\Cbernstein a^2/\subgauss^2}
\eenn
for each~$j\in[\px]$, where~$\subgauss=\maxjpx\subgauss_{j}$.
To deduce the original claim, take the union bound over such~$j$ and observe that
\begin{align*}
  \subgauss \equals \maxjpx \maxin\|\dxij^2-\E[\dxij^2]\|_{\sen}
  \smallpad{&\le}
  \maxjpx \maxin \big\{ \|\dxij^2\|_\sen + \E[\dxij^2]/\log2 \big\} \\
  \aequals
  \maxjpx \maxin \big\{ \|\dxij\|_\sgn^2 + \E[\dxij^2]/\log2 \big\} \\
  \smallpad{&\le}
  \maxjpx \maxin \big\{ (\|\ddij\|_\sgn+\|\dvij\|_\sgn)^2 \\
  &\qquad + \E[\dxij^2]/\log2 \big\} \\
  \smallpad{&\le}
  \maxjpx \big\{(\Mra\sgnz + \sgnvj)^2  \\
  &\qquad + \maxin\E[\dxij^2]/\log2\big\} \\
  \aequals \senxscb\,,
\end{align*}
where we infer that~$\|\ddij\|_\sgn\le\Mra\sgnz$ as in the proof of Lemma~\ref{concentration}.
\end{proof}

%% file: script.bbl
\begin{thebibliography}{56}
% BibTex style file: imsart-number.bst, 2013-01-28
% Default style options (sort=1,type=number).
% Used options (sort=1,type=number).

\bibitem{A74}
\begin{barticle}[author]
\bauthor{\bsnm{Amemiya},~\bfnm{T.}\binits{T.}}
(\byear{1974}).
\btitle{The nonlinear two-stage least-squares estimator}.
\bjournal{J. Econometrics}
\bvolume{2}
\bpages{105--110}.
\end{barticle}
\endbibitem

\bibitem{A77}
\begin{barticle}[author]
\bauthor{\bsnm{Amemiya},~\bfnm{T.}\binits{T.}}
(\byear{1977}).
\btitle{The maximum likelihood and the nonlinear three-stage least squares
  estimator in the general nonlinear simultaneous equation model}.
\bjournal{Econometrica}
\bvolume{45}
\bpages{955--968}.
\end{barticle}
\endbibitem

\bibitem{AP09}
\begin{bbook}[author]
\bauthor{\bsnm{Angrist},~\bfnm{J.}\binits{J.}} \AND
  \bauthor{\bsnm{Pischke},~\bfnm{J-S.}\binits{J.-S.}}
(\byear{2009}).
\btitle{Mostly harmless econometrics: An empiricist's companion}.
\bpublisher{Princeton University Press}.
\end{bbook}
\endbibitem

\bibitem{mosek}
\begin{bmanual}[author]
\bauthor{\bsnm{{MOSEK ApS}}}
(\byear{2017}).
\btitle{MOSEK Rmosek Package 8.1.0.34}.
\end{bmanual}
\endbibitem

\bibitem{BCCH12}
\begin{barticle}[author]
\bauthor{\bsnm{Belloni},~\bfnm{A.}\binits{A.}},
  \bauthor{\bsnm{Chen},~\bfnm{D.}\binits{D.}},
  \bauthor{\bsnm{Chernozhukov},~\bfnm{V.}\binits{V.}} \AND
  \bauthor{\bsnm{Hansen},~\bfnm{C.}\binits{C.}}
(\byear{2012}).
\btitle{Sparse models and methods for optimal instruments with an application
  to eminent domain}.
\bjournal{Econometrica}
\bvolume{80}
\bpages{2369--2429}.
\end{barticle}
\endbibitem

\bibitem{BC13}
\begin{barticle}[author]
\bauthor{\bsnm{Belloni},~\bfnm{A.}\binits{A.}} \AND
  \bauthor{\bsnm{Chernozhukov},~\bfnm{V.}\binits{V.}}
(\byear{2013}).
\btitle{Least squares after model selection in high-dimensional sparse models}.
\bjournal{Bernoulli}
\bvolume{19}
\bpages{521--547}.
\end{barticle}
\endbibitem

\bibitem{BCH11}
\begin{bincollection}[author]
\bauthor{\bsnm{Belloni},~\bfnm{A.}\binits{A.}},
  \bauthor{\bsnm{Chernozhukov},~\bfnm{V.}\binits{V.}} \AND
  \bauthor{\bsnm{Hansen},~\bfnm{C.}\binits{C.}}
(\byear{2011}).
\btitle{Inference for high-dimensional sparse econometric models}.
In \bbooktitle{Advances in Economics and Econometrics: Tenth World Congress
  Volume 3, Econometrics}.
\end{bincollection}
\endbibitem

\bibitem{belloni2018}
\begin{barticle}[author]
\bauthor{\bsnm{Belloni},~\bfnm{Alexandre}\binits{A.}},
  \bauthor{\bsnm{Chernozhukov},~\bfnm{Victor}\binits{V.}},
  \bauthor{\bsnm{Hansen},~\bfnm{Christian}\binits{C.}} \AND
  \bauthor{\bsnm{Newey},~\bfnm{Whitney}\binits{W.}}
(\byear{2018}).
\btitle{Simultaneous Confidence Intervals for High-dimensional Linear Models
  with Many Endogenous Variables}.
\bjournal{arXiv preprint arXiv:1712.08102}.
\end{barticle}
\endbibitem

\bibitem{BKRW98}
\begin{bbook}[author]
\bauthor{\bsnm{Bickel},~\bfnm{P.}\binits{P.}},
  \bauthor{\bsnm{Klaassen},~\bfnm{C.}\binits{C.}},
  \bauthor{\bsnm{Ritov},~\bfnm{Y.}\binits{Y.}} \AND
  \bauthor{\bsnm{Wellner},~\bfnm{J.}\binits{J.}}
(\byear{1998}).
\btitle{Efficient and adaptive inference in semiparametric models}.
\bpublisher{Spring-Verlag, New York}.
\end{bbook}
\endbibitem

\bibitem{BRT09}
\begin{barticle}[author]
\bauthor{\bsnm{Bickel},~\bfnm{P.}\binits{P.}},
  \bauthor{\bsnm{Ritov},~\bfnm{Y.}\binits{Y.}} \AND
  \bauthor{\bsnm{Tsybakov},~\bfnm{A.}\binits{A.}}
(\byear{2009}).
\btitle{Simultaneous analysis of Lasso and Dantzig selector}.
\bjournal{Ann. Statist.}
\bvolume{37}
\bpages{1705--1732}.
\end{barticle}
\endbibitem

\bibitem{BvdG11}
\begin{bbook}[author]
\bauthor{\bsnm{B{\"u}hlmann},~\bfnm{P.}\binits{P.}} \AND
  \bauthor{\bparticle{van~de} \bsnm{Geer},~\bfnm{S.}\binits{S.}}
(\byear{2011}).
\btitle{Statistics for High-Dimensional Data}.
\bpublisher{Springer}.
\end{bbook}
\endbibitem

\bibitem{CLL11}
\begin{barticle}[author]
\bauthor{\bsnm{Cai},~\bfnm{T.}\binits{T.}},
  \bauthor{\bsnm{Liu},~\bfnm{W.}\binits{W.}} \AND
  \bauthor{\bsnm{Luo},~\bfnm{X.}\binits{X.}}
(\byear{2011}).
\btitle{A constrained $\ell_1$ minimization approach to aparse precision matrix
  estimation}.
\bjournal{J. Amer. Statist. Assoc.}
\bvolume{106}
\bpages{594--607}.
\end{barticle}
\endbibitem

\bibitem{Caner&Fan2015hybrid}
\begin{barticle}[author]
\bauthor{\bsnm{Caner},~\bfnm{Mehmet}\binits{M.}} \AND
  \bauthor{\bsnm{Fan},~\bfnm{Qingliang}\binits{Q.}}
(\byear{2015}).
\btitle{Hybrid generalized empirical likelihood estimators: Instrument
  selection with adaptive lasso}.
\bjournal{Journal of Econometrics}
\bvolume{187}
\bpages{256--274}.
\end{barticle}
\endbibitem

\bibitem{Caner&Kock2018}
\begin{barticle}[author]
\bauthor{\bsnm{Caner},~\bfnm{Mehmet}\binits{M.}} \AND
  \bauthor{\bsnm{Kock},~\bfnm{Anders~Bredahl}\binits{A.~B.}}
(\byear{2018}).
\btitle{High Dimensional Linear GMM}.
\bjournal{arXiv preprint arXiv:1811.08779}.
\end{barticle}
\endbibitem

\bibitem{CanerKock18}
\begin{barticle}[author]
\bauthor{\bsnm{Caner},~\bfnm{Mehmet}\binits{M.}} \AND
  \bauthor{\bsnm{Kock},~\bfnm{Anders~Bredahl}\binits{A.~B.}}
(\byear{2018}).
\btitle{Asymptotically honest confidence regions for high dimensional
  parameters by the desparsified conservative Lasso}.
\bjournal{Journal of Econometrics}
\bvolume{203}
\bpages{143 - 168}.
\bdoi{https://doi.org/10.1016/j.jeconom.2017.11.005}
\end{barticle}
\endbibitem

\bibitem{chamberlain1987asymptotic}
\begin{barticle}[author]
\bauthor{\bsnm{Chamberlain},~\bfnm{Gary}\binits{G.}}
(\byear{1987}).
\btitle{Asymptotic efficiency in estimation with conditional moment
  restrictions}.
\bjournal{Journal of Econometrics}
\bvolume{34}
\bpages{305--334}.
\end{barticle}
\endbibitem

\bibitem{Cheng&Liao2015}
\begin{barticle}[author]
\bauthor{\bsnm{Cheng},~\bfnm{Xu}\binits{X.}} \AND
  \bauthor{\bsnm{Liao},~\bfnm{Zhipeng}\binits{Z.}}
(\byear{2015}).
\btitle{Select the valid and relevant moments: An information-based lasso for
  gmm with many moments}.
\bjournal{Journal of Econometrics}
\bvolume{186}
\bpages{443--464}.
\end{barticle}
\endbibitem

\bibitem{Chetelat17}
\begin{barticle}[author]
\bauthor{\bsnm{Ch\'etelat},~\bfnm{Didier}\binits{D.}},
  \bauthor{\bsnm{Lederer},~\bfnm{Johannes}\binits{J.}} \AND
  \bauthor{\bsnm{Salmon},~\bfnm{Joseph}\binits{J.}}
(\byear{2017}).
\btitle{Optimal two-step prediction in regression}.
\bjournal{Electron. J. Statist.}
\bvolume{11}
\bpages{2519--2546}.
\bdoi{10.1214/17-EJS1287}
\end{barticle}
\endbibitem

\bibitem{CLW16}
\begin{barticle}[author]
\bauthor{\bsnm{Chichignoud},~\bfnm{M.}\binits{M.}},
  \bauthor{\bsnm{Lederer},~\bfnm{J.}\binits{J.}} \AND
  \bauthor{\bsnm{Wainwright},~\bfnm{M.}\binits{M.}}
(\byear{2016}).
\btitle{A practical scheme and fast algorithm to tune the Lasso with optimality
  guarantees}.
\bjournal{J. Mach. Learn. Res.}
\bvolume{17}
\bpages{1--17}.
\end{barticle}
\endbibitem

\bibitem{Dalalyan17}
\begin{barticle}[author]
\bauthor{\bsnm{Dalalyan},~\bfnm{Arnak}\binits{A.}},
  \bauthor{\bsnm{Hebiri},~\bfnm{Mohamed}\binits{M.}} \AND
  \bauthor{\bsnm{Lederer},~\bfnm{Johannes}\binits{J.}}
(\byear{2017}).
\btitle{On the prediction performance of the Lasso}.
\bjournal{Bernoulli}
\bvolume{23}
\bpages{552--581}.
\bdoi{10.3150/15-BEJ756}
\end{barticle}
\endbibitem

\bibitem{Fan&Li2001}
\begin{barticle}[author]
\bauthor{\bsnm{Fan},~\bfnm{Jianqing}\binits{J.}} \AND
  \bauthor{\bsnm{Li},~\bfnm{Runze}\binits{R.}}
(\byear{2001}).
\btitle{Variable selection via nonconcave penalized likelihood and its oracle
  properties}.
\bjournal{Journal of the American statistical Association}
\bvolume{96}
\bpages{1348--1360}.
\end{barticle}
\endbibitem

\bibitem{FL14}
\begin{barticle}[author]
\bauthor{\bsnm{Fan},~\bfnm{J.}\binits{J.}} \AND
  \bauthor{\bsnm{Liao},~\bfnm{Y.}\binits{Y.}}
(\byear{2014}).
\btitle{Endogeneity in high dimensions}.
\bjournal{Ann. Statist.}
\bvolume{42}
\bpages{872--917}.
\end{barticle}
\endbibitem

\bibitem{Fan&Zhong2018}
\begin{barticle}[author]
\bauthor{\bsnm{Fan},~\bfnm{Qingliang}\binits{Q.}} \AND
  \bauthor{\bsnm{Zhong},~\bfnm{Wei}\binits{W.}}
(\byear{2018}).
\btitle{Nonparametric additive instrumental variable estimator: A group
  shrinkage estimation perspective}.
\bjournal{Journal of Business \& Economic Statistics}
\bvolume{36}
\bpages{388--399}.
\end{barticle}
\endbibitem

\bibitem{glmnet}
\begin{barticle}[author]
\bauthor{\bsnm{Friedman},~\bfnm{Jerome}\binits{J.}},
  \bauthor{\bsnm{Hastie},~\bfnm{Trevor}\binits{T.}} \AND
  \bauthor{\bsnm{Tibshirani},~\bfnm{Robert}\binits{R.}}
(\byear{2010}).
\btitle{Regularization paths for generalized linear models via coordinate
  descent}.
\bjournal{J. Stat. Softw.}
\bvolume{33}
\bpages{1--22}.
\end{barticle}
\endbibitem

\bibitem{GT11}
\begin{barticle}[author]
\bauthor{\bsnm{{Gautier}},~\bfnm{E.}\binits{E.}} \AND
  \bauthor{\bsnm{{Tsybakov}},~\bfnm{A.}\binits{A.}}
(\byear{2011}).
\btitle{{High-dimensional instrumental variables regression and confidence
  sets}}.
\bjournal{ArXiv e-prints}.
\end{barticle}
\endbibitem

\bibitem{GT14}
\begin{barticle}[author]
\bauthor{\bsnm{Gautier},~\bfnm{E.}\binits{E.}} \AND
  \bauthor{\bsnm{Tsybakov},~\bfnm{A.}\binits{A.}}
(\byear{2014}).
\btitle{High-dimensional instrumental variables regression and confidence
  sets}.
\bjournal{ArXiv:1105.2454v4}.
\end{barticle}
\endbibitem

\bibitem{GT18}
\begin{barticle}[author]
\bauthor{\bsnm{{Gautier}},~\bfnm{E.}\binits{E.}} \AND
  \bauthor{\bsnm{{Tsybakov}},~\bfnm{A.}\binits{A.}}
(\byear{2018}).
\btitle{{High-dimensional instrumental variables regression and confidence
  sets}}.
\bjournal{ArXiv e-prints}.
\end{barticle}
\endbibitem

\bibitem{G14}
\begin{bbook}[author]
\bauthor{\bsnm{Giraud},~\bfnm{C.}\binits{C.}}
(\byear{2014}).
\btitle{Introduction to high-dimensional statistics}.
\bseries{Monographs on statistics and applied probability (Series); 139}.
\bpublisher{CRC Press, Taylor \& Francis Group.}
\end{bbook}
\endbibitem

\bibitem{Hansen&Kozbur2014RJIVE}
\begin{barticle}[author]
\bauthor{\bsnm{Hansen},~\bfnm{Christian}\binits{C.}} \AND
  \bauthor{\bsnm{Kozbur},~\bfnm{Damian}\binits{D.}}
(\byear{2014}).
\btitle{Instrumental variables estimation with many weak instruments using
  regularized JIVE}.
\bjournal{Journal of Econometrics}
\bvolume{182}
\bpages{290--308}.
\end{barticle}
\endbibitem

\bibitem{H82}
\begin{barticle}[author]
\bauthor{\bsnm{Hansen},~\bfnm{L.}\binits{L.}}
(\byear{1982}).
\btitle{Large sample properties of generalized method of moments estimators}.
\bjournal{Econometrica}
\bvolume{50}
\bpages{1029--1054}.
\end{barticle}
\endbibitem

\bibitem{HTW15}
\begin{bbook}[author]
\bauthor{\bsnm{Hastie},~\bfnm{T.}\binits{T.}},
  \bauthor{\bsnm{Tibshirani},~\bfnm{R.}\binits{R.}} \AND
  \bauthor{\bsnm{Wainwright},~\bfnm{M.}\binits{M.}}
(\byear{2015}).
\btitle{Statistical learning with sparsity: The Lasso and generalizations}.
\bseries{Monographs on statistics and applied probability (Series); 143}.
\bpublisher{Boca Raton: CRC Press, Taylor \& Francis Group.}
\end{bbook}
\endbibitem

\bibitem{hebiri2013correlations}
\begin{barticle}[author]
\bauthor{\bsnm{Hebiri},~\bfnm{Mohamed}\binits{M.}} \AND
  \bauthor{\bsnm{Lederer},~\bfnm{Johannes}\binits{J.}}
(\byear{2013}).
\btitle{How correlations influence lasso prediction}.
\bjournal{IEEE Transactions on Information Theory}
\bvolume{59}
\bpages{1846--1854}.
\end{barticle}
\endbibitem

\bibitem{IDN03}
\begin{barticle}[author]
\bauthor{\bsnm{Imbens},~\bfnm{G.}\binits{G.}},
  \bauthor{\bsnm{Donald},~\bfnm{S.}\binits{S.}} \AND
  \bauthor{\bsnm{Newey},~\bfnm{W.}\binits{W.}}
(\byear{2003}).
\btitle{Empirical likelihood estimation and consistent tests with conditional
  moment restrictions}.
\bjournal{J. Econometrics}
\bvolume{117}
\bpages{55--93}.
\end{barticle}
\endbibitem

\bibitem{Jankova16}
\begin{barticle}[author]
\bauthor{\bsnm{Jankova},~\bfnm{Jana}\binits{J.}} \AND
  \bauthor{\bparticle{van~de} \bsnm{Geer},~\bfnm{Sara}\binits{S.}}
(\byear{2016}).
\btitle{Semi-parametric efficiency bounds and efficient estimation for
  high-dimensional models}.
\bjournal{To appear in Ann. Statist.}
\end{barticle}
\endbibitem

\bibitem{JM14}
\begin{barticle}[author]
\bauthor{\bsnm{Javanmard},~\bfnm{A.}\binits{A.}} \AND
  \bauthor{\bsnm{Montanari},~\bfnm{A.}\binits{A.}}
(\byear{2014}).
\btitle{Confidence intervals and hypothesis testing for high-dimensional
  regression}.
\bjournal{J. Mach. Learn. Res.}
\bvolume{15}
\bpages{2869--2909}.
\end{barticle}
\endbibitem

\bibitem{Knight00}
\begin{barticle}[author]
\bauthor{\bsnm{Knight},~\bfnm{Keith}\binits{K.}} \AND
  \bauthor{\bsnm{Fu},~\bfnm{Wenjiang}\binits{W.}}
(\byear{2000}).
\btitle{Asymptotics for lasso-type estimators}.
\bjournal{Ann. Statist.}
\bvolume{28}
\bpages{1356--1378}.
\end{barticle}
\endbibitem

\bibitem{lederer2019}
\begin{barticle}[author]
\bauthor{\bsnm{Lederer},~\bfnm{Johannes}\binits{J.}},
  \bauthor{\bsnm{Yu},~\bfnm{Lu}\binits{L.}} \AND
  \bauthor{\bsnm{Gaynanova},~\bfnm{Irina}\binits{I.}}
(\byear{2019}).
\btitle{Oracle inequalities for high-dimensional prediction}.
\bjournal{Bernoulli}
\bvolume{25}
\bpages{1225--1255}.
\bdoi{10.3150/18-BEJ1019}
\end{barticle}
\endbibitem

\bibitem{lin2015regularization}
\begin{barticle}[author]
\bauthor{\bsnm{Lin},~\bfnm{Wei}\binits{W.}},
  \bauthor{\bsnm{Feng},~\bfnm{Rui}\binits{R.}} \AND
  \bauthor{\bsnm{Li},~\bfnm{Hongzhe}\binits{H.}}
(\byear{2015}).
\btitle{Regularization methods for high-dimensional instrumental variables
  regression with an application to genetical genomics}.
\bjournal{Journal of the American Statistical Association}
\bvolume{110}
\bpages{270--288}.
\end{barticle}
\endbibitem

\bibitem{LWZ15}
\begin{barticle}[author]
\bauthor{\bsnm{Liu},~\bfnm{H.}\binits{H.}},
  \bauthor{\bsnm{Wang},~\bfnm{L.}\binits{L.}} \AND
  \bauthor{\bsnm{Zhao},~\bfnm{T.}\binits{T.}}
(\byear{2015}).
\btitle{Calibrated multivariate regression with application to neural semantic
  basis discovery}.
\bjournal{J. Mach. Learn. Res.}
\bvolume{16}
\bpages{1579--1606}.
\end{barticle}
\endbibitem

\bibitem{N90a}
\begin{barticle}[author]
\bauthor{\bsnm{Newey},~\bfnm{W.}\binits{W.}}
(\byear{1990}).
\btitle{Efficient instrumental variables estimation of nonlinear models}.
\bjournal{Econometrica}
\bvolume{58}
\bpages{809--837}.
\end{barticle}
\endbibitem

\bibitem{newey1990efficient}
\begin{barticle}[author]
\bauthor{\bsnm{Newey},~\bfnm{Whitney~K}\binits{W.~K.}}
(\byear{1990}).
\btitle{Efficient instrumental variables estimation of nonlinear models}.
\bjournal{Econometrica: Journal of the Econometric Society}
\bpages{809--837}.
\end{barticle}
\endbibitem

\bibitem{neykov15}
\begin{barticle}[author]
\bauthor{\bsnm{Neykov},~\bfnm{M.}\binits{M.}},
  \bauthor{\bsnm{Ning},~\bfnm{Y.}\binits{Y.}},
  \bauthor{\bsnm{Liu},~\bfnm{J.}\binits{J.}} \AND
  \bauthor{\bsnm{Liu.},~\bfnm{H.}\binits{H.}}
(\byear{2015}).
\btitle{A unified theory of confidence regions and testing for high dimensional
  estimating equations}.
\bjournal{arXiv:1510.08986}.
\end{barticle}
\endbibitem

\bibitem{Ning17}
\begin{barticle}[author]
\bauthor{\bsnm{Ning},~\bfnm{Yang}\binits{Y.}} \AND
  \bauthor{\bsnm{Liu},~\bfnm{Han}\binits{H.}}
(\byear{2017}).
\btitle{A general theory of hypothesis tests and confidence regions for sparse
  high dimensional models}.
\bjournal{Ann. Statist.}
\bvolume{45}
\bpages{158--195}.
\bdoi{10.1214/16-AOS1448}
\end{barticle}
\endbibitem

\bibitem{Potscher09}
\begin{barticle}[author]
\bauthor{\bsnm{P\"{o}tscher},~\bfnm{Benedikt}\binits{B.}} \AND
  \bauthor{\bsnm{Leeb},~\bfnm{Hannes}\binits{H.}}
(\byear{2009}).
\btitle{On the distribution of penalized maximum likelihood estimators: The
  {LASSO}, {SCAD}, and thresholding}.
\bjournal{J. Multivar. Anal.}
\bvolume{100}
\bpages{2065--2082}.
\end{barticle}
\endbibitem

\bibitem{RZ12}
\begin{barticle}[author]
\bauthor{\bsnm{Rudelson},~\bfnm{M.}\binits{M.}} \AND
  \bauthor{\bsnm{Zhou},~\bfnm{S.}\binits{S.}}
(\byear{2012}).
\btitle{Reconstruction from anisotropic random measurements}.
\bjournal{JMLR Workshop Conf. Proc.}
\bvolume{23}
\bpages{10.1--10.28}.
\end{barticle}
\endbibitem

\bibitem{stock2003}
\begin{barticle}[author]
\bauthor{\bsnm{Stock},~\bfnm{J.}\binits{J.}} \AND
  \bauthor{\bsnm{Trebbi},~\bfnm{F.}\binits{F.}}
(\byear{2003}).
\btitle{Retrospectives: Who invented instrumental variable regression?}
\bjournal{J. Econ. Perspect.}
\bvolume{17}
\bpages{177--194}.
\end{barticle}
\endbibitem

\bibitem{vdGB09}
\begin{barticle}[author]
\bauthor{\bparticle{van~de} \bsnm{Geer},~\bfnm{S.}\binits{S.}} \AND
  \bauthor{\bsnm{B{\"u}hlmann},~\bfnm{P.}\binits{P.}}
(\byear{2009}).
\btitle{On the conditions used to prove oracle results for the Lasso}.
\bjournal{Electron. J. Stat.}
\bvolume{3}
\bpages{1360--1392}.
\end{barticle}
\endbibitem

\bibitem{vdGetal14}
\begin{barticle}[author]
\bauthor{\bparticle{van~de} \bsnm{Geer},~\bfnm{S.}\binits{S.}},
  \bauthor{\bsnm{B{\"u}hlmann},~\bfnm{P.}\binits{P.}},
  \bauthor{\bsnm{Ritov},~\bfnm{Y.}\binits{Y.}} \AND
  \bauthor{\bsnm{Dezeure},~\bfnm{R.}\binits{R.}}
(\byear{2014}).
\btitle{On asymptotically optimal confidence regions and tests for
  high-dimensional models}.
\bjournal{Ann. Statist.}
\bvolume{42}
\bpages{1166--1202}.
\end{barticle}
\endbibitem

\bibitem{vandegeer}
\begin{binbook}[author]
\bauthor{\bparticle{van~de} \bsnm{Geer},~\bfnm{Sara}\binits{S.}} \AND
  \bauthor{\bsnm{Lederer},~\bfnm{Johannes}\binits{J.}}
(\byear{2013}).
\btitle{The Lasso, correlated design, and improved oracle inequalities}.
In \bbooktitle{From Probability to Statistics and Back: High-Dimensional Models
  and Processes -- A Festschrift in Honor of Jon A. Wellner}.
\bseries{Collections}
\bvolume{Volume 9}
\bpages{303--316}.
\bpublisher{Institute of Mathematical Statistics}, \baddress{Beachwood, Ohio,
  USA}.
\bdoi{10.1214/12-IMSCOLL922}
\end{binbook}
\endbibitem

\bibitem{vdGM14}
\begin{barticle}[author]
\bauthor{\bparticle{van~de} \bsnm{Geer},~\bfnm{S.}\binits{S.}} \AND
  \bauthor{\bsnm{Muro},~\bfnm{A.}\binits{A.}}
(\byear{2014}).
\btitle{On higher order isotropy conditions and lower bounds for sparse
  quadratic forms}.
\bjournal{Electron. J. Stat.}
\bvolume{8}
\bpages{3031-3061}.
\end{barticle}
\endbibitem

\bibitem{V96}
\begin{bbook}[author]
\bauthor{\bparticle{van~der} \bsnm{Vaart},~\bfnm{A.}\binits{A.}} \AND
  \bauthor{\bsnm{Wellner},~\bfnm{J.}\binits{J.}}
(\byear{1996}).
\btitle{Weak convergence and empirical processes: With applications to
  statistics}.
\bpublisher{Springer-Verlag New York}.
\end{bbook}
\endbibitem

\bibitem{V12}
\begin{bincollection}[author]
\bauthor{\bsnm{Vershynin},~\bfnm{R.}\binits{R.}}
(\byear{2012}).
\btitle{Introduction to the non-asymptotic analysis of random matrices}.
In \bbooktitle{Compressed Sensing: Theory and Applications}
(\beditor{\bfnm{Yonina~C.}\binits{Y.~C.}~\bsnm{Eldar}} \AND
  \beditor{\bfnm{Gitta}\binits{G.}~\bsnm{Kutyniok}}, eds.)
\bchapter{5},
\bpages{210--268}.
\bpublisher{Cambridge University Press}.
\end{bincollection}
\endbibitem

\bibitem{zhang2010nearly}
\begin{barticle}[author]
\bauthor{\bsnm{Zhang},~\bfnm{Cun-Hui}\binits{C.-H.}} \betal{et~al.}
(\byear{2010}).
\btitle{Nearly unbiased variable selection under minimax concave penalty}.
\bjournal{The Annals of statistics}
\bvolume{38}
\bpages{894--942}.
\end{barticle}
\endbibitem

\bibitem{ZZ14}
\begin{barticle}[author]
\bauthor{\bsnm{Zhang},~\bfnm{C.~H.}\binits{C.~H.}} \AND
  \bauthor{\bsnm{Zhang},~\bfnm{S.}\binits{S.}}
(\byear{2014}).
\btitle{Confidence intervals for low dimensional parameters in high dimensional
  linear models}.
\bjournal{J. R. Stat. Soc. Ser. B. Stat. Methodol.}
\bvolume{76}
\bpages{217--242}.
\bdoi{10.1111/rssb.12026}
\end{barticle}
\endbibitem

\bibitem{Zhu15}
\begin{barticle}[author]
\bauthor{\bsnm{Zhu},~\bfnm{Ying}\binits{Y.}}
(\byear{2018}).
\btitle{Sparse linear models and $\ell_1$-regularized {2SLS} with
  high-dimensional endogenous regressors and instruments}.
\bjournal{To appear in J. Econometrics}.
\end{barticle}
\endbibitem

\bibitem{doi:10.1002/sta4.186}
\begin{barticle}[author]
\bauthor{\bsnm{Zhuang},~\bfnm{Rui}\binits{R.}} \AND
  \bauthor{\bsnm{Lederer},~\bfnm{Johannes}\binits{J.}}
(\byear{2018}).
\btitle{Maximum regularized likelihood estimators: A general prediction theory
  and applications}.
\bjournal{Stat}
\bvolume{7}
\bpages{e186}.
\bnote{e186 sta4.186}.
\bdoi{10.1002/sta4.186}
\end{barticle}
\endbibitem

\end{thebibliography}
